\documentclass[10pt]{article}
\usepackage{setspace}
  
\usepackage{amsmath} 
\usepackage{amssymb}
\usepackage{latexsym}
\usepackage{xcolor}
\usepackage{fancyhdr}
\allowdisplaybreaks 

 \usepackage{comment}

\def\XXint#1#2#3{{\setbox0=\hbox{$#1{#2#3}{\int}$} 
\vcenter{\hbox{$#2#3$}}\kern-.5\wd0}}   

 \numberwithin{equation}{section}
\newtheorem{theorem}[equation]{Theorem}
\newtheorem{proposition}[equation]{Proposition}
\newtheorem{definition}[equation]{Definition}
\newtheorem{remark}[equation]{Remark}
\newtheorem{lemma}[equation]{Lemma}
\newtheorem{corollary}[equation]{Corollary}

\title{
Representation theorems for  nonvariational solutions of  the Helmholtz equation
 } 
 
\author{  
Massimo Lanza de Cristoforis
\\
Dipartimento di Matematica `Tullio Levi-Civita', 
\\
Universit\`a degli Studi di Padova, 
\\
Via Trieste 63, Padova 35121, 
Italy. 
\\
E-mail: mldc@math.unipd.it   }

\date{\ }

\begin{document}
 
 \maketitle

\noindent
{\bf Abstract:}  We consider a possibly multiply connected bounded open subset $\Omega$ of ${\mathbb{R}}^n$ of class $C^{\max\{1,m\},\alpha}$ for some $m\in {\mathbb{N}}$, 
$\alpha\in]0,1[$ and we plan to solve both the Dirichlet and the Neumann problem  for the Helmholtz equation in   $\Omega$ and  in the exterior  of $\Omega$ in terms of acoustic  layer potentials.   Then we turn to 
prove an integral representation theorem    solutions of the Helmholtz equation in terms of an  acoustic  single layer   potential. 
The main focus of the paper is on $\alpha$-H\"{o}lder continuous solutions which  may  not have a classical normal derivative at the boundary points of $\Omega$ and  that may have an infinite Dirichlet integral around the boundary of $\Omega$\, \textit{i.e.}, case $m=0$. Namely for solutions  that  do not belong to the classical variational setting.

 \vspace{\baselineskip}



\section{Introduction}  

 Our starting point are the classical examples of Prym \cite{Pr1871} and  Hadamard \cite{Ha1906} of harmonic functions in the unit ball of the plane that are continuous up to the boundary and have  infinite Dirichlet integral, \textit{i.e.},  whose
 gradient is not square-summable. Such functions solve the classical Dirichlet problem in the unit ball, but not the corresponding weak (variational) formulation. For a discussion on this point we refer to 
Maz’ya and Shaposnikova \cite{MaSh98}, Bottazzini and Gray \cite{BoGr13} and
 Bramati,  Dalla Riva and  Luczak~\cite{BrDaLu23} which contains examples of H\"{o}lder continuous harmonic functions with infinite Dirichlet integral.

In the papers \cite{La24b}, \cite{La24c}, the author has analyzed the Neumann problem for the Laplace and the Poisson equation. In particular, the author has introduced a distributional normal derivative on the boundary  for $\alpha$-H\"{o}lder continuous functions that have a Laplace operator in a Schauder space of negative exponent and has characterized the corresponding space for the first order traces.

In \cite{La25}, \cite{La25a}, \cite{La25b}, the author has carried out the corresponding analysis for the Helmholtz equation. 

In the present paper, we plan to consider both the interior and the exterior Dirichlet and   Neumann problems for the Helmholtz equation in an open subset $\Omega$ of ${\mathbb{R}}^n$ of class $C^{\max\{1,m\},\alpha}$ for some $m\in {\mathbb{N}}$, $\alpha\in]0,1[$ and to prove corresponding existence theorems for solutions that can be written as sums of  acoustic  layer potentials. Here we note that $\Omega$ may be multiply connected.\par 

Case $m=0$ corresponds the case in which solutions may  not have a classical normal derivative at the boundary points of $\Omega$ and that may have an infinite Dirichlet integral around the boundary of $\Omega$. Namely for solutions  that  do not belong to the classical variational setting and that we address to as `nonvariational'.   Case $m\geq 1$ instead corresponds to the classical case.

In the classical case, $m\geq 1$ existence theorems for the interior and exterior Dirichlet and Neumann problems for the Helmholtz equation are well known in case the exterior of $\Omega$ is connected. Here we mention major monographs as that of Colton and Kress \cite[p.~32 line 4, Thms.~3.20, 3.21, 3.24, 3.25]{CoKr92}, Kirsch and  Hettlich \cite{KiHe15}. In the present paper we present a solution in terms of  acoustic  layer potentials also in case 
the exterior of $\Omega$ is not connected. Moreover, we prove  corresponding results also in the `nonvariational' case $m=0$, which is the case of major focus of the present paper:  see Theorem \ref{mldc_thm:exintdirhe} for the interior Dirichlet problem, 
Theorem \ref{mldc_thm:exextdirhe} for the exterior Dirichlet problem, Theorem \ref{mldc_thm:exintneuhe} for the interior Neumann problem and Theorem \ref{mldc_thm:exextneuhe} for the exterior Neumann  problem.

Then we turn to the problem of representing solutions of the Helmholtz equation  in terms of  acoustic  single layer potentials. The problem of  representing
solutions of second order elliptic boundary value problems as single layer potentials  is by now classic and we wish to mention the contributions for  high order  elliptic equations and systems in two variables of  Fichera \cite{Fi60},  Ricci \cite{Ri74},  Fichera and Ricci \cite{FiRi76}, the corresponding applications to numerical analysis of
  MacCamy \cite{Mac66}, 
 Hsiao and  MacCamy  \cite{HsMac73},  
 Hsiao and Wendland \cite{HsWe77}  and the contribution of  Cialdea \cite{Ci91} in $n$ variables. In the present paper,  the emphasis is not on the generality of the elliptic operator or the number of variables. Instead, we look for representations for solution of the Helmholtz equation  which   may  not have a classical normal derivative at the boundary points of $\Omega$ and that may have an infinite Dirichlet integral around the boundary of $\Omega$  in terms of   acoustic  single layer potentials with a   density  which may be a distribution on the boundary, and with a one-to-one correspondence between the density and the  acoustic  single layer, possibly adding to the  acoustic  single layer a term
 in a finite dimensional space. The analysis of the present paper and of 
 \cite{La24b}, \cite{La24c},  \cite{La25}, \cite{La25a}, \cite{La25b}  is motivated by applications
 to  singularly perturbed problems as in the monograph of Dalla Riva, the author and Musolino \cite{DaLaMu21}. More specifically to problems of the electromagnetism. Thus we prove the representation Theorem \ref{mldc_thm:reprhelfun} for solutions of the Helmholtz equation  in terms of a single layer (which may have a distributional density), up to an additive term in the form of a double layer with a density in a finite dimensional space. Then we also note the validity of Corollary \ref{mldc_corol:reprhelfun}, that generalizes to case $m=0$ the variational case $m\geq1$ that has been exploited to solve a singularly perturbed problem for the Helmholtz equation by Akyel and the author in 
   \cite[Thm.~A1 of the Appendix]{AkLa22}. At the end of the paper we have included an Appendix
  on a consequence of a Holmgren uniqueness theorem that has been pointed out by  Lupo and Micheletti \cite[Proof of Theorem 2]{LuMi93}.

 \section{Preliminaries and notation}\label{mldc_sec:prelnot} Unless otherwise specified,  we assume  throughout the paper that
\[
n\in {\mathbb{N}}\setminus\{0,1\}\,,
\]
where ${\mathbb{N}}$ denotes the set of natural numbers including $0$. 
If $X$, $Y$ and $Z$ are normed spaces, then ${\mathcal{L}}(X,Y)$ denotes the space of linear and continuous maps from $X$ to $Y$ and ${\mathcal{L}}^{(2)}(X\times Y, Z)$ 
denotes the space of bilinear and continuous maps from $X\times Y$ to $Z$ with their usual operator norm (cf.~\textit{e.g.}, \cite[pp.~16, 621]{DaLaMu21}). $I$ denotes the identity map. ${\mathrm{Im}}$  and ${\mathrm{Ker}}$ denote  the image (or range) and the kernel (or null space) of an operator, respectively.
  $|A|$ denotes the operator norm of a matrix $A$ with real (or complex) entries, 
       $A^{t}$ denotes the transpose matrix of $A$.	 
 
 Let $\Omega$ be an open subset of ${\mathbb{R}}^n$. $C^{1}(\Omega)$ denotes the set of continuously differentiable functions from $\Omega$ to ${\mathbb{C}}$. 
 If $s\in {\mathbb{N}}\setminus\{0\}$, $f\in \left(C^{1}(\Omega)\right)^{s} $, then   $Df$ denotes the Jacobian matrix of $f$.   
 
 For the (classical) definition of   open Lipschitz subset of ${\mathbb{R}}^n$ and of   open subset of ${\mathbb{R}}^n$ 
   of class $C^{m}$ or of class $C^{m,\alpha}$
  and of the H\"{o}lder and Schauder spaces $C^{m,\alpha}(\overline{\Omega})$
  on the closure $\overline{\Omega}$ of  an open set $\Omega$ and 
  of the H\"{o}lder and Schauder spaces
   $C^{m,\alpha}(\partial\Omega)$ 
on the boundary $\partial\Omega$ of an open set $\Omega$ for some $m\in{\mathbb{N}}$, $\alpha\in ]0,1]$, we refer for example to
    Dalla Riva, the author and Musolino  \cite[\S 2.3, \S 2.6, \S 2.7, \S 2.9, \S 2.11, \S 2.13,   \S 2.20]{DaLaMu21}.  If $m\in {\mathbb{N}}$, 
 $C^{m}_b(\overline{\Omega})$ denotes the space of $m$-times continuously differentiable functions from $\Omega$ to ${\mathbb{C}}$ such that all 
the partial derivatives up to order $m$ have a bounded continuous extension to    $\overline{\Omega}$ and we set
\[
\|f\|_{   C^{m}_{b}(
\overline{\Omega} )   }\equiv
\sum_{|\eta|\leq m}\, \sup_{x\in \overline{\Omega}}|D^{\eta}f(x)|
\qquad\forall f\in C^{m}_{b}(
\overline{\Omega} )\,.
\]
If $\alpha\in ]0,1]$, then 
$C^{m,\alpha}_b(\overline{\Omega})$ denotes the space of functions of $C^{m}_{b}(
\overline{\Omega}) $  such that the  partial derivatives of order $m$ are $\alpha$-H\"{o}lder continuous in $\Omega$. Then we equip $C^{m,\alpha}_{b}(\overline{\Omega})$ with the norm
\[
\|f\|_{  C^{m,\alpha}_{b}(\overline{\Omega})  }\equiv 
\|f\|_{  C^{m }_{b}(\overline{\Omega})  }
+\sum_{|\eta|=m}|D^{\eta}f|_{\alpha}\qquad\forall f\in C^{m,\alpha}_{b}(\overline{\Omega})\,,
\]
where $|D^{\eta}f|_{\alpha}$ denotes the $\alpha$-H\"{o}lder constant of the partial derivative $D^{\eta}f$ of order $\eta$ of $f$ in $\Omega$. If $\Omega$ is bounded, we obviously have $C^{m }_{b}(\overline{\Omega})=C^{m } (\overline{\Omega})$ and $C^{m,\alpha}_{b}(\overline{\Omega})=C^{m,\alpha} (\overline{\Omega})$. 
  Then $C^{m,\alpha}_{{\mathrm{loc}}}(\overline{\Omega }) $  denotes 
the space  of those functions $f\in C^{m}(\overline{\Omega} ) $ such that $f_{|\overline{\Omega'}} $ belongs to $
C^{m,\alpha}(   \overline{ \Omega' }   )$ for all bounded open subsets $\Omega'$ of ${\mathbb{R}}^n$ such that $\overline{\Omega'}\subseteq\overline{\Omega}$. 
The space of complex valued functions of class $C^m$ with compact support in an open set $\Omega$ of ${\mathbb{R}}^n$ is denoted $C^m_c(\Omega)$ and similarly for $C^\infty_c(\Omega)$. We also set ${\mathcal{D}}(\Omega)\equiv C^\infty_c(\Omega)$. Then the dual ${\mathcal{D}}'(\Omega)$ is known to be the space of distributions in $\Omega$. The support either of a function or of a distribution is denoted by the abbreviation `${\mathrm{supp}}$'.  We also set
\begin{equation}\label{mldc_eq:exto}
\Omega^-\equiv {\mathbb{R}}^n\setminus\overline{\Omega}\,,
\end{equation} 
for the exterior of $\Omega$. If $\Omega$ is a bounded open Lipschitz subset of ${\mathbb{R}}^n$, then $\Omega$  is known to have has at most a finite number of pairwise disjoint  connected components, which we denote by $\Omega_{1}$,\dots, $\Omega_{\varkappa^{+}}$ and which are open and $\Omega^-$ is known to have  at most a finite number of  pairwise disjoint  connected components, which we denote by $(\Omega^{-})_{0}$,
$(\Omega^{-})_{1}$, \dots, $(\Omega^{-})_{\varkappa^{-}}$ and which are open.   One and only one of such connected components is unbounded. We denote it by $(\Omega^{-})_{0}$  (cf.~\textit{e.g.}, \cite[Lemma~2.38]{DaLaMu21}).\par  

We denote by $\nu_\Omega$ or simply by $\nu$ the outward unit normal of $\Omega$ on $\partial\Omega$. Then $\nu_{\Omega^-}=-\nu_\Omega$ is the outward unit normal of $\Omega^-$ on $\partial\Omega=\partial\Omega^-$.

Morever, we retain the standard notation for the Lebesgue spaces $L^p$ for $p\in [1,+\infty]$ (cf.~\textit{e.g.}, Folland \cite[Chapter~6]{Fo99}, \cite[\S 2.1]{DaLaMu21}) and
$m_n$ denotes the
$n$ dimensional Lebesgue measure.\par

We now summarize the definition and some elementary properties of the Schau\-der space $C^{-1,\alpha}(\overline{\Omega})$ 
by following the presentation of Dalla Riva, the author and Musolino \cite[\S 2.22]{DaLaMu21}.
\begin{definition} 
\label{mldc_defn:sch-1}\index{Schauder space!with negative exponent}
 Let $\alpha\in]0,1]$. Let $\Omega$ be a bounded open subset of ${\mathbb{R}}^{n}$. We denote by $C^{-1,\alpha}(\overline{\Omega})$ the subspace 
 \[
 \left\{
 f_{0}+\sum_{j=1}^{n}\frac{\partial}{\partial x_{j}}f_{j}:\,f_{j}\in 
 C^{0,\alpha}(\overline{\Omega})\ \forall j\in\{0,\dots,n\}
 \right\}\,,
 \]
 of the space of distributions ${\mathcal{D}}'(\Omega)$  in $\Omega$ and we set
 \begin{eqnarray}
\label{mldc_defn:sch-2}
\lefteqn{
\|f\|_{  C^{-1,\alpha}(\overline{\Omega})  }
\equiv\inf\biggl\{\biggr.
\sum_{j=0}^{n}\|f_{j}\|_{ C^{0,\alpha}(\overline{\Omega})  }
:\,
}
\\ \nonumber
&&\qquad\qquad\qquad\qquad
f=f_{0}+\sum_{j=1}^{n}\frac{\partial}{\partial x_{j}}f_{j}\,,\ 
f_{j}\in C^{0,\alpha}(\overline{\Omega})\ \forall j\in \{0,\dots,n\}
\biggl.\biggr\}\,.
\end{eqnarray}
\end{definition}
$(C^{-1,\alpha}(\overline{\Omega}), \|\cdot\|_{  C^{-1,\alpha}(\overline{\Omega})  })$ is known to be a Banach space and  is continuously embedded into ${\mathcal{D}}'(\Omega)$. Also, the definition of the norm $\|\cdot\|_{  C^{-1,\alpha}(\overline{\Omega})  }$ implies that $C^{0,\alpha}(\overline{\Omega})$ is continuously embedded into $C^{-1,\alpha}(\overline{\Omega})$ and that the partial differentiation $\frac{\partial}{\partial x_{j}}$ is continuous from 
$C^{0,\alpha}(\overline{\Omega})$ to $C^{-1,\alpha}(\overline{\Omega})$ for all $j\in\{1,\dots,n\}$.  Generically, the  elements of $C^{-1,\alpha}(\overline{\Omega})$ for $\alpha\in]0,1[$ are not integrable functions, but distributions in $\Omega$.  Then we have the following statement of \cite[Prop.~3.1]{La24c} that shows that the elements of $C^{-1,\alpha}(\overline{\Omega}) $, which belong to the dual of ${\mathcal{D}}(\Omega)$, can be extended to elements of the dual of $C^{1,\alpha}(\overline{\Omega})$.  
\begin{proposition}\label{mldc_prop:nschext}
 Let $\alpha\in]0,1[$. Let $\Omega$ be a bounded open Lipschitz subset of ${\mathbb{R}}^{n}$.  There exists one and only one  linear and continuous extension operator $E^\sharp_\Omega$ from $C^{-1,\alpha}(\overline{\Omega})$ to $\left(C^{1,\alpha}(\overline{\Omega})\right)'$ such that
 \begin{eqnarray}\label{mldc_prop:nschext2}
\lefteqn{
\langle E^\sharp_\Omega[f],v\rangle 
}
\\ \nonumber
&&\ \
 =
\int_{\Omega}f_{0}v\,dx+\int_{\partial\Omega}\sum_{j=1}^{n} (\nu_{\Omega})_{j}f_{j}v\,d\sigma
 -\sum_{j=1}^{n}\int_{\Omega}f_{j}\frac{\partial v}{\partial x_j}\,dx
\quad \forall v\in C^{1,\alpha}(\overline{\Omega})
\end{eqnarray}
for all $f=  f_{0}+\sum_{j=1}^{n}\frac{\partial}{\partial x_{j}}f_{j}\in C^{-1,\alpha}(\overline{\Omega}) $. Moreover, 
\begin{equation}\label{mldc_prop:nschext1}
E^\sharp_\Omega[f]_{|\Omega}=f\,, \ i.e.,\ 
\langle E^\sharp_\Omega[f],v\rangle =\langle f,v\rangle \qquad\forall v\in {\mathcal{D}}(\Omega)
\end{equation}
for all $f\in C^{-1,\alpha}(\overline{\Omega})$ and
\begin{equation}\label{mldc_prop:nschext3}
\langle E^\sharp_\Omega[f],v\rangle =\langle f,v\rangle =\int_\Omega fv\,dx
\qquad\forall v\in C^{1,\alpha}(\overline{\Omega})
\end{equation}
for all $f\in C^{0,\alpha}(\overline{\Omega})$.
\end{proposition}
When no ambiguity can arise, we simply write $E^\sharp$ instead of $E^\sharp_\Omega$. To see why the extension operator $E^\sharp$ can be considered as `canonical', we refer to \cite[Prop.~7]{La24d}. Next we introduce the following multiplication lemma. For a proof, we refer to \cite[Lem.~2.3]{La25}
\begin{lemma}\label{mldc_lem:multc1ac-1a}
  Let   $\alpha\in ]0,1[$. Let $\Omega$ be a bounded open Lipschitz subset of ${\mathbb{R}}^{n}$. Then the  product is bilinear and continuous from 
  $C^{1,\alpha}(\overline{\Omega})\times C^{-1,\alpha}(\overline{\Omega})$ to $C^{-1,\alpha}(\overline{\Omega})$. 
\end{lemma}
 Next we introduce our space for the solutions.

\begin{definition}\label{mldc_defn:c0ade}
 Let   $\alpha\in ]0,1]$. Let $\Omega$ be a bounded open  subset of ${\mathbb{R}}^{n}$. Let
 \begin{eqnarray}\label{mldc_defn:c0ade1}
C^{0,\alpha}(\overline{\Omega})_\Delta
&\equiv&\biggl\{u\in C^{0,\alpha}(\overline{\Omega}):\,\Delta u\in C^{-1,\alpha}(\overline{\Omega})\biggr\}\,,
\\ \nonumber
\|u\|_{ C^{0,\alpha}(\overline{\Omega})_\Delta }
&\equiv& \|u\|_{ C^{0,\alpha}(\overline{\Omega})}
+\|\Delta u\|_{C^{-1,\alpha}(\overline{\Omega})}
\qquad\forall u\in C^{0,\alpha}(\overline{\Omega})_\Delta\,.
\end{eqnarray}
\end{definition}
Since $C^{0,\alpha}(\overline{\Omega})$ and $C^{-1,\alpha}(\overline{\Omega})$ are Banach spaces,   $\left(\|u\|_{ C^{0,\alpha}(\overline{\Omega})_\Delta }, \|\cdot \|_{ C^{0,\alpha}(\overline{\Omega})_\Delta }\right)$ is a Banach space. For subsets $\Omega$ that are not necessarily bounded, we introduce the following   statement of 
\cite[Lem.~2.5]{La25}.
\begin{lemma}\label{mldc_lem:c1alcof}
 Let   $\alpha\in ]0,1]	$. Let $\Omega$ be an open  subset of ${\mathbb{R}}^{n}$. Then the space
\begin{eqnarray}\label{mldc_lem:c1alcof1}
\lefteqn{
 C^{0,\alpha}_{	{\mathrm{loc}}	}(\overline{\Omega})_\Delta\equiv\biggl\{
 f\in C^{0}(\overline{\Omega}):\, f_{|\overline{\Omega'}}\in C^{0,\alpha}(\overline{\Omega'})_\Delta\ \text{for\ all\  } 
 }
\\ \nonumber
&&\qquad\qquad\qquad\qquad 
 \text{bounded\ open\ subsets}\ \Omega'\ \text{of}\ {\mathbb{R}}^n\ \text{such\ that} \  \overline{\Omega'}\subseteq\overline{\Omega}  \biggr\}\,,
\end{eqnarray}
 with the family of seminorms
\begin{eqnarray}\label{mldc_lem:c1alcof2}
\lefteqn{
{\mathcal{P}}_{C^{0,\alpha}_{	{\mathrm{loc}}	}(\overline{\Omega})_\Delta}\equiv
\biggl\{
\|\cdot\|_{C^{0,\alpha}(\overline{\Omega'})_\Delta}:\,  
}
\\ \nonumber
&& \qquad\qquad\qquad\qquad 
 \Omega'\ \text{is\ a\ bounded\ open  subset\ of\ }  {\mathbb{R}}^n \ \text{such\ that} \  \overline{\Omega'}\subseteq\overline{\Omega}
\biggr\}
\end{eqnarray}
is a Fr\'{e}chet space.
\end{lemma}
 Next we   introduce a subspace of  $C^{1,\alpha}(\overline{\Omega})$ that we need in the sequel as in \cite[Defn.~2.9]{La25}.  
\begin{definition}\label{mldc_defn:c1ade}
 Let   $\alpha\in ]0,1]		$. Let $\Omega$ be a bounded open  subset of ${\mathbb{R}}^{n}$. Let
 \begin{eqnarray}\label{mldc_defn:c1ade1}
C^{1,\alpha}(\overline{\Omega})_\Delta
&\equiv&\biggl\{u\in C^{1,\alpha}(\overline{\Omega}):\,\Delta u\in C^{0,\alpha}(\overline{\Omega})\biggr\}\,,
\\ \nonumber
\|u\|_{ C^{1,\alpha}(\overline{\Omega})_\Delta }
&\equiv& \|u\|_{ C^{1,\alpha}(\overline{\Omega})}
+\|\Delta u\|_{C^{0,\alpha}(\overline{\Omega})}
\qquad\forall u\in C^{1,\alpha}(\overline{\Omega})_\Delta\,.
\end{eqnarray}
\end{definition}
 If $\Omega$ is not necessarily bounded, then $C^{1,\alpha}_{	{\mathrm{loc}}	}(\overline{\Omega})_\Delta$ denotes the set of $f\in C^{0}(\overline{\Omega})$ such that $ f_{|\overline{\Omega'}}\in C^{1,\alpha}(\overline{\Omega'})_\Delta $ for all 
 bounded  open  subsets $\Omega'$ of ${\mathbb{R}}^n$ such that $\overline{\Omega'}\subseteq\overline{\Omega}$.

Next we introduce the following classical result on the Green operator for the interior Dirichlet problem. For a proof, we refer for example to \cite[Thm.~4.8]{La24b}.
 \begin{theorem}\label{mldc_thm:idwp}
Let $m\in {\mathbb{N}}$, $\alpha\in ]0,1[$. Let $\Omega$ be a bounded open  subset of ${\mathbb{R}}^{n}$ of class $C^{\max\{m,1\},\alpha}$.
Then the map ${\mathcal{G}}_{\Omega,d,+}$  from $C^{m,\alpha}(\partial\Omega)$ to the closed subspace 
 \begin{equation}\label{mldc_thm:idwp1}
C^{m,\alpha}_h(\overline{\Omega}) \equiv \{
u\in C^{m,\alpha}(\overline{\Omega}), u\ \text{is\ harmonic\ in}\ \Omega\}
\end{equation}
of $ C^{m,\alpha}(\overline{\Omega})$ that takes $v$ to the only solution $v^\sharp$ of the Dirichlet problem
\begin{equation}\label{mldc_defn:cinspo3}
\left\{
\begin{array}{ll}
 \Delta v^\sharp=0 & \text{in}\ \Omega\,,
 \\
v^\sharp_{|\partial\Omega} =v& \text{on}\ \partial\Omega 
\end{array}
\right.
\end{equation}
is a linear homeomorphism.
\end{theorem}
Next we introduce the following two approximation lemmas.
\begin{lemma}\label{mldc_lem:apr0ad}
 Let $\alpha\in]0,1[$. 
 Let $\Omega$ be a bounded open subset of ${\mathbb{R}}^n$ of class $C^{1,\alpha}$.  If $u\in C^{0,\alpha}(\overline{\Omega})_\Delta$, then there exists a sequence $\{u_j\}_{j\in {\mathbb{N}}}$ in 
 $C^{1,\alpha}(\overline{\Omega})_\Delta$ such that
 \begin{equation}\label{mldc_lem:apr1ad1}
 \sup_{j\in {\mathbb{N}}}\|u_j\|_{
 C^{0,\alpha}(\overline{\Omega})_\Delta
 }<+\infty\,,\qquad
 \lim_{j\to\infty}u_j=u\quad\text{in}\ C^{0,\beta}(\overline{\Omega})_\Delta\quad \forall\beta\in]0,\alpha[\,.
 \end{equation}
\end{lemma}
For a proof, we refer to \cite[Lem.~5.14]{La24c}.
\begin{lemma}\label{mldc_lem:apr1a}
 Let $\alpha\in]0,1]$. 
 Let $\Omega$ be a bounded open subset of ${\mathbb{R}}^n$ of class $C^{1,\alpha}$.  If $g\in C^{1,\alpha}(\overline{\Omega})$, then there exists a sequence $\{g_j\}_{j\in {\mathbb{N}}}$ in 
 $C^{\infty}(\overline{\Omega})$ such that
 \begin{equation}\label{mldc_lem:apr1a1}
 \sup_{j\in {\mathbb{N}}}\|g_j\|_{
 C^{1,\alpha}(\overline{\Omega}) 
 }<+\infty\,,\qquad
 \lim_{j\to\infty}g_j=g\quad\text{in}\ C^{1,\beta}(\overline{\Omega})\quad \forall\beta\in]0,\alpha[\,.
 \end{equation}
\end{lemma}
For a proof, we refer to \cite[Lem.~A.3]{La24c}.  Next we introduce the  following Lemma that is  well known and is an immediate consequence of the H\"{o}lder inequality.
\begin{lemma}\label{mldc_lem:caincl}
 Let $m\in {\mathbb{N}}$, $\alpha\in ]0,1[$. Let $\Omega$ be a bounded open  subset of ${\mathbb{R}}^{n}$ of class $C^{\max\{m,1\},\alpha}$.  Then the canonical inclusion  ${\mathcal{J}}$ from the Lebesgue space $L^1(\partial\Omega)$ of integrable functions in $\partial\Omega$ to $(C^{m,\alpha}(\partial\Omega))'$ that takes $\mu$ to the functional ${\mathcal{J}}[\mu]$ defined by 
 \begin{equation}\label{mldc_lem:caincl1}
\langle {\mathcal{J}}[\mu],v\rangle \equiv \int_{\partial\Omega}\mu v\,d\sigma\qquad\forall v\in C^{
 m,
\alpha}(\partial\Omega)\,,
\end{equation}
is linear continuous and injective.
\end{lemma}
As customary, we say  that ${\mathcal{J}}[\mu]$ is the `distribution' that is canonically associated to $\mu$ and we omit the indication of the inclusion map ${\mathcal{J}}$. By Lemma \ref{mldc_lem:caincl}, the space $C^{0,\alpha}(\partial\Omega)$ is continuously embedded into $(C^{m,\alpha}(\partial\Omega))'$.  Next we plan to introduce the normal derivative of the functions in  $C^{0,\alpha}(\overline{\Omega})_\Delta$ as in \cite{La24c}. To do so,    we introduce the (classical) interior Steklov-Poincar\'{e} operator (or interior  Dirichlet-to-Neumann map).
\begin{definition}\label{mldc_defn:cinspo}
 Let $\alpha\in]0,1[$.  Let  $\Omega$ be a  bounded open subset of ${\mathbb{R}}^{n}$ of class $C^{1,\alpha}$. The classical interior Steklov-Poincar\'{e} operator is defined to be the operator $S_{\Omega,+}$ from
 \begin{equation}\label{mldc_defn:cinspo1}
C^{1,\alpha}(\partial\Omega)\quad\text{to}\quad C^{0,\alpha}(\partial\Omega)
\end{equation}
that takes $v\in C^{1,\alpha}(\partial\Omega)$ to the function 
 \begin{equation}\label{mldc_defn:cinspo2}
S_{\Omega,+}[v](x)\equiv \frac{\partial  }{\partial\nu}{\mathcal{G}}_{\Omega,d,+}[v](x)\qquad\forall x\in\partial\Omega\,.
\end{equation}
 \end{definition}
   Since   the classical normal derivative is continuous from $C^{1,\alpha}(\overline{\Omega})$ to $C^{0,\alpha}(\partial\Omega)$, the continuity of ${\mathcal{G}}_{\Omega,d,+}$ implies  that $S_{\Omega,+}[\cdot]$ is linear and continuous from 
  $C^{1,\alpha}(\partial\Omega)$ to $C^{0,\alpha}(\partial\Omega)$. Then we have the following definition of \cite[(41)]{La24c}.
  \begin{definition}\label{mldc_defn:conoderdedu}
 Let $\alpha\in]0,1[$.  Let  $\Omega$ be a  bounded open subset of ${\mathbb{R}}^{n}$ of class $C^{1,\alpha}$. If  $u\in C^{0}(\overline{\Omega})$ and $\Delta u\in  C^{-1,\alpha}(\overline{\Omega})$, then we define the distributional  normal derivative $\partial_\nu u$
 of $u$ to be the only element of the dual $(C^{1,\alpha}(\partial\Omega))'$ that satisfies the following equality
 \begin{equation}\label{mldc_defn:conoderdedu1}
\langle \partial_\nu u ,v\rangle \equiv\int_{\partial\Omega}uS_{\Omega,+}[v]\,d\sigma
+\langle E^\sharp_\Omega[\Delta u],{\mathcal{G}}_{\Omega,d,+}[v]\rangle 
\qquad\forall v\in C^{1,\alpha}(\partial\Omega)\,.
\end{equation}
\end{definition}
 The normal derivative of Definition \ref{mldc_defn:conoderdedu} extends the classical one in the sense that if $u\in C^{1,\alpha}(\overline{\Omega})$, then under the assumptions on $\alpha$ and $\Omega$  of Definition \ref{mldc_defn:conoderdedu}, we have 
 \begin{equation}\label{mldc_lem:conoderdeducl1}
 \langle \partial_\nu u ,v\rangle \equiv\int_{\partial\Omega}\frac{\partial u}{\partial\nu}v\,d\sigma
  \quad\forall v\in C^{1,\alpha}(\partial\Omega)\,,
\end{equation}
where $\frac{\partial u}{\partial\nu}$ in the right hand side denotes the classical normal derivative of $u$ on $\partial\Omega$ (cf.~\cite[Lem.~5.5]{La24c}). In the sequel, we use the classical symbol $\frac{\partial u}{\partial\nu}$ also for  $\partial_\nu u$ when no ambiguity can arise.\par  
  
Next we introduce the  function space $V^{-1,\alpha}(\partial\Omega)$ on the boundary of $\Omega$ for the normal derivatives of the functions of $C^{0,\alpha}(\overline{\Omega})_\Delta$ as in  \cite[Defn.~13.2, 15.10, Thm.~18.1]{La24b}.
\begin{definition}\label{mldc_defn:v-1a}
Let   $\alpha\in ]0,1[$. Let $\Omega$ be a bounded open  subset of ${\mathbb{R}}^{n}$ of class $C^{1,\alpha}$. Let 
\begin{eqnarray}\label{mldc_defn:v-1a1}
 \lefteqn{V^{-1,\alpha}(\partial\Omega)\equiv \biggl\{\mu_0+S_{\Omega,+}^t[\mu_1]:\,\mu_0, \mu_1\in C^{0,\alpha}(\partial\Omega)
\biggr\}\,,
}
\\ \nonumber
\lefteqn{
\|\tau\|_{  V^{-1,\alpha}(\partial\Omega) }
\equiv\inf\biggl\{\biggr.
 \|\mu_0\|_{ C^{0,\alpha}(\partial\Omega)  }+\|\mu_1\|_{ C^{0,\alpha}(\partial\Omega)  }
:\,
 \tau=\mu_0+S_{\Omega,+}^t[\mu_1]\biggl.\biggr\}\,,
 }
 \\ \nonumber
 &&\qquad\qquad\qquad\qquad\qquad\qquad\qquad\qquad\qquad
 \forall \tau\in  V^{-1,\alpha}(\partial\Omega)\,,
\end{eqnarray}
where $S_{\Omega,+}^t$ is the transpose map of $S_{\Omega,+}$.
\end{definition}
As shown in \cite[\S 13]{La24b},  $(V^{-1,\alpha}(\partial\Omega), \|\cdot\|_{  V^{-1,\alpha}(\partial\Omega)  })$ is a Banach space. By definition of the norm, $C^{0,\alpha}(\partial\Omega)$ is continuously embedded into $V^{-1,\alpha}(\partial\Omega)$. Moreover, we have the following statement  of \cite[Prop.~6.6]{La24c}  on the continuity of the normal derivative on $C^{0,\alpha}(\overline{\Omega})_\Delta$. 
\begin{proposition}\label{mldc_prop:ricodnu}
 Let   $\alpha\in ]0,1[$. Let $\Omega$ be a bounded open  subset of 
 ${\mathbb{R}}^{n}$ of class $C^{1,\alpha}$. Then the distributional normal derivative   $\partial_\nu$ is a continuous surjection of $C^{0,\alpha}(\overline{\Omega})_\Delta$ onto $V^{-1,\alpha}(\partial\Omega)$ and there exists $Z\in {\mathcal{L}}\left(V^{-1,\alpha}(\partial\Omega),C^{0,\alpha}(\overline{\Omega})_\Delta\right)$ such that
 \begin{equation}\label{mldc_prop:ricodnu1}
\partial_\nu Z[g]=g\qquad\forall g\in V^{-1,\alpha}(\partial\Omega)\,,
\end{equation}
\textit{i.e.}, $Z$ is a right inverse of  $\partial_\nu$.  
\end{proposition}
We also mention that the condition of Definition \ref{mldc_defn:conoderdedu} admits a different formulation at least in case $u\in C^{0,\alpha}(\overline{\Omega})_\Delta$ (see \cite[Prop.~5.15]{La24c}).
\begin{proposition}\label{mldc_prop:node1adeq}
 Let   $\alpha\in ]0,1[$. Let $\Omega$ be a bounded open  subset of 
 ${\mathbb{R}}^{n}$ of class $C^{1,\alpha}$. Let $E_\Omega$ be a linear map from $C^{1,\alpha}(\partial\Omega)$ to $ C^{1,\alpha}(\overline{\Omega})_\Delta$ such that
 \begin{equation}\label{mldc_prop:node1adeq0}
 E_\Omega[f]_{|\partial\Omega}=f\qquad\forall f\in C^{1,\alpha}(\partial\Omega)\,.
 \end{equation}
 If $u\in C^{0,\alpha}(\overline{\Omega})_\Delta$, then the distributional  normal derivative $\partial_\nu u$
 of $u$  is characterized by the validity of the following equality
 \begin{eqnarray}\label{mldc_prop:node1adeq1}
 \lefteqn{
\langle \partial_\nu u ,v\rangle =\int_{\partial\Omega}u
\frac{\partial}{\partial\nu}E_\Omega[v]
\,d\sigma
}
\\ \nonumber
&&\qquad\qquad
+\langle E^\sharp_\Omega[\Delta u],E_\Omega[v]\rangle -\int_\Omega\Delta (E_\Omega[v]) u\,dx
\qquad\forall v\in C^{1,\alpha}(\partial\Omega)\,.
\end{eqnarray}
\end{proposition}
\begin{remark}\label{rem:exopalex}
   We note that in equality (\ref{mldc_prop:node1adeq1}) we can use the `extension' operator $E_\Omega$ that we prefer and that accordingly equality (\ref{mldc_prop:node1adeq1}) is independent of the specific choice of $E_\Omega$. When we deal with problems for the Laplace operator, a good choice is
 $E_\Omega={\mathcal{G}}_{\Omega,d,+}$, so that the last term in the right hand side of (\ref{mldc_prop:node1adeq1}) disappears.  
 In particular, under the assumptions of Proposition \ref{mldc_prop:node1adeq}, an extension operator as $E_\Omega$ always exists.
\end{remark}

 \section{A  distributional outward unit normal derivative for the functions of $C^{0,\alpha}_{
{\mathrm{loc}}	}(\overline{\Omega^-})_\Delta$}
\label{mldc_sec:doun}
 
 We now plan to define the normal derivative on the boundary the functions of $C^{0,\alpha}_{
{\mathrm{loc}}	}(\overline{\Omega^-})_\Delta$ with $\alpha\in]0,1[$ in case $\Omega$ is a bounded open subset of ${\mathbb{R}}^n$ of class $C^{1,\alpha}$ as in \cite[\S 3]{La25}.  To do so, we  choose $r\in]0,+\infty[$ such that $\overline{\Omega}\subseteq {\mathbb{B}}_n(0,r)$ and we observe that the linear operator from $C^{1,\alpha}(\partial\Omega)$ to $C^{1,\alpha}((\partial\Omega)\cup(\partial{\mathbb{B}}_n(0,r)))$ that is defined by the equality
\begin{equation}\label{mldc_eq:zerom}
\stackrel{o}{E}_{\partial\Omega,r}[v]\equiv\left\{
\begin{array}{ll}
 v(x) & \text{if}\ x\in \partial\Omega\,,
 \\
 0 & \text{if}\ x\in \partial{\mathbb{B}}_n(0,r)\,,
\end{array}
\right.\qquad\forall v\in C^{1,\alpha}(\partial\Omega) 
\end{equation}
is linear and continuous and that accordingly the transpose map $\stackrel{o}{E}_{\partial\Omega,r}^t$ is linear and continuous
\[
\text{from}\ 
\left(C^{1,\alpha}((\partial\Omega)\cup(\partial{\mathbb{B}}_n(0,r)))\right)'\ 
\text{to}\ \left(C^{1,\alpha}(\partial\Omega)\right)'\,.
\]
Then we introduce the following definition as in \cite[Defn.~3.1]{La25}.
\begin{definition}\label{mldc_defn:endedr}
 Let   $\alpha\in ]0,1[$. Let $\Omega$ be a bounded open  subset of 
 ${\mathbb{R}}^{n}$ of class $C^{1,\alpha}$. Let $r\in]0,+\infty[$ be such that $\overline{\Omega}\subseteq {\mathbb{B}}_n(0,r)$. If $u\in C^{0,\alpha}_{
{\mathrm{loc}}	}(\overline{\Omega^-})_\Delta$, then we set
\begin{equation}\label{mldc_defn:endedr1}
 \partial_{\nu_{\Omega^-}}u
 \equiv 
 \stackrel{o}{E}_{\partial\Omega,r}^t\left[
\partial_{\nu_{{\mathbb{B}}_n(0,r)\setminus\overline{\Omega}}}u
\right]\,.
\end{equation}
 \end{definition}
 Under the assumptions of Definition \ref{mldc_defn:endedr},   $\partial_{\nu_{\Omega^-}}u$ is an element of $\left(C^{1,\alpha}(\partial\Omega)\right)'$. 
Under the assumptions of Definition \ref{mldc_defn:endedr}, there exists  a linear (extension) map $E_{{\mathbb{B}}_n(0,r)\setminus\overline{\Omega}}$ from $C^{1,\alpha}((\partial\Omega)\cup(\partial
{\mathbb{B}}_n(0,r)
))$ to $ C^{1,\alpha}(\overline{{\mathbb{B}}_n(0,r)}\setminus\Omega)_\Delta$ such that
 \begin{equation}\label{mldc_prop:node1adeq0r}
 E_{{\mathbb{B}}_n(0,r)\setminus\overline{\Omega}}[f]_{|(\partial\Omega)\cup(\partial
{\mathbb{B}}_n(0,r)
)}=f\qquad\forall f\in C^{1,\alpha}((\partial\Omega)\cup(\partial
{\mathbb{B}}_n(0,r)
))\,,
 \end{equation}
 (cf.~Remark \ref{rem:exopalex}).
Then the definition of normal derivative on the boundary of 
${\mathbb{B}}_n(0,r)\setminus\overline{\Omega}$ implies that
\begin{eqnarray}\label{mldc_defn:endedr2}
\lefteqn{
\langle \partial_{\nu_{\Omega^-}}u,v\rangle 
=-\int_{\partial\Omega}u\frac{\partial}{\partial\nu_{\Omega}}E_{{\mathbb{B}}_n(0,r)\setminus\overline{\Omega}}[\stackrel{o}{E}_{\partial\Omega,r}[v]]\,d\sigma
}
\\ \nonumber
&&\qquad\qquad 
+\int_{\partial{\mathbb{B}}_n(0,r)}u\frac{\partial}{\partial\nu_{{\mathbb{B}}_n(0,r)}}E_{{\mathbb{B}}_n(0,r)\setminus\overline{\Omega}}[\stackrel{o}{E}_{\partial\Omega,r}[v]]\,d\sigma
\\ \nonumber
&&\qquad\qquad 
+\langle E^\sharp_{{\mathbb{B}}_n(0,r)\setminus\overline{\Omega}}[\Delta u],E_{{\mathbb{B}}_n(0,r)\setminus\overline{\Omega}}[\stackrel{o}{E}_{\partial\Omega,r}[v]]\rangle 
\\ \nonumber
&&\qquad\qquad 
-\int_{{\mathbb{B}}_n(0,r)\setminus\overline{\Omega}}\Delta (E_{{\mathbb{B}}_n(0,r)\setminus\overline{\Omega}}[\stackrel{o}{E}_{\partial\Omega,r}[v]]) u\,dx
\qquad\forall v\in C^{1,\alpha}(\partial\Omega)\,.
\end{eqnarray}
 (cf.~(\ref{mldc_prop:node1adeq1})).   
As shown in \cite[Prop.~3.2]{La25},  Definition  \ref{mldc_defn:endedr} is independent of the choice of $r\in]0,+\infty[$ such that $\overline{\Omega}\subseteq {\mathbb{B}}_n(0,r)$.  

\begin{remark}\label{mldc_rem:conoderdeducl-} 
 Let   $\alpha\in ]0,1[$. Let $\Omega$ be a bounded open  subset of 
 ${\mathbb{R}}^{n}$ of class $C^{1,\alpha}$, $u\in C^{1,\alpha}_{
{\mathrm{loc}}	}(\overline{\Omega^-})$. Then equality (\ref{mldc_defn:endedr1}) of  the  
 of Definition \ref{mldc_defn:endedr} of normal derivative and equality (\ref{mldc_lem:conoderdeducl1}) imply that
 \begin{equation}\label{mldc_rem:conoderdeducl-1}
 \langle \partial_{\nu_{\Omega^-}} u ,v\rangle \equiv\int_{\partial\Omega}\frac{\partial u}{\partial\nu_{\Omega^-}}v\,d\sigma
  \quad\forall v\in C^{1,\alpha}(\partial\Omega)\,,
\end{equation}
where $\frac{\partial u}{\partial\nu_{\Omega^-}}$ in the right hand side denotes the classical $\nu_{\Omega^-}$-normal derivative of $u$ on $\partial\Omega$.  In the sequel, we use the classical symbol $\frac{\partial u}{\partial\nu_{\Omega^-}}$ also for  $\partial_{\nu_{\Omega^-}} u$ when no ambiguity can arise.\par 
\end{remark}

Next we note that if  $u\in C^{0,\alpha}_{{\mathrm{loc}}}(\overline{\Omega^-})$ and $u$ is both harmonic in $\Omega^-$ and harmonic at infinity, then $u\in C^{0,\alpha}_{
{\mathrm{loc}}	}(\overline{\Omega^-})_\Delta$. As shown in \cite[Prop.~3.4]{La25}, 
  the distribution $\partial_{\nu_{\Omega^-}} u $ of Definition \ref{mldc_defn:endedr} coincides with the normal derivative that has been introduced in \cite[Defn.~6.4]{La24b} for harmonic functions in $\Omega^-$ that are harmonic at infinity.   By \cite[Thm.~1.1]{La25}, we have the following continuity statement for the distributional normal derivative $\partial_{\nu_{\Omega^-}}$ of
Definition \ref{mldc_defn:endedr}.

\begin{proposition}\label{mldc_prop:recodnu}
 Let   $\alpha\in ]0,1[$. Let $\Omega$ be a bounded open  subset of 
 ${\mathbb{R}}^{n}$ of class $C^{1,\alpha}$.  Then the distributional normal derivative   $\partial_{\nu_{\Omega^-}}$ is a continuous surjection of $C^{0,\alpha}_{
{\mathrm{loc}}	}(\overline{\Omega^-})_\Delta$ onto $V^{-1,\alpha}(\partial\Omega)$ and there exists a linear and continuous  operator   $Z_-$ from  $ V^{-1,\alpha}(\partial\Omega)$ to $C^{0,\alpha}_{
{\mathrm{loc}}	}(\overline{\Omega^-})_\Delta $ such that
 \begin{equation}\label{mldc_prop:recodnu1}
\partial_{\nu_{\Omega^-}} Z_-[g]=g\qquad\forall g\in V^{-1,\alpha}(\partial\Omega)\,,
\end{equation}
\textit{i.e.}, $Z_-$ is a right inverse of  the operator $\partial_{\nu_{\Omega^-}}$. (See Lemma \ref{mldc_lem:c1alcof} for the topology of $C^{0,\alpha}_{
{\mathrm{loc}}	}(\overline{\Omega^-})_\Delta $).	 
\end{proposition}

\section{Preliminaries on the  acoustic potentials}
Let  $\alpha\in]0,1]$. Let  $\Omega$ be a bounded open subset of ${\mathbb{R}}^{n}$ of class $C^{1,\alpha}$. Let  $r_{|\partial\Omega}$  be the restriction map  from ${\mathcal{D}}({\mathbb{R}}^n)$ to $C^{1,\alpha}(\partial\Omega)$. Let $\lambda\in {\mathbb{C}}$. If $S_{n,\lambda} $ is a fundamental solution 
of the operator $\Delta+\lambda$ and $\mu\in (C^{1,\alpha}(\partial\Omega))'$, then the  (distributional)  acoustic  single layer  potential relative to $S_{n,\lambda} $ and $\mu$ is the distribution
\[
v_\Omega[S_{n,\lambda} ,\mu]=(r_{|\partial\Omega}^t\mu)\ast S_{n,\lambda}  \in {\mathcal{D}}'({\mathbb{R}}^n) 
\]
 and we  set
\begin{eqnarray}\label{mldc_eq:dsila}
v_\Omega^+[S_{n,\lambda} ,\mu]&\equiv&\left((r_{|\partial\Omega}^t\mu)\ast S_{n,\lambda}  \right)_{|\Omega}
\qquad\text{in}\ \Omega\,,
\\ \nonumber
v_\Omega^-[S_{n,\lambda} ,\mu] &\equiv&
\left((r_{|\partial\Omega}^t\mu)\ast S_{n,\lambda}  \right)_{|\Omega^-}
\qquad\text{in}\ \Omega^-\,.
\end{eqnarray}
It is also known that the restriction of $v_\Omega[S_{n,\lambda} ,\mu]$ to ${\mathbb{R}}^n\setminus\partial\Omega$ equals the (distribution that is associated to) the function
\[
 \langle (r_{|\partial\Omega}^t\mu)(y),S_{n,\lambda} (\cdot-y)\rangle \,.
\]
In the cases in which both $v_\Omega^+[S_{n,\lambda} ,\mu]$ and $v_\Omega^-[S_{n,\lambda} ,\mu]$ admit a continuous extension to $\overline{\Omega}$ and to $\overline{\Omega^-}$, respectively, we still use the symbols $v_\Omega^+[S_{n,\lambda} ,\mu]$ and $v_\Omega^-[S_{n,\lambda} ,\mu]$ for the continuous extensions and if the values of $v_\Omega^\pm[S_{n,\lambda} ,\mu](x)$ coincide for each $x\in\partial\Omega$, then we set
\[
V_\Omega[S_{n,\lambda} ,\mu](x)\equiv v_\Omega^+[S_{n,\lambda} ,\mu](x)=v_\Omega^-[S_{n,\lambda} ,\mu]^-(x)\qquad\forall x\in\partial\Omega\,.
\]
If $\mu$ is continuous, then it is known that $v_\Omega[S_{n,\lambda} ,\mu ]$ is continuous in ${\mathbb{R}}^n$ (cf.~\textit{e.g.}, \cite[Lem.~4.2 (i), Lem.~6.2]{DoLa17}). For the classical  properties of the acoustic single layer potential, we refer for example to  the paper \cite{DoLa17} of the author and Dondi and to \cite[\S 5]{La25a}. Next we introduce the following (classical) technical statement on the  acoustic  double layer potential (cf.~\textit{e.g.}, \cite[Thm.~5.1]{La25a}).
\begin{theorem}\label{mldc_thm:dlay}
 Let  $\alpha\in]0,1[$. Let $\Omega$ be a bounded open subset of  ${\mathbb{R}}^n$ of class $C^{1,\alpha}$. Let $\lambda\in {\mathbb{C}}$.  Let $S_{n,\lambda} $ be a fundamental solution 
of the operator $\Delta+\lambda$. If $\mu\in C^{0,\alpha}(\partial\Omega)$, and
\begin{equation}\label{mldc_thm:dlay1}
w_\Omega[S_{n,\lambda} ,\mu](x)\equiv\int_{\partial\Omega}\frac{\partial}{\partial\nu_{\Omega,y}}\left(S_{n,\lambda} (x-y)\right)\mu(y)\,d\sigma_y\qquad\forall x\in {\mathbb{R}}^n\,,
\end{equation}
where
\[
\frac{\partial}{\partial \nu_{\Omega,y} }
\left(S_{n,\lambda}(x-y)\right)\equiv
  -DS_{n,\lambda}(x-y) \nu_{\Omega}(y) \qquad\forall (x,y)\in\mathbb{R}^n\times\partial\Omega\,, x\neq y\,,
 \]
then the restriction 
$w_\Omega[S_{n,\lambda} ,\mu]_{|\Omega}$ can be extended uniquely to a function  
$w^{+}_\Omega [S_{n,\lambda} ,\mu]$ of class $ C^{0,\alpha}(\overline{\Omega})$ and $w_\Omega[S_{n,\lambda} ,\mu]_{|{\mathbb{R}}^n\setminus\overline{\Omega}}$ can be extended uniquely to a function  
$w^{-}_\Omega [S_{n,\lambda} ,\mu]$ of class  $
C^{0,\alpha}_{ {\mathrm{loc}} }(\overline{\Omega^{-}})$. Moreover, we have the following jump relation
\begin{equation}\label{mldc_thm:dlay1a}
w^{\pm}_\Omega [S_{n,\lambda},\mu](x)
=\pm\frac{1}{2}\mu(x)+w_\Omega[S_{n,\lambda},\mu](x)
\qquad\forall x\in\partial\Omega\,.
\end{equation}
 \end{theorem}
For the classical  properties of the acoustic double layer potential, we refer for example to  the paper \cite{DoLa17} of the author and Dondi and to \cite[\S 5]{La25a}. We also set
\[
W_\Omega[S_{n,\lambda} ,\mu](x)\equiv w_\Omega[S_{n,\lambda} ,\mu](x)\qquad\forall x\in\partial\Omega\,.
\]
In order to shorten our notation, we introduce the following abbreviation. Let $m\in {\mathbb{N}}$,   $\alpha\in]0,1[$. Let $\Omega$ be a bounded open subset of ${\mathbb{R}}^{n}$ of class $C^{\max\{1,m\},\alpha}$. Then we set
\begin{equation}\label{mldc_eq:vm-1a}
V^{m-1,\alpha}(\partial\Omega)\equiv\left\{
\begin{array}{ll}
 C^{m-1,\alpha}(\partial\Omega)& \text{if}\ m\geq 1\,,
 \\
 V^{-1,\alpha}(\partial\Omega)& \text{if}\ m=0\,.
\end{array}
\right.
\end{equation}
\begin{remark}\label{mldc_rem:wtnotation}
  Let  $m\in{\mathbb{N}}$, $\alpha\in]0,1[$. Let $\Omega$ be a bounded open subset of ${\mathbb{R}}^{n}$ of class $C^{\max\{1,m\},\alpha}$. Let $\lambda\in {\mathbb{C}}$. Let $S_{n,\lambda}$ be a fundamental solution of $\Delta+\lambda$. Then
 $W_\Omega[S_{n,\lambda},\cdot]$ is known to be compact in $C^{\max\{1,m\},\alpha}(\partial\Omega)$ (cf.~\textit{e.g.},   \cite[Cor.~9.1]{DoLa17}) and $W_\Omega^t[S_{n,\lambda},\cdot]$ denotes the transpose map to the operator $W_\Omega[S_{n,\lambda},\cdot]$ with respect to the natural duality pairing
\begin{equation}\label{mldc_rem:wtnotation1}
 \left(
 \left(C^{\max\{1,m\},\alpha}(\partial\Omega)\right)',C^{\max\{1,m\},\alpha}(\partial\Omega)
 \right)\,.
 \end{equation}
 Via the identification map ${\mathcal{J}}$ of Lemma \ref{mldc_lem:caincl} in caqse $m\geq 1$ and by definition in case $m=0$, $V^{m-1,\alpha}(\partial\Omega)$
 can be regarded as a subspace of $\left(C^{\max\{1,m\},\alpha}(\partial\Omega)\right)'$ and 
 it is   known that $W_\Omega^t[S_{n,\lambda},\cdot]$ is compact in $V^{m-1,\alpha}(\partial\Omega)$ (cf.~\textit{e.g.}, \cite[Cor.~10.1]{DoLa17} in case $m\geq 1$  and \cite[Cor.~8.5]{La25a} in case $m=0$). 
 Then    $W_\Omega^t[S_{n,\lambda},\cdot]_{|V^{m-1,\alpha}(\partial\Omega)}$ denotes the transpose map to the operator $W_\Omega[S_{n,\lambda},\cdot]$ with respect to the  duality pairing
\begin{equation}\label{mldc_rem:wtnotation2}
 \left(
 V^{m-1,\alpha}(\partial\Omega),C^{\max\{1,m\},\alpha}(\partial\Omega)
 \right)\,,
\end{equation}
 (cf.~\textit{e.g.}, Kress~\cite[Defn.~4.5]{Kr14} for the definition of transpose map with respect to a duality pairing).
\end{remark}

\section{A representation theorem for the interior Diri\-chlet eigenfunc\-tions}\label{mldc_sec:ridiefh}

Let $k\in {\mathbb{C}}\setminus]-\infty,0]$, ${\mathrm{Im}}\, k\geq 0$. As  well-known	  in scattering theory, a function $u$ of class $C^1$ in the complement of a compact subset of ${\mathbb{R}}^{n}$ satisfies the outgoing $k$-radiation condition  provided that
\begin{equation}
\label{mldc_rad1}
\lim_{x\to\infty}|x|^{\frac{n-1}{2}}(Du(x)\frac{x}{|x|}-iku(x))=0\,.
\end{equation}
Then we set 
\begin{equation}\label{mldc_eq:sofura}
\tilde{S}_{n,k;r}(x)=C_n|x|^{-\frac{n-2}{2}}H^{(1)}_{ \frac{n-2}{2}}(k|x|)
\qquad\forall x\in {\mathbb{R}}^n\setminus\{0\}\,,
\end{equation}
where $H^{(1)}_{ \frac{n-2}{2}}$ is the first Hankel function,
\begin{equation}\label{mldc_eq:sofura1}
C_{n}\equiv\left\{
\begin{array}{ll}
\frac{	k^{ \frac{n-2}{2} } }{i4(2\pi)^{\frac{n-2}{2}}}
& {\text{if}}\ n\ {\text{is\ even}}\,,
\\
\frac{e^{  \frac{n-2}{2}  \log (k)}}{i4(2\pi)^{\frac{n-2}{2}}}
& {\text{if}}\ n\ {\text{is\ odd}}\,,
\end{array}
\right.
\end{equation}
and $\log$  denotes the principal branch of  the logarithm in $ {\mathbb{C}}\setminus]-\infty,0]$
(see \cite[\S 4]{La25b},  Mitrea \cite[(7.6.11), p.~266]{Mit18}, Mitrea, Mitrea and Mitrea \cite[(6.1.7), p.~882]{MitMitMit23}). Thus $\tilde{S}_{n,k;r}$ is the fundamental solution of 
$\Delta+k^2$ that satisfies the outgoing $k$-radiation condition  and  that is classically used in scattering theory. The subscript $r$ stands for `radiation'. 
Next, we introduce the following elementary remark (cf.~\cite[Rem.~9.1]{La25b}). 
\begin{remark}\label{mldc_rem:helpo}
 Let $m\in {\mathbb{N}}$, $\alpha\in]0,1[$, $\lambda\in {\mathbb{C}}$. 	 	Let $\Omega$ be a bounded open subset of ${\mathbb{R}}^{n}$ of class $C^{\max\{1,m\},\alpha}$.Then the following statements hold. 
 \begin{enumerate}
\item[(i)]  If $m\in \{0,1\}$ and  $u\in C^{m,\alpha}(\overline{\Omega})$ satisfies equation $\Delta u+ \lambda 	u=0$ in $\Omega$, then $u\in C^{m,\alpha}(\overline{\Omega})_\Delta$.
 \item[(ii)] If $m\geq 2$, then $C^{m,\alpha}(\overline{\Omega})$ is continuously embedded into  $C^{1,\alpha}(\overline{\Omega})_\Delta$.
 \item[(iii)] If $m\in \{0,1\}$ and  $u\in C^{m,\alpha}_{{\mathrm{loc}}}(\overline{\Omega^-})$ satisfies equation $\Delta u+ \lambda	 	u=0$ in $\Omega^-$, then $u\in C^{m,\alpha}_{{\mathrm{loc}}}(\overline{\Omega^-})_\Delta$.
 \item[(iv)] If $m\geq 2$, then $C^{m,\alpha}_{{\mathrm{loc}}}(\overline{\Omega^-})$  is continuously embedded into $C^{1,\alpha}_{{\mathrm{loc}}}(\overline{\Omega^-})_\Delta$.\par 
\end{enumerate}
\end{remark}

Next we show the following representation theorem for the interior Dirichlet eigenfunctions that is classical for $m\geq 1$. Here we note that  the normal derivative is to be interpreted in the sense of Definition \ref{mldc_defn:conoderdedu} if  $m=0$.

\begin{theorem}
\label{mldc_thm:eipd}
Let  $m\in{\mathbb{N}}$, $\alpha\in]0,1[$. Let $\Omega$ be a bounded open subset of ${\mathbb{R}}^{n}$ of class $C^{\max\{1,m\},\alpha}$. Then the following statements hold.
\begin{enumerate}
\item[(i)] Let  $\lambda\in {\mathbb{C}}$. If $S_{n,\lambda}$ is a fundamental solution of $\Delta+\lambda$ and if $u\in C^{m,\alpha}(\overline{\Omega}) $ satisfies the problem
\begin{equation}
\label{mldc_thm:eipd1}
\left\{
\begin{array}{ll}
\Delta u+\lambda u=0 &{\text{in}}\ \Omega\,,
\\
u=0
&{\text{on}}\ \partial\Omega\,
\end{array}
\right.
\end{equation}
then $u\in C^{\max\{1,m\},\alpha}(\overline{\Omega})$, 
\begin{equation}
\label{mldc_thm:eipd2}
u(x)\equiv 
-v_\Omega^+[S_{n,\lambda},\frac{\partial u}{\partial\nu_{\Omega}}]
\qquad\text{on}\ \overline{\Omega} \,,
\end{equation}
and 
\begin{equation}
\label{mldc_thm:eipd3}
\frac{1}{2}\frac{\partial u}{\partial\nu_{\Omega}} 
+
 W_\Omega^t[S_{n,\lambda} ,\frac{\partial u}{\partial\nu_{\Omega}}]
=0\qquad\text{on}\ \partial\Omega\,.
\end{equation}
\item[(ii)] Let $k\in {\mathbb{C}}\setminus ]-\infty,0]$, ${\mathrm{Im}}\,k\geq 0$. If 
 for each $j\in {\mathbb{N}}\setminus\{0\}$ such that $j\leq \varkappa^{-}$, $k^{2}$ is not a Neumann eigenvalue for $-\Delta$ in $(\Omega^-)_{j}$, and if $\mu\in V^{m-1,\alpha}(\partial\Omega)$ (cf.~(\ref{mldc_eq:vm-1a}))
  satisfies the equation
\begin{equation}
\label{mldc_thm:eipd4}
\frac{1}{2}\mu(x)
+
 W_\Omega^t[\tilde{S}_{n,k;r},\mu]
 =0\qquad\text{on}\ \partial\Omega\,,
\end{equation}
then the function $u$ from $\overline{\Omega}$ to ${\mathbb{C}}$ defined by 
\begin{equation}
\label{mldc_thm:eipd5}
u(x)= -v_\Omega^+[\tilde{S}_{n,k;r},\mu] (x) \qquad\forall x\in \overline{\Omega} \,,
\end{equation}
belongs to $C^{\max\{1,m\},\alpha}(\overline{\Omega})$ and satisfies the Dirichlet problem (\ref{mldc_thm:eipd1}) with $\lambda=k^{2}$, and 
\begin{equation}
\label{mldc_thm:eipd6}
\frac{\partial u}{\partial\nu_{\Omega}}=\mu\qquad{\text{on}}\ \partial
\Omega\,.
\end{equation}
In particular,  $\mu\in C^{\max\{1,m\}-1,\alpha}(\partial\Omega) $. Moreover, $-v_\Omega^-[\tilde{S}_{n,k;r},\mu]=0$ in $\Omega^{-}$.  
\end{enumerate}
\end{theorem}
{\bf Proof.} (i) The membership of $u$ in $C^{\max\{1,m\},\alpha}(\overline{\Omega})$ holds by assumption in case $m\geq 1$ and is a consequence of Remark \ref{mldc_rem:helpo} and of
the regularization theorem of  \cite[Thm.~1.2 (iii)]{La25} in case $m=0$. Equality (\ref{mldc_thm:eipd2}) is an immediate consequence of the classical third Green identity (cf.~\textit{e.g.}, \cite[
Thm.~9.2 (ii)]{La25a})  and of the continuity in $\overline{\Omega}$ of the  acoustic  simple layer potential  (cf.~\textit{e.g.},   \cite[Thm.~7.1]{DoLa17}). By taking the normal derivative of (\ref{mldc_thm:eipd2}) on $\partial\Omega$ and by exploiting the  jump formulas for the acoustic simple layer potential, we obtain
\[
\frac{\partial u}{\partial\nu_{\Omega}} 
=\frac{1}{2}\frac{\partial u}{\partial\nu_{\Omega}} 
-
W_\Omega^t[S_{n,\lambda} ,\frac{\partial u}{\partial\nu_{\Omega}}]
\qquad\qquad\text{on}\ \partial\Omega\,,
\]
(cf.~\textit{e.g.},  \cite[Thm.~7.1 (ii)]{DoLa17}) and thus equality (\ref{mldc_thm:eipd3}) follows. 

We now prove statement (ii). We set   $u^{\pm}\equiv -v^{\pm}_\Omega[\tilde{S}_{n,k;r},\mu]$. 
If $m=0$, then \cite[Thm.~5.3]{La25a} implies that $u^+\in C^{0,\alpha}(\overline{\Omega})_\Delta $, 
$u^-\in C^{0,\alpha}_{{\mathrm{loc}}}(\overline{\Omega^-})_\Delta$. If $m\geq 1$, then it is classically known that $u^+\in C^{m,\alpha}(\overline{\Omega}) $, 
$u^-\in C^{m,\alpha}_{{\mathrm{loc}}}(\overline{\Omega^-})$.
 (cf.~\textit{e.g.},  \cite[Thm.~7.1]{DoLa17}).  By equation (\ref{mldc_thm:eipd4}) and by the jump formulas for the normal derivative of the acoustic simple layer potential of \cite[Thm.~6.3]{La25a}, we have 
 $\frac{\partial u^{-}}{\partial\nu_{\Omega^-}}=0$ on $\partial\Omega$. Then $u^-$ solves the Neumann problem for the Helmholtz equation in $\Omega^{-}_{j}$  and our assumption implies that
   $u^{-}=0$ on $(\Omega^{-})_{j}$ for each $j\in {\mathbb{N}}\setminus\{0\}$ such that $j\leq \varkappa^{-}$. On the other hand, $u^{-}$ satisfies the outgoing $k$-radiation condition  (cf.~\cite[Thm.~6.20 (i)]{La25b}) and thus $u^{-}=0$ on $(\Omega^{-})_{0}$ (cf.~\cite[Thm.~1.3]{La25}). Accordingly, $u^{-}=0$ on $\Omega^{-}$ and thus $u^{+}_{|\partial\Omega}=u^{-}_{|\partial\Omega}=0$ on $\partial\Omega$, and thus $u^{+}$ satisfies the Dirichlet problem (\ref{mldc_thm:eipd1}) with $\lambda=k^{2}$ (cf.~\cite[Thm.~6.3]{La25a}). Then statement (i) implies that $u^+\in C^{\max\{1,m\},\alpha}(\overline{\Omega})$. By equation (\ref{mldc_thm:eipd4}) and by the jump formulas for the acoustic simple layer potential of 
   \cite[Thm.~6.3]{La25a}, we also have 
\[
\mu=\frac{\partial u^{-}}{\partial\nu_{\Omega^-}}+\frac{\partial u^{+}}{\partial\nu_{\Omega}}=
    \frac{\partial u^{+}}{\partial\nu_{\Omega}}
   \]
    and thus $\mu\in C^{\max\{1,m\}-1,\alpha}(\partial\Omega) $ and equality (\ref{mldc_thm:eipd6}) holds true.\hfill  $\Box$ 

\vspace{\baselineskip}

\begin{corollary}
\label{mldc_corol:eipdd}
Let  $m\in{\mathbb{N}}$, $\alpha\in]0,1[$. Let $\Omega$ be a bounded open subset of ${\mathbb{R}}^{n}$ of class $C^{\max\{1,m\},\alpha}$. Then the following statements hold. 
\begin{enumerate}
\item[(i)] Let  $\lambda\in {\mathbb{C}}$. Let $S_{n,\lambda}$ is a fundamental solution of $\Delta+\lambda$. The operator $\frac{\partial}{\partial\nu_{\Omega}}$ is a linear injection from the space 
\begin{eqnarray}
\label{mldc_corol:eipdd1}
\lefteqn{
\{
u\in C^{m,\alpha}(\overline{\Omega}):\,\Delta u+\lambda u=0, u_{|\partial\Omega}=0\}
}
\\ \nonumber
&&\qquad
=\{
u\in C^{\max\{1,m\},\alpha}(\overline{\Omega}):\,\Delta u+\lambda u=0, u_{|\partial\Omega}=0\}
\end{eqnarray}
into the space
 \begin{equation}
\label{mldc_corol:eipdd2}
\{
\mu\in C^{\max\{1,m\}-1,\alpha}(\partial\Omega):\,\frac{1}{2}\mu 
+
 W_\Omega^t[S_{n,\lambda},\mu]
 =0\}\,.
 \end{equation}
Moreover,
\begin{eqnarray}\label{mldc_corol:eipdd3}
\lefteqn{
\{
\mu\in C^{\max\{1,m\}-1,\alpha}(\partial\Omega):\,\frac{1}{2}\mu 
+
 W_\Omega^t[S_{n,\lambda},\mu]
 =0\}
}
\\ \nonumber
&&\qquad\qquad 
=
 \{\mu\in V^{m-1,\alpha}(\partial\Omega):\,\frac{1}{2}\mu 
+
 W_\Omega^t[S_{n,\lambda},\mu]
 =0\}
\end{eqnarray}
(cf.~(\ref{mldc_eq:vm-1a})).
\item[(ii)] Let $k\in {\mathbb{C}}\setminus ]-\infty,0]$, ${\mathrm{Im}}\,k\geq 0$. Assume that
 for each $j\in {\mathbb{N}}\setminus\{0\}$ such that $j\leq \varkappa^{-}$, $k^{2}$ is not a Neumann eigenvalue for $-\Delta$ in $(\Omega^-)_{j}$. Then the operator $\frac{\partial}{\partial\nu_{\Omega}}$ induces a linear isomorphism from the space 
 in (\ref{mldc_corol:eipdd1}) with $\lambda=k^2$ onto the space
 \begin{equation}\label{mldc_corol:eipdd4}
\{
\mu\in C^{\max\{1,m\}-1,\alpha}(\partial\Omega):\,\frac{1}{2}\mu 
+
 W_\Omega^t[\tilde{S}_{n,k;r},\mu]
 =0\}
\end{equation}
 and the inverse of such an isomorphism is the map $-v_\Omega^+[\tilde{S}_{n,k;r},\cdot]$. In particular, the dimension of the space in (\ref{mldc_corol:eipdd4}) equals the geometric multiplicity of $k^{2}$ as a Dirichlet eigenvalue of $-\Delta$ in $\Omega$. 
 \end{enumerate}
\end{corollary}
{\bf Proof.} (i) By Theorem \ref{mldc_thm:eipd} (i), equality (\ref{mldc_corol:eipdd1}) holds true.  
Then again Theorem \ref{mldc_thm:eipd} (i) implies that 
 $\frac{\partial u}{\partial\nu_{\Omega}}$ is a linear operator from the space in (\ref{mldc_corol:eipdd1}) to the space in (\ref{mldc_corol:eipdd2}). By the Holmgren Uniqueness Corollary \ref{mldc_corol:holmgr}  of the Appendix, $\frac{\partial u}{\partial\nu_{\Omega}}$ is injective in the space in (\ref{mldc_corol:eipdd1}). In order to show the equality in (\ref{mldc_corol:eipdd3}), it suffices to assume that $m=0$, that $\mu\in V^{-1,\alpha}(\partial\Omega)$ solves
equation 
\begin{equation}\label{mldc_corol:eipdd5}
\frac{1}{2}\mu 
+
 W_\Omega^t[S_{n,\lambda},\mu]
 =0
\end{equation}
 and to show that $\mu\in C^{0,\alpha}(\partial\Omega)$. If $\mu\in V^{-1,\alpha}(\partial\Omega)$ solves equation (\ref{mldc_corol:eipdd5}), then the jump formula  for the distributional normal derivative of the  acoustic  single layer potential of \cite[Thm.~6.3]{La25a} implies that 
 \begin{equation}\label{mldc_corol:eipdd6}
 \frac{\partial}{\partial\nu_{\Omega^-}}v_\Omega^-[S_{n,\lambda},\mu]=0\,.
\end{equation}
 By \cite[Thm.~5.3]{La25a} on the regularity of the  acoustic  single layer potential,  $v_\Omega^-[S_{n,\lambda},\mu]$ belongs to $ C^{0,\alpha}_{{\mathrm{loc}}}(\overline{\Omega^-})_\Delta$. Then  equality (\ref{mldc_corol:eipdd6}) and the regularization Theorem of \cite[Thm.~1.2  (iv)]{La25} for solutions of the Helmholtz equation implies that $v_\Omega^-[S_{n,\lambda},\mu]\in C^{1,\alpha}_{{\mathrm{loc}}}(\overline{\Omega^-})$. Since
 \[
 v_\Omega^+[S_{n,\lambda},\mu]_{|\partial\Omega}= v_\Omega^-[S_{n,\lambda},\mu]_{|\partial\Omega}\in C^{1,\alpha}(\partial\Omega)\,,
 \]
 and   $v_\Omega^+[S_{n,\lambda},\mu]$ in $C^{0,\alpha}(\overline{\Omega})_\Delta$
 (cf.~\cite[Thm.~5.3]{La25a}), then the regularization  Theorem of \cite[Thm.~1.2 (iii)]{La25} for solutions of the Helmholtz equation implies that $v_\Omega^+[S_{n,\lambda},\mu]$ belongs to $ C^{1,\alpha}(\overline{\Omega})$. Thus the the jump formulas for the distributional normal derivative of the  acoustic  single layer potential of \cite[Thm.~6.3]{La25a} imply that 
\[
\mu=-\frac{\partial}{\partial\nu_{\Omega^-}}v_\Omega^-[S_{n,\lambda},\mu]-
\frac{\partial}{\partial\nu_{\Omega}}v_\Omega^+[S_{n,\lambda},\mu]\in C^{0,\alpha}(\partial\Omega)\,.
\]
 Statement (ii) is an immediate consequence of statement (i) and of Theorem \ref{mldc_thm:eipd} (ii).\hfill  $\Box$ 

\vspace{\baselineskip}

Next we note that under appropriate assumptions, the space of solutions of the integral equation (\ref{mldc_thm:eipd4}) is stable under conjugation. 

\begin{remark}\label{mldc_rem:+i+wtconj}  
Let $m\in{\mathbb{N}}$, $\alpha\in]0,1[$. Let $\Omega$ be a bounded open subset of ${\mathbb{R}}^{n}$ of class $C^{\max\{1,m\},\alpha}$. Let 
$k\in {\mathbb{C}}\setminus ]-\infty,0]$, ${\mathrm{Im}}\,k\geq 0$. Assume that 
 for each $j\in {\mathbb{N}}\setminus\{0\}$ such that $j\leq \varkappa^{-}$, $k^{2}$ is not a Neumann eigenvalue for $-\Delta$ in $(\Omega^-)_{j}$.\par
 
 If  $\mu $ belongs to  $V^{m-1,\alpha}(\partial\Omega)$ and satisfies equation (\ref{mldc_thm:eipd4}), then $\mu$ belongs to    $C^{\max\{1,m\}-1,\alpha}(\partial\Omega)$
 and
$\overline{\mu}$   satisfies the same equation. Indeed, if $k^2$ is not real, then Corollary \ref{mldc_corol:eipdd} ensures that equation (\ref{mldc_thm:eipd4}) has only the trivial solution. 
If instead $k^2$ is real and if $\mu $ in $V^{m-1,\alpha}(\partial\Omega)$ satisfies equation (\ref{mldc_thm:eipd4}), then Corollary \ref{mldc_corol:eipdd} (i) implies that $\mu\in C^{\max\{1,m\}-1,\alpha}(\partial\Omega)$, 
Corollary \ref{mldc_corol:eipdd} (ii) implies that $v^{+}_\Omega[\tilde{S}_{n,k;r},-\mu]$ belongs to the set in (\ref{mldc_corol:eipdd1}) with $\lambda=k^2$ and that $\mu=\frac{\partial}{\partial\nu}v^{+}[\tilde{S}_{n,k;r},-\mu]$. Since $k^2$ is real, $\overline{v^{+}_\Omega[\tilde{S}_{n,k;r},-\mu]}$ belongs to the set in (\ref{mldc_corol:eipdd1}) with $\lambda=k^2$ too. In particular, $\overline{v^{+}_\Omega[\tilde{S}_{n,k;r},-\mu]}\in C^{\max\{1,m\},\alpha}(\overline{\Omega})$. Then
 Corollary \ref{mldc_corol:eipdd} (i) implies that $\frac{\partial}{\partial\nu}\overline{v^{+}[\tilde{S}_{n,k;r},-\mu]}$
satisfies equation (\ref{mldc_thm:eipd4}). Since
\[
 \frac{\partial}{\partial\nu}\overline{v^{+}[\tilde{S}_{n,k;r},-\mu] } 
 =
 \overline{ \frac{\partial}{\partial\nu}v^{+}[\tilde{S}_{n,k;r},-\mu] }=\overline{\mu}  \,,
\]
then $\overline{\mu}$ satisfies equation (\ref{mldc_thm:eipd4}).
\end{remark}
As is well known, if      $\alpha\in]0,1[$ and $\Omega$ is a bounded open subset of ${\mathbb{R}}^{n}$ of class $C^{1,\alpha}$, then problem (\ref{mldc_thm:eipd1}) has nontrivial solutions precisely  when $\lambda$ belongs to an increasing sequence $(\lambda_{j}^{({\mathcal{D}})}[\Omega])_{j\in{\mathbb{N}} }$ of   eigenvalues in $]0,+\infty[$, and we can write such a sequence so that each eigenvalue $\lambda_{j}^{({\mathcal{D}})}[\Omega]$ is repeated as many times as its multiplicity. 

\section{A nonvariational form of the interior Dirichlet problem for the Helmholtz equation}\label{mldc_sec:inpddh}

We now consider the interior  Dirichlet problem in the case in which the Dirichlet datum is in the space $C^{m,\alpha}(\partial\Omega)$ (nonvariational in case $m=0$) . For the corresponding classical result in case $\varkappa^-=0$, we refer to  Colton and Kress \cite[Thm.~3.24]{CoKr92}. We first introduce the following classical necessary condition for the existence.

\begin{theorem}\label{mldc_thm:nexintdirhe}
 Let $m\in {\mathbb{N}}$,   $\alpha\in]0,1[$. Let $\Omega$ be a bounded open subset of ${\mathbb{R}}^{n}$ of class $C^{\max\{1,m\},\alpha}$. 
Let  $\lambda\in {\mathbb{C}}$.  Let $g\in C^{m,\alpha}(\partial\Omega)$. If  the interior Dirichlet problem 
\begin{equation}
\label{mldc_thm:nexintdirhe1}
\left\{
\begin{array}{ll}
\Delta u+\lambda u=0 &{\text{in}}\ \Omega\,,
\\
u=g
&{\text{on}}\ \partial\Omega\,
\end{array}
\right.
\end{equation}
has a solution $u\in C^{m,\alpha}(\overline{\Omega})$, then
\begin{equation}\label{mldc_thm:nexintdirhe2}
\int_{\partial\Omega}g\frac{\partial v}{\partial\nu_\Omega}\,d\sigma=0  \,,
\end{equation}
for all $v\in C^{\max\{1,m\},\alpha}(\overline{\Omega})$  such that
\begin{equation}\label{mldc_thm:nexintdirhe3}
\Delta v+\lambda v=0 \qquad\text{in}\ \Omega\,,\qquad v=0\qquad\text{on}\ \partial\Omega\,.
\end{equation}
\end{theorem}
{\bf Proof.} Assume that a solution $u\in C^{m,\alpha}(\overline{\Omega})$  exists and that $v\in C^{\max\{1,m\},\alpha}(\overline{\Omega})$ satisfies the homogeneous Dirichlet problem of (\ref{mldc_thm:nexintdirhe3}). Then
Remark \ref{mldc_rem:helpo}, the second Green Identity in distributional form of \cite[Thm.~8.1]{La25b} and equality (\ref{mldc_prop:nschext3})   imply that 
 \begin{eqnarray*}
\lefteqn{
\int_{\partial\Omega}g\frac{\partial v}{\partial\nu_\Omega}\,d\sigma=
\int_{\partial\Omega}u\frac{\partial v}{\partial\nu_\Omega}\,d\sigma
}
\\ \nonumber
&&\qquad
=\langle  \frac{\partial u}{\partial\nu_\Omega}, v_{|\partial\Omega}\rangle
+\int_{\Omega}u(\Delta v+\lambda v)\,dx
-\langle E^\sharp_\Omega[\Delta u+\lambda u],v\rangle
=0\,.
\end{eqnarray*}\hfill  $\Box$ 

\vspace{\baselineskip}
By Theorem \ref{mldc_thm:nexintdirhe}, a necessary condition for the solvability of the Dirichlet problem is that condition (\ref{mldc_thm:nexintdirhe2}) holds true. Thus we now introduce the following notation. If $m\in{\mathbb{N}}$, $\alpha\in]0,1[$,  $\Omega$ is a bounded open subset of ${\mathbb{R}}^{n}$ of class $C^{\max\{1,m\},\alpha}$, $\lambda\in {\mathbb{C}}$, then we set
\begin{equation}\label{mldc_eq:xmal}
X^{m,\alpha}_\lambda\equiv\left\{
\frac{\partial v}{\partial\nu_\Omega}:\, v\in C^{\max\{1,m\},\alpha}(\overline{\Omega})\,,\ 
\Delta v+\lambda v=0 \ \text{in}\ \Omega\,,\  v=0\ \text{on}\ \partial\Omega\right\}
\end{equation}
and
\begin{equation}\label{mldc_eq:xmal1}
(X^{m,\alpha}_\lambda)^\perp\equiv
\left\{
g\in C^{m,\alpha}(\partial\Omega):\, \int_{\partial\Omega}gf\,d\sigma=0\ \forall f\in X^{m,\alpha}_\lambda
\right\}\,.
\end{equation}
\begin{lemma}\label{mldc_lem:xmalco}
 Let $m\in{\mathbb{N}}$, $\alpha\in]0,1[$. Let $\Omega$ be a bounded open subset of ${\mathbb{R}}^{n}$ of class $C^{\max\{1,m\},\alpha}$. Let $\lambda\in {\mathbb{C}}$. If
 $\mu\in  X^{m,\alpha}_\lambda$, then $\overline{\mu}\in  X^{m,\alpha}_\lambda$. Moreover, if $g\in (X^{m,\alpha}_\lambda)^\perp$, then $\overline{g}\in  (X^{m,\alpha}_\lambda)^\perp$. 
\end{lemma}
{\bf Proof.} If $\lambda$ is not real, then $X^{m,\alpha}_\lambda=\{0\}$ and both statements are true. Let $\lambda$ be real. If $\mu\in  X^{m,\alpha}_\lambda$, then there exists $v$ as in (\ref{mldc_eq:xmal}) such that $\mu=\frac{\partial v}{\partial\nu_\Omega}$.  Since $\lambda$ is real, then    $\Delta \overline{v}+\lambda\overline{v}=0$ in $\Omega$, $\overline{v}=0$ on $\partial\Omega$ and accordingly $\frac{\partial \overline{v}}{\partial\nu_\Omega}\in X^{m,\alpha}_\lambda$. Since
\[
\overline{\mu}=\overline{\frac{\partial v}{\partial\nu_\Omega}}
=\frac{\partial \overline{v}}{\partial\nu_\Omega}\,,
\]
we have $\overline{\mu}\in X^{m,\alpha}_\lambda$. Then by expoiting the definition of $(X^{m,\alpha}_\lambda)^\perp$, we see that if $g\in (X^{m,\alpha}_\lambda)^\perp$, then $\overline{g}\in  (X^{m,\alpha}_\lambda)^\perp$.\hfill  $\Box$ 

\vspace{\baselineskip}

\begin{proposition}\label{mldc_prop:1/2+w,vxmap}
 Let $m\in{\mathbb{N}}$, $\alpha\in]0,1[$. Let $\Omega$ be a bounded open subset of ${\mathbb{R}}^{n}$ of class $C^{\max\{1,m\},\alpha}$. Let $\lambda\in {\mathbb{C}}$. Let $S_{n,\lambda}$ is a fundamental solution of $\Delta+\lambda$. Then the following statements hold true. 
 \begin{enumerate}
\item[(i)] $ \frac{1}{2}\phi+W_\Omega[S_{n,\lambda},\phi]\in (X^{m,\alpha}_\lambda)^\perp$ for all $\phi\in  C^{m,\alpha}(\partial\Omega)$.
 \item[(ii)] $ V_\Omega[S_{n,\lambda},\psi] \in (X^{m,\alpha}_\lambda)^\perp$ for all $\psi\in C^{\max\{1,m\}-1,\alpha}(\partial\Omega)$.
\end{enumerate}
\end{proposition}
{\bf Proof.} (i) If $\phi\in C^{m,\alpha}(\partial\Omega)$ and $f\in X^{m,\alpha}_\lambda$, then 
\[
\int_{\partial\Omega}\left( \frac{1}{2}\phi+W_\Omega[S_{n,\lambda},\phi]\right)f\,d\sigma
=\int_{\partial\Omega}\phi \left( \frac{1}{2}f+W_\Omega^t[S_{n,\lambda},f]\right)\,d\sigma
\]
and Corollary \ref{mldc_corol:eipdd} (i) ensures that the right hand side equals $0$. 

(ii) If $\psi\in C^{\max\{1,m\}-1,\alpha}(\partial\Omega)$, $f\in X^{m,\alpha}_\lambda$, then there exists $v$ as in (\ref{mldc_eq:xmal}) such that $f=\frac{\partial v}{\partial\nu_\Omega}$ and thus a classical regularity result
for the  acoustic  single layer potential  implies that $v_\Omega^+[S_{n,\lambda},\psi]\in C^{1,\alpha}(\overline{\Omega})$ (cf. \textit{e.g.}, \cite[Thm.~7.1]{DoLa17}) and the second Green Identity implies that
\[
\int_{\partial\Omega}v_\Omega^+[S_{n,\lambda},\psi]f\,d\sigma=
\int_{\partial\Omega}v_\Omega^+[S_{n,\lambda},\psi]\frac{\partial v}{\partial\nu_\Omega}\,d\sigma=
\int_{\partial\Omega}\frac{\partial }{\partial\nu_\Omega}v_\Omega^+[S_{n,\lambda},\psi]v\,d\sigma=0\,,
\]
(see Remark \ref{mldc_rem:helpo} and 
\cite[Thm.~4.3]{DaLaMu21}).
\hfill  $\Box$ 

\vspace{\baselineskip}

\begin{proposition}\label{mldc_prop:silamokeri-wt}
Let $m\in{\mathbb{N}}$, $\alpha\in]0,1[$. Let $\Omega$ be a bounded open subset of ${\mathbb{R}}^{n}$ of class $C^{\max\{1,m\},\alpha}$.  
 Then the following statements hold.
\begin{enumerate}
\item[(i)] Let $\lambda\in {\mathbb{C}}$. Let $S_{n,\lambda}$ be a fundamental solution of $\Delta+\lambda$. If 
\[
\phi\in {\mathrm{Ker}}\left(
\frac{1}{2}I+W_\Omega^t[S_{n,\lambda},\cdot] 
\right)_{|C^{\max\{1,m\}-1,\alpha}(\partial\Omega)} 
\]
 (cf.~Remark \ref{mldc_rem:wtnotation}), then $v^-\equiv v_\Omega^-[S_{n,\lambda},\phi]$ belongs to  $C^{\max\{1,m\},\alpha}_{{\mathrm{loc}} }(\overline{\Omega^-})$ and solves the Neumann problem
 \begin{equation}\label{mldc_prop:silamokeri-wt0}
 \Delta v^-+\lambda v^-=0\qquad\text{in}\ \Omega^-\,,
 \qquad\frac{\partial v^-}{\partial\nu_{\Omega^-}}=0\qquad\text{on}\ \partial\Omega\,.
 \end{equation}
\item[(ii)] Let $\lambda\in {\mathbb{C}}$. Let $S_{n,\lambda}$ be a fundamental solution of $\Delta+\lambda$. If 
\[
\phi\in {\mathrm{Ker}}\left(
\frac{1}{2}I+W_\Omega^t[S_{n,\lambda},\cdot] \right)_{|C^{\max\{1,m\}-1,\alpha}(\partial\Omega)} 
\]
(cf.~Remark \ref{mldc_rem:wtnotation})   and $v_\Omega^-[S_{n,\lambda},\phi]_{|\partial\Omega}=0$ on $\partial\Omega$, then 
$\phi\in X^{m,\alpha}_{\lambda}$.
\item[(iii)]  Let $k\in {\mathbb{C}}\setminus ]-\infty,0]$, ${\mathrm{Im}}\,k\geq 0$. If 
\[
\phi\in {\mathrm{Ker}}\left(
\frac{1}{2}I+W_\Omega^t[\tilde{S}_{n,k;r},\cdot]
\right)_{|C^{\max\{1,m\}-1,\alpha}(\partial\Omega)} 
\]
(cf.~Remark \ref{mldc_rem:wtnotation}),  and $\phi\in X^{m,\alpha}_{k^2}$, then $v_\Omega^-[\tilde{S}_{n,k;r},\phi]=0$ in $\overline{\Omega^-}$.
\end{enumerate}
\end{proposition} 
{\bf Proof.} (i)   Classical regularity results on the  acoustic  single layer potential and the classical jump formulas for the normal derivative of the  acoustic  single layer potential imply that 
$v_\Omega^-[S_{n,\lambda},\phi]$ belongs to $ C^{\max\{1,m\},\alpha}_{{\mathrm{loc}} }(\overline{\Omega^-}) $ and satisfies the exterior Neumann problem (\ref{mldc_prop:silamokeri-wt0}) 
(cf.~\textit{e.g.}, \cite[Thm.~7.1]{DoLa17}). (ii) By classical regularity results on the  acoustic  single layer potential,   
$v_\Omega^+[S_{n,\lambda},\phi]$ belongs to $ C^{\max\{1,m\},\alpha}(\overline{\Omega}) $ 
(cf.~\textit{e.g.}, \cite[Thm.~7.1]{DoLa17}). 
Also, the classical jump conditions for the  acoustic  single layer potential imply that
\[v_\Omega^+[S_{n,\lambda},\phi]_{|\partial\Omega}=v_\Omega^-[S_{n,\lambda},\phi]_{|\partial\Omega}=0\,,
\]
 (cf.~\textit{e.g.}, \cite[Lem.~4.2 (i), Lem.~6.2]{DoLa17}). Then
$v_\Omega^+[S_{n,\lambda},\phi]$ satisfies the Dirichlet problem
\[
\Delta v_\Omega^+[S_{n,\lambda},\phi]+\lambda v_\Omega^+[S_{n,\lambda},\phi]=0, v_\Omega^+[S_{n,\lambda},\phi]_{|\partial\Omega}=0\,,
\]
and accordingly the definition (\ref{mldc_eq:xmal}) of $X^{m,\alpha}_{\lambda}$ 
  implies that  
$
\frac{\partial}{\partial\nu_{\Omega}}v_\Omega^+[S_{n,\lambda},\phi]\in X^{m,\alpha}_{\lambda}$. On the other hand, statement (i) implies that $\frac{\partial}{\partial\nu_{\Omega^-}}v_\Omega^-[S_{n,\lambda},\phi]=0$ on $\partial\Omega$. Then the classical jump conditions for the normal derivative of the  acoustic  single layer potential imply that
\[
\phi=-\frac{\partial}{\partial\nu_{\Omega^-}}v_\Omega^-[S_{n,\lambda},\phi]
-\frac{\partial}{\partial\nu_{\Omega}}v_\Omega^+[S_{n,\lambda},\phi]
=-\frac{\partial}{\partial\nu_{\Omega}}v_\Omega^+[S_{n,\lambda},\phi]\in X^{m,\alpha}_{\lambda} 
\]
(cf.~\textit{e.g.}, \cite[Thm.~7.1]{DoLa17}).  
(iii) Since 
\[
\phi\in {\mathrm{Ker}}\left(
\frac{1}{2}I+W_\Omega^t[\tilde{S}_{n,k;r},\cdot] 
\right)_{|C^{\max\{1,m\}-1,\alpha}(\partial\Omega)}\,,
\]
classical regularity results on the  acoustic  single layer potential and the classical jump formulas for the normal derivative of the  acoustic  single layer potential imply that 
$v_\Omega^-[\tilde{S}_{n,k;r},\phi]$ belongs to $ C^{\max\{1,m\},\alpha}_{{\mathrm{loc}} }(\overline{\Omega^-}) $  and satisfies the exterior Neumann problem 
\[
\Delta v_\Omega^-[\tilde{S}_{n,k;r},\phi]+k^2 v_\Omega^-[\tilde{S}_{n,k;r},\phi]=0\qquad\text{in}\ \Omega^-\,,
 \quad\frac{\partial v_\Omega^-[\tilde{S}_{n,k;r},\phi]}{\partial\nu_{\Omega^-}}=0\qquad\text{on}\ \partial\Omega\,.
\]
 Then it is also known that  $v_\Omega^-\left[\tilde{S}_{n,k;r},\phi\right]$ satisfies 
the outgoing $k$-radiation condition (cf.~\textit{e.g.}, \cite[Thm.~6.20 (i)]{La25b}).   By the uniqueness Theorem for the Neumann problem on the unbounded connected component of an exterior domain, we know that $v_\Omega^-\left[\tilde{S}_{n,k;r},\phi\right]$ equals zero on the unbounded connected component $(\Omega^-)_0$ of the exterior $\Omega^-$ (cf.~\cite[Thm.~1.3]{La25}). Thus 
by the Holmgren Uniqueness  Corollary \ref{mldc_corol:holmgr}   of the Appendix on the bounded connected components of the exterior $\Omega^-$,  it suffices to show that $v_\Omega^-\left[\tilde{S}_{n,k;r},\phi\right]$ equals $0$ on $\partial\Omega$.
Since  $\phi\in X^{m,\alpha}_{k^2}$, then there exists $v\in C^{\max\{1,m\},\alpha}(\overline{\Omega})$ such that
\[
\Delta v+k^2 v=0 \ \text{in}\ \Omega\,,\  v=0\ \text{on}\ \partial\Omega\,,
\]
and $\phi=\frac{\partial v}{\partial\nu_{\Omega}}$. 
Then the classical third Green Identity implies that 
\[
v=w_\Omega^+[\tilde{S}_{n,k;r},v_{|\partial\Omega}]-v_\Omega^+\left[\tilde{S}_{n,k;r},\frac{\partial v}{\partial\nu_{\Omega}}\right]=-v_\Omega^+\left[\tilde{S}_{n,k;r},\frac{\partial v}{\partial\nu_{\Omega}}\right]\quad\text{in}\ \Omega
\]
(cf.~\textit{e.g.}, \cite[Thm.~9.2 (ii)]{La25a}) and the classical jump conditions for the   acoustic  single layer potential imply that
\[
v_\Omega^-\left[\tilde{S}_{n,k;r},\phi\right]_{|\partial\Omega}=
v_\Omega^-\left[\tilde{S}_{n,k;r},\frac{\partial v}{\partial\nu_{\Omega}}\right]_{|\partial\Omega}=
v_\Omega^+\left[\tilde{S}_{n,k;r},\frac{\partial v}{\partial\nu_{\Omega}}\right]_{|\partial\Omega}
=-v_{|\partial\Omega}=0
\]
(cf.~\textit{e.g.}, \cite[Lem.~4.2 (i), Lem.~6.2]{DoLa17}).   Hence, the proof is complete. \hfill  $\Box$ 

\vspace{\baselineskip}

\begin{proposition}\label{mldc_prop:isomquotd}
 Let $m\in{\mathbb{N}}$, $\alpha\in]0,1[$. Let $\Omega$ be a bounded open subset of ${\mathbb{R}}^{n}$ of class $C^{\max\{1,m\},\alpha}$.  
 Let $\lambda\in {\mathbb{C}}$. Let $S_{n,\lambda}$ be a fundamental solution of $\Delta+\lambda$. Then the following statements hold true (cf.~Remark \ref{mldc_rem:wtnotation}). 
\begin{enumerate}
\item[(i)] ${\mathrm{Im}}\, \left(
\frac{1}{2}I+W_\Omega[S_{n,\lambda},\cdot]\right)_{|C^{m,\alpha}(\partial\Omega)}$ is a closed subspace of
$C^{m,\alpha}(\partial\Omega)$ and 
the quotient space $\frac{
C^{m,\alpha}(\partial\Omega)
}{
{\mathrm{Im}}\, \left(
\frac{1}{2}I+W_\Omega[S_{n,\lambda},\cdot]
\right)_{|C^{m,\alpha}(\partial\Omega)}}$ has dimension equal to the finite dimension of  ${\mathrm{Ker}}\, \left(
\frac{1}{2}I+W_\Omega^t[S_{n,\lambda},\cdot]\right)_{|C^{\max\{1,m\}-1,\alpha}(\partial\Omega)}$. 
\item[(ii)] ${\mathrm{Im}}\, \left(
\frac{1}{2}I+W_\Omega[S_{n,\lambda},\cdot]\right)_{|C^{m,\alpha}(\partial\Omega)}$ is a closed subspace of
$(X^{m,\alpha}_\lambda)^\perp$. 
\item[(iii)] The quotient space $\frac{
C^{m,\alpha}(\partial\Omega)
}{
(X^{m,\alpha}_\lambda)^\perp
}$ is linearly homeomorphic to $X^{m,\alpha}_\lambda$ and has finite dimension.
\item[(iv)] The following formula holds
\begin{eqnarray*}
\lefteqn{
{\mathrm{dim}}\,\left(
\frac{(X^{m,\alpha}_\lambda)^\perp}{
{\mathrm{Im}}\, \left(
\frac{1}{2}I+W_\Omega[S_{n,\lambda},\cdot] \right)_{|C^{m,\alpha}(\partial\Omega)}
}\right)
}
\\ \nonumber
&&\qquad
={\mathrm{dim}}\,{\mathrm{Ker}}\, \left(
\frac{1}{2}I+W_\Omega^t[S_{n,\lambda},\cdot]\right)_{|C^{\max\{1,m\}-1,\alpha}(\partial\Omega)}
-{\mathrm{dim}}\,X^{m,\alpha}_\lambda\,.
\end{eqnarray*}
\end{enumerate}
\end{proposition}
{\bf Proof.} (i) Since $W_\Omega[S_{n,\lambda},\cdot]_{|C^{m,\alpha}(\partial\Omega)} $ is compact in $C^{m,\alpha}(\partial\Omega)$ (cf.~\textit{e.g.}, \cite[Thm.~7.4]{DoLa17} in case $m=0$ and  \cite[Cor.~9.1]{DoLa17} in case $m\geq 1$) and
$W_\Omega^t[S_{n,\lambda},\cdot]_{|C^{\max\{1,m\}-1,\alpha}(\partial\Omega)} $ is compact in $C^{\max\{1,m\}-1,\alpha}(\partial\Omega)$ (cf.~\textit{e.g.}, \cite[Cor.~10.1]{DoLa17}), then the Fredholm Alternative Theorem of Wendland \cite{We67}, \cite{We70} in the duality pairing
\[
\left(C^{\max\{1,m\}-1,\alpha} (\partial\Omega),C^{m,\alpha} (\partial\Omega)\right) 
\]
(cf.~\textit{e.g.}, \cite[Thm.~4.17]{Kr14}) implies the validity of the following equality
\begin{eqnarray*}
\lefteqn{
{\mathrm{Im}}\, \left(
\frac{1}{2}I+W_\Omega[S_{n,\lambda},\cdot] \right)_{|C^{m,\alpha}(\partial\Omega)}
}
\\ \nonumber
&&\qquad\qquad\qquad\qquad
=\left(
{\mathrm{Ker}}\, \left(
\frac{1}{2}I+W_\Omega^t[S_{n,\lambda},\cdot] \right)_{|C^{\max\{1,m\}-1,\alpha}(\partial\Omega)}
\right)^\perp\,,
\end{eqnarray*}
where the orthogonality has to be understood according to the above duality pairing
and that the finite dimension $s$ of ${\mathrm{Ker}}\, \left(
\frac{1}{2}I+W_\Omega^t[S_{n,\lambda},\cdot] \right)_{|C^{\max\{1,m\}-1,\alpha}(\partial\Omega)}$
equals that of ${\mathrm{Ker}}\, \left(
\frac{1}{2}I+W_\Omega[S_{n,\lambda},\cdot] \right)_{|C^{m,\alpha}(\partial\Omega)}$. If $s=0$, then 
\[
\left(
{\mathrm{Ker}}\, \left(
\frac{1}{2}I+W_\Omega^t[S_{n,\lambda},\cdot] \right)_{|C^{\max\{1,m\}-1,\alpha}(\partial\Omega)}
\right)^\perp=C^{m,\alpha}(\partial\Omega)
\]
and statement (i) is obviously true. Thus we can assume that $s\geq 1$. 
Let $\{\mu_j\}_{j=1}^s$ be a basis of ${\mathrm{Ker}}\, \left(
\frac{1}{2}I+W_\Omega^t[S_{n,\lambda},\cdot] \right)_{|C^{\max\{1,m\}-1,\alpha}(\partial\Omega)}$. 

Then we consider the map $\Lambda$ from $C^{m,\alpha}(\partial\Omega)$ to ${\mathbb{C}}^s$ that is defined by the equality
\[
\Lambda[g]\equiv\left(\int_{\partial\Omega}g\mu_1\,d\sigma, 
\dots,
\int_{\partial\Omega}g\mu_s\,d\sigma
\right)\qquad\forall g\in C^{m,\alpha}(\partial\Omega)\,.
\]
If $\Lambda$ is not surjective, then there exists $(c_1,\dots,c_s)\in  {\mathbb{C}}^s$ such that
\[
0=\sum_{j=1}^s c_j\int_{\partial\Omega}g\mu_j\,d\sigma=\int_{\partial\Omega}g\sum_{j=1}^sc_j\mu_j\,d\sigma
\qquad\forall g\in C^{m,\alpha}(\partial\Omega)\,.
\]
Hence, $\sum_{j=1}^sc_j\mu_j=0$, a contradiction. Then $\Lambda$ is surjective. Since
\begin{eqnarray*}
\lefteqn{
{\mathrm{Ker}}\, \Lambda=\left(
{\mathrm{Ker}}\, \left(
\frac{1}{2}I+W_\Omega^t[S_{n,\lambda},\cdot] \right)_{|C^{\max\{1,m\}-1,\alpha}(\partial\Omega)}
\right)^\perp
}
\\ \nonumber
&&\qquad
=
{\mathrm{Im}}\, \left(
\frac{1}{2}I+W_\Omega[S_{n,\lambda},\cdot] \right)_{|C^{m,\alpha}(\partial\Omega)}\,,
\end{eqnarray*}
 the Homomorphism Theorem implies that 
\[
{\mathbb{C}}^s \quad \text{is\ isomorphic\ to}\quad  
\frac{
C^{m,\alpha}(\partial\Omega)
}{
{\mathrm{Im}}\, \left(
\frac{1}{2}I+W_\Omega[S_{n,\lambda},\cdot]
\right)_{|C^{m,\alpha}(\partial\Omega)}
}
\]
and accordingly, statement (i) holds true. Statement (ii) holds by statement (i) and by  Proposition \ref{mldc_prop:1/2+w,vxmap} (i).
Since $X^{m,\alpha}_\lambda$ can be identified with a linear subspace of the dual of 
$C^{m,\alpha}(\partial\Omega)$, it is known that the dual
\[
\left(\frac{
C^{m,\alpha}(\partial\Omega)
}{
(X^{m,\alpha}_\lambda)^\perp
}\right)'\quad\text{is\ isomorphic\ to}\ (X^{m,\alpha}_\lambda)^{\perp\perp}
\]
(cf.~\textit{e.g.}, Rudin, \cite[Thm.~4.9, Chap. IV]{Ru91}). On the other hand, $X^{m,\alpha}_\lambda$ is finite dimensional and accordingly, $X^{m,\alpha}_\lambda$ is weak$^\ast$ closed and
\[
(X^{m,\alpha}_\lambda)^{\perp\perp}=X^{m,\alpha}_\lambda
\]
(cf.~\textit{e.g.}, Rudin, \cite[Thm.~4.7, Chap. IV]{Ru91}). Hence, $\left(\frac{
C^{m,\alpha}(\partial\Omega)
}{
(X^{m,\alpha}_\lambda)^\perp
}\right)'$ is linearly  isomorphic to $X^{m,\alpha}_\lambda$ and has finite dimension. Hence, 
$\frac{
C^{m,\alpha}(\partial\Omega)
}{
(X^{m,\alpha}_\lambda)^\perp
}$ is also linearly isomorphic to $X^{m,\alpha}_\lambda$. On the other hand linear isomorphisms between normed spaces of finite dimension are also homeomorphisms and thus statement (iii) holds true. 

 By the Homomorphism Theorem (cf.~\textit{e.g.}, Greub \cite[7. b), p.~49]{Gr75}), the quotient $\left(
 \frac{C^{m,\alpha}(\partial\Omega)}{(X^{m,\alpha}_\lambda)^\perp}
  \right)$ is linearly isomorphic to the quotient
\[
\left(\frac{
C^{m,\alpha}(\partial\Omega)
}{
{\mathrm{Im}}\, \left(
\frac{1}{2}I+W_\Omega[S_{n,\lambda},\cdot]
\right)_{|C^{m,\alpha}(\partial\Omega)}
}\right)
{\displaystyle{\diagup}}
\left(
\frac{(X^{m,\alpha}_\lambda)^\perp}
{{\mathrm{Im}}\, \left(
\frac{1}{2}I+W_\Omega[S_{n,\lambda},\cdot]
\right)_{|C^{m,\alpha}(\partial\Omega)}}
\right)
\]
that has dimension equal to 
\begin{eqnarray*}
\lefteqn{
{\mathrm{dim}}\,{\mathrm{Ker}}\, \left(
\frac{1}{2}I+W_\Omega^t[S_{n,\lambda},\cdot]
\right)_{|C^{\max\{1,m\}-1,\alpha}(\partial\Omega)}
}
\\ \nonumber
&&\qquad
-{\mathrm{dim}}\,\left(
\frac{(X^{m,\alpha}_\lambda)^\perp}
{{\mathrm{Im}}\, \left(
\frac{1}{2}I+W_\Omega[S_{n,\lambda},\cdot]
\right)_{|C^{m,\alpha}(\partial\Omega)}}
\right)
\end{eqnarray*}
by statement (i). Hence, statement (iii) implies the validity of the formula in statement (iv) and the proof is complete.\hfill  $\Box$ 

\vspace{\baselineskip}

  Then we can prove the following statement that says that the space of compatible data $(X^{m,\alpha}_\lambda)^\perp$ can be written as an algebraic sum of data that equal the boundary values of an  acoustic  double   layer potential and of an  acoustic  single  layer potential. 
\begin{theorem}\label{mldc_thm:codasumwv}
Let $m\in{\mathbb{N}}$, $\alpha\in]0,1[$. Let $\Omega$ be a bounded open subset of ${\mathbb{R}}^{n}$ of class $C^{\max\{1,m\},\alpha}$.  
 Let $\lambda\in {\mathbb{C}}$. Let $S_{n,\lambda}$ be a fundamental solution of $\Delta+\lambda$. Then
\begin{equation}\label{mldc_thm:codasumwv1}
(X^{m,\alpha}_\lambda)^\perp={\mathrm{Im}}\, \left(
\frac{1}{2}I+W_\Omega[S_{n,\lambda},\cdot]
\right)_{|C^{m,\alpha}(\partial\Omega)}
+
V_\Omega[S_{n,\lambda},C^{\max\{1,m\}-1,\alpha}(\partial\Omega)]\,.
\end{equation}
\end{theorem}
{\bf Proof.} By Proposition \ref{mldc_prop:isomquotd} (i), (ii), (iii),  ${\mathrm{Im}}\, \left(
\frac{1}{2}I+W_\Omega[S_{n,\lambda},\cdot]
\right)_{|C^{m,\alpha}(\partial\Omega)}$ is a closed subspace of
$(X^{m,\alpha}_\lambda)^\perp$ of finite codimension. Let $\pi$ be the canonical projection of 
$(X^{m,\alpha}_\lambda)^\perp$ onto the finite dimensional quotient space 
\[
Q^{m,\alpha}_\lambda\equiv 
\frac{(X^{m,\alpha}_\lambda)^\perp}{
{\mathrm{Im}}\, \left(
\frac{1}{2}I+W_\Omega[S_{n,\lambda},\cdot]
\right)_{|C^{m,\alpha}(\partial\Omega)}
} \,.
\]
 By Proposition \ref{mldc_prop:1/2+w,vxmap}  (ii), the image $V_\Omega[S_{n,\lambda},C^{\max\{1,m\}-1,\alpha}(\partial\Omega)]$ is a subspace of $(X^{m,\alpha}_\lambda)^\perp$ and thus  the $\pi$-image $\pi\left[Y^{m,\alpha}_\lambda\right]$
of
\[
Y^{m,\alpha}_\lambda\equiv  {\mathrm{Im}}\, \left(
\frac{1}{2}I+W_\Omega[S_{n,\lambda},\cdot]
\right)_{|C^{m,\alpha}(\partial\Omega)}
+
V_\Omega[S_{n,\lambda},C^{\max\{1,m\}-1,\alpha}(\partial\Omega)] 
\]
is a   subspace of the finite dimensional Banach quotient space $ 
Q^{m,\alpha}_\lambda$ and is accordingly closed in $Q^{m,\alpha}_\lambda$. Since the canonical projection $\pi$ is continuous, then the preimage
\[
\pi^{\leftarrow}\left[\pi\left[Y^{m,\alpha}_\lambda \right]\right]
\]
is closed in $(X^{m,\alpha}_\lambda)^\perp$. Now a simple computation based on the definition of $\pi$  shows that
\[
\pi^{\leftarrow}\left[\pi\left[Y^{m,\alpha}_\lambda \right]\right]=Y^{m,\alpha}_\lambda
\]
and accordingly $Y^{m,\alpha}_\lambda$ is a closed subspace of $(X^{m,\alpha}_\lambda)^\perp$. Assume by contradiction that there exists
\[
\tau\in (X^{m,\alpha}_\lambda)^\perp\setminus Y^{m,\alpha}_\lambda\,.
\]
Then the Hahn-Banach Theorem implies the existence of   $\tilde{\mu}\in \left(C^{m,\alpha}(\partial\Omega)\right)'$ such that
\begin{equation}\label{mldc_thm:codasumwv2}
\langle \tilde{\mu},\tau\rangle=1\,,\qquad \langle \tilde{\mu},\eta\rangle=0\qquad
\forall \eta \in Y^{m,\alpha}_\lambda\,.
\end{equation}
Since $\tilde{\mu}$ annihilates ${\mathrm{Im}}\, \left(
\frac{1}{2}I+W_\Omega[S_{n,\lambda},\cdot]\right)_{|C^{m,\alpha}(\partial\Omega)}$, then
\[
\tilde{\mu}\in {\mathrm{Ker}}\, \left(
\frac{1}{2}I+W_\Omega^t[S_{n,\lambda},\cdot]\right)\,,
\]
where $W_\Omega^t[S_{n,\lambda},\cdot]$ is the transpose operator of $W_\Omega[S_{n,\lambda},\cdot]$ with respect to the canonical duality pairing
\[
\left(\left(C^{m,\alpha}(\partial\Omega)\right)',C^{m,\alpha}(\partial\Omega)\right)\,.
\]
Since $W_\Omega[S_{n,\lambda},\cdot]$ is compact in $C^{m,\alpha}(\partial\Omega)$ (cf.~\textit{e.g.}, \cite[Thm.~7.4]{DoLa17} in case $m=0$ and  \cite[Cor.~9.1]{DoLa17} in case $m\geq 1$), the Fredholm Alternative Theorem implies that
\[
{\mathrm{dim}}\,{\mathrm{Ker}}\, \left(
\frac{1}{2}I+W_\Omega^t[S_{n,\lambda},\cdot]\right)
={\mathrm{dim}}\,{\mathrm{Ker}}\,
\left(
\frac{1}{2}I+W_\Omega[S_{n,\lambda},\cdot]
\right)_{|C^{m,\alpha}(\partial\Omega)}
 \,.
\]
On the other hand, also the operator  $W_\Omega^t[S_{n,\lambda},\cdot]_{|C^{\max\{1,m\}-1,\alpha}(\partial\Omega)} $ is compact in $C^{\max\{1,m\}-1,\alpha}(\partial\Omega)$ (cf.~\textit{e.g.}, \cite[Cor.~10.1]{DoLa17}) and thus the Fredholm Alternative Theorem of Wendland \cite{We67}, \cite{We70} in the duality pairing
\begin{eqnarray}\label{mldc_thm:codasumwv3}
\left(C^{\max\{1,m\}-1,\alpha} (\partial\Omega),C^{m,\alpha} (\partial\Omega)\right) 
\end{eqnarray}
(cf.~\textit{e.g.}, Kress~\cite[Thm.~4.17]{Kr14}) implies that
\begin{eqnarray*}
\lefteqn{
{\mathrm{dim}}\,{\mathrm{Ker}}\, \left(
\frac{1}{2}I+W_\Omega^t[S_{n,\lambda},\cdot]
\right)_{|C^{\max\{1,m\}-1,\alpha}(\partial\Omega)}
}
\\ \nonumber
&&\qquad\qquad\qquad
={\mathrm{dim}}\,{\mathrm{Ker}}\,
\left(
\frac{1}{2}I+W_\Omega[S_{n,\lambda},\cdot]
\right)_{|C^{m,\alpha}(\partial\Omega)}
 \,,
\end{eqnarray*}
where now $W_\Omega^t[S_{n,\lambda},\cdot]_{|C^{\max\{1,m\}-1,\alpha}(\partial\Omega)}$ is the transpose with respect to the duality pairing in (\ref{mldc_thm:codasumwv3}). Since the canonical embedding ${\mathcal{J}}$ of Lemma \ref{mldc_lem:caincl} maps 
\begin{equation}\label{mldc_thm:codasumwv4}
{\mathrm{Ker}}\,\left(
\frac{1}{2}I+W_\Omega^t[S_{n,\lambda},\cdot]
\right)_{|C^{\max\{1,m\}-1,\alpha}(\partial\Omega)}
\quad\text{into}\quad {\mathrm{Ker}}\, \left(
\frac{1}{2}I+W_\Omega^t[S_{n,\lambda},\cdot]\right)\,,
\end{equation}
the above finite dimensional equalities imply that ${\mathcal{J}}$ restricts a bijection between the spaces in (\ref{mldc_thm:codasumwv4}). Hence, there exists 
\[
\mu\in {\mathrm{Ker}}\, \left(
\frac{1}{2}I+W_\Omega^t[S_{n,\lambda},\cdot]
\right)_{|C^{\max\{1,m\}-1,\alpha}(\partial\Omega)}
\]
 such that $\tilde{\mu}={\mathcal{J}}[\mu]$. Since $\tilde{\mu}$ annihilates $Y^{m,\alpha}_\lambda$, then $\tilde{\mu}$ annihilates its subspace $V_\Omega[S_{n,\lambda},C^{\max\{1,m\}-1,\alpha}(\partial\Omega)] $ and accordingly,
\begin{eqnarray*}
\lefteqn{
0=\langle \tilde{\mu},V_\Omega[S_{n,\lambda},\eta]\rangle=\langle {\mathcal{J}}[\mu], V_\Omega[S_{n,\lambda},\eta]\rangle
}
\\ \nonumber
&& 
=\int_{\partial\Omega}\mu  V_\Omega[S_{n,\lambda},\eta]\,d\sigma
=\int_{\partial\Omega} V_\Omega[S_{n,\lambda},\mu]\eta\,d\sigma
\qquad\forall\eta
\in C^{\max\{1,m\}-1,\alpha}(\partial\Omega)\,.
\end{eqnarray*}
Hence, $V_\Omega[S_{n,\lambda},\mu]=0$. Then Proposition \ref{mldc_prop:silamokeri-wt} (ii) implies that $\mu\in X^{m,\alpha}_\lambda$. Since $\tau\in (X^{m,\alpha}_\lambda)^\perp$, we must have 
\[
\langle \tilde{\mu},\tau\rangle=\int_{\partial\Omega}\mu\tau\,d\sigma=0\,,
\]
a contradiction (see (\ref{mldc_thm:codasumwv2})). Hence, $(X^{m,\alpha}_\lambda)^\perp=Y^{m,\alpha}_\lambda$ and the proof is complete.\hfill  $\Box$ 

\vspace{\baselineskip}

Then by combining the necessary condition of Theorem \ref{mldc_thm:nexintdirhe}, Theorem \ref{mldc_thm:codasumwv} and the jump formulas for the single and the  acoustic  double layer potential, we deduce the validity of the following statement.
\begin{theorem}[of existence for the interior Dirichlet problem]\label{mldc_thm:exintdirhe}
  Let $m\in {\mathbb{N}}$,   $\alpha\in]0,1[$. Let $\Omega$ be a bounded open subset of ${\mathbb{R}}^{n}$ of class $C^{\max\{1,m\},\alpha}$. Let  $\lambda\in {\mathbb{C}}$. Let $S_{n,\lambda}$ be a fundamental solution of $\Delta+\lambda$.  Let $g\in C^{m,\alpha}(\partial\Omega)$. Then the interior Dirichlet problem  
\begin{equation}
\label{mldc_thm:exintdirhe1}
\left\{
\begin{array}{ll}
\Delta u+\lambda u=0 &{\text{in}}\ \Omega\,,
\\
u=g
&{\text{on}}\ \partial\Omega\,
\end{array}
\right.
\end{equation}
has a solution $u\in C^{m,\alpha}(\overline{\Omega})$ if and only if the compatibility condition
\begin{equation}\label{mldc_thm:exintdirhe2}
\int_{\partial\Omega}g\frac{\partial v}{\partial\nu_\Omega}\,d\sigma=0  \,,
\end{equation}
is satisfied for all $v\in C^{\max\{1,m\},\alpha}(\overline{\Omega})$  such that
\begin{equation}\label{mldc_thm:exintdirhe3}
\Delta v+\lambda v=0 \qquad\text{in}\ \Omega\,,\qquad v=0\qquad\text{on}\ \partial\Omega\,.
\end{equation}
\end{theorem}
{\bf Proof.} The necessity follows by Theorem \ref{mldc_thm:nexintdirhe}. If $g$ satisfies condition (\ref{mldc_thm:exintdirhe2}), then $g\in (X^{m,\alpha}_{\lambda})^\perp$ and Theorem \ref{mldc_thm:codasumwv}
implies the existence of $\phi\in C^{m,\alpha}(\partial\Omega)$ and of $\psi\in C^{\max\{1,m\}-1,\alpha}(\partial\Omega)$ such that
\[
g=\frac{1}{2}\phi+W_\Omega[S_{n,\lambda},\phi]+V_\Omega[S_{n,\lambda},\psi]\,.
\]
Thus if we set
\[
u=w_\Omega^+[S_{n,\lambda},\phi] +v_\Omega^+[S_{n,\lambda},\psi]\,,
\]
then $u\in  C^{m,\alpha}(\overline{\Omega})$ by known  regularity results for the  acoustic  single   layer potential (cf.~\textit{e.g.}, \cite[Thm.~7.1]{DoLa17}) and by known regularity results for the  acoustic  double layer   potential (cf.~\textit{e.g.}, 
 \cite[Thm.~7.3]{DoLa17} in case $m\geq 1$, 
 \cite[Thm.~5.1]{La25a} in case $m=0$). Also, the known  jump formulas for the acoustic single and double  layer potential  (cf.~\textit{e.g.},  \cite[Thms.~7.1, 7.3]{DoLa17}) imply that $u_{|\partial\Omega}=g$.\hfill  $\Box$ 

\vspace{\baselineskip}
 
\section{A representation theorem for the exterior Diri\-chlet eigenfunctions}\label{mldc_sec:rediefh}

We first   show the following representation theorem for the exterior Dirichlet eigenfunctions that is classical for $m\geq 1$. Here we note that the normal derivative is to be interpreted in the sense of Definition \ref{mldc_defn:endedr}  if  $m=0$.
\begin{theorem}
\label{mldc_thm:eiped}
Let $m\in{\mathbb{N}}$, $\alpha\in]0,1[$. Let $\Omega$ be an open subset of ${\mathbb{R}}^{n}$ of class $C^{\max\{1,m\},\alpha}$. Then the following statements hold.
\begin{enumerate}
\item[(i)] Let  $k\in {\mathbb{C}}\setminus ]-\infty,0]$, ${\mathrm{Im}}\,k\geq 0$. 
  If   $u\in C^{m,\alpha}_{{\mathrm{loc}} }(\overline{\Omega^-}) $ satisfies the problem
\begin{equation}
\label{mldc_thm:eiped1}
\left\{
\begin{array}{ll}
\Delta u+k^2 u=0 &{\text{in}}\ \Omega^-\,,
\\
u=0 
&{\text{on}}\ \partial\Omega\,,
\\
u \ \text{satisfies\ the\ outgoing}\   \text{$k$-radiation\ condition}\,,
\end{array}
\right.
\end{equation}
then  $u\in C^{\max\{1,m\},\alpha}_{{\mathrm{loc}} }(\overline{\Omega^-})$,
\begin{equation}
\label{mldc_thm:eiped2}
u(x)= -v_\Omega^-[\tilde{S}_{n,k;r},\frac{\partial u}{\partial\nu_{\Omega^-}}](x)\qquad\forall x\in\overline{\Omega^-}\,,
\end{equation}
and 
\begin{equation}
\label{mldc_thm:eiped3}
-\frac{1}{2}\frac{\partial u}{\partial\nu_{\Omega^-}} 
+
 W_\Omega^t[\tilde{S}_{n,k;r},
\frac{\partial u}{\partial\nu_{\Omega^-}}]=0\qquad\text{on}\ \partial\Omega\,.
\end{equation}
\item[(ii)] Let $k\in {\mathbb{C}}\setminus ]-\infty,0]$, ${\mathrm{Im}}\,k\geq 0$. If $k^{2}$ is not a Neumann eigenvalue for $-\Delta$ in $\Omega$, and if $\mu\in V^{m-1,\alpha}(\partial\Omega)$ (cf.~(\ref{mldc_eq:vm-1a})) 
satisfies the equation
\begin{equation}
\label{mldc_thm:eiped4}
-\frac{1}{2}\mu 
+W_\Omega^t[\tilde{S}_{n,k;r},
\mu]=0\qquad\text{on}\ \partial\Omega\,,
\end{equation}
then the function $u$ from $\overline{\Omega^-}$ to ${\mathbb{C}}$ defined by 
\begin{equation}
\label{mldc_thm:eiped5}
u(x)= -v_\Omega^-[\tilde{S}_{n,k;r},\mu](x)\qquad\forall x\in\overline{\Omega^-}\,,
\end{equation}
 belongs to $C^{\max\{1,m\},\alpha}_{{\mathrm{loc}} }(\overline{\Omega^-}) $,  satisfies the Dirichlet problem (\ref{mldc_thm:eiped1})  
  and 
\begin{equation}
\label{mldc_thm:eiped6}
\frac{\partial u}{\partial\nu_{\Omega^-}}=\mu\qquad{\text{on}}\ \partial
\Omega\,.
\end{equation}
In particular,  $\mu\in C^{\max\{1,m\}-1,\alpha}(\overline{\Omega}) $. 
Moreover, $v_\Omega^+[\tilde{S}_{n,k;r},\mu]=0$ in $\Omega$. 
\end{enumerate}
\end{theorem}
{\bf Proof.} (i) The membership of $u$ in $C^{\max\{1,m\},\alpha}_{{\mathrm{loc}} }(\overline{\Omega^-}) $ holds by assumption in case $m\geq 1$ and is a consequence of Remark \ref{mldc_rem:helpo} and of the regularization theorem of  \cite[Thm~1.2 (iv)]{La25} in case $m=0$. 

Equality (\ref{mldc_thm:eiped2}) is an immediate consequence of the  third Green identity  for exterior domains   (cf.~\cite[
Thm.~8.24]{La25b}) and of the continuity in $\overline{\Omega^-}$ of the  acoustic   simple layer potential  (cf.~\textit{e.g.},  \cite[Thm.~7.1]{DoLa17}). By taking the $\nu_{\Omega^-}$-normal derivative of (\ref{mldc_thm:eiped2}) on $\partial\Omega$, we obtain
\[
\frac{\partial u}{\partial\nu_{\Omega^-}}(x)
=\frac{1}{2}\frac{\partial u}{\partial\nu_{\Omega^-}}(x)
+
W_\Omega^t[\tilde{S}_{n,k;r},\frac{\partial u}{\partial\nu_{\Omega^-}}](x)
\qquad\forall x\in\partial\Omega\,,
\]
(cf.~\textit{e.g.},  \cite[Thm.~7.1]{DoLa17}) and thus equality (\ref{mldc_thm:eiped3}) follows. We now prove statement (ii). We set $u^{\pm}\equiv  -v_\Omega^\pm[\tilde{S}_{n,k;r},\mu]$. If $m=0$, then \cite[Thm.~5.3]{La25a} implies that $u^+\in C^{0,\alpha}(\overline{\Omega})_\Delta $ and  $u^-\in C^{0,\alpha}_{{\mathrm{loc}} }(\overline{\Omega^-})_\Delta$. If $m\geq 1$, then it is classically known that $u^+\in C^{m,\alpha}(\overline{\Omega}) $ and $u^-\in C^{m,\alpha}_{{\mathrm{loc}} }(\overline{\Omega^-})$ (cf.~\textit{e.g.},  \cite[Thm.~7.1]{DoLa17}). 

By equation (\ref{mldc_thm:eiped4}) and by the jump formulas for the normal derivative of the acoustic simple layer potential of \cite[Thm.~6.3]{La25a}, we have 
 $\frac{\partial u^{+}}{\partial\nu_{\Omega}}=0$ on $\partial\Omega$. Then $u^+$ solves the Neumann problem for the Helmholtz equation in $\Omega$  and our assumption on $k^2$ implies that
   $u^{+}$ equals $0$ in $\Omega$.   Also, $u^{-}_{|\partial\Omega}=u^{+}_{|\partial\Omega}=0$ on $\partial\Omega$
   (cf.~\cite[Thm.~6.3]{La25a}). On the other hand, $u^{-}$ satisfies the outgoing $k$-radiation condition  (cf.~\cite[Thm.~6.20 (i)]{La25b})  and thus  $u^{-}$ satisfies the exterior Dirichlet problem (\ref{mldc_thm:eiped1}). Then statement (i) implies that $u\in C^{\max\{1,m\},\alpha}_{{\mathrm{loc}}}(\overline{\Omega^-})$. Then the jump formula for the normal derivative of  the  acoustic  single layer potential implies that 
   \[
\mu=\frac{\partial u^-}{\partial \nu_{\Omega^-}}  + \frac{\partial u^+}{\partial \nu_\Omega}= \frac{\partial u^-}{\partial \nu_{\Omega^-}}
   \]
   (cf. \textit{e.g.}, \cite[Thm.~6.3]{La25a}). In particular, equality (\ref{mldc_thm:eiped6}) holds true.\hfill  $\Box$ 

\vspace{\baselineskip}

\begin{corollary}
\label{mldc_corol:eipedd}
Let  $m\in{\mathbb{N}}$, $\alpha\in]0,1[$. Let $\Omega$ be a bounded open subset of ${\mathbb{R}}^{n}$ of class $C^{\max\{1,m\},\alpha}$.  Let  $k\in {\mathbb{C}}\setminus ]-\infty,0]$, ${\mathrm{Im}}\,k\geq 0$. 
\begin{enumerate}
\item[(i)] The  operator $\frac{\partial}{\partial\nu_{\Omega^-}}$ is a linear injection from the space 
\begin{eqnarray}
\label{mldc_corol:eipedd1}
\lefteqn{
\{
u\in C^{m,\alpha}_{{\mathrm{loc}}}(\overline{\Omega^-}):\,\Delta u+k^{2}u=0\ 
\text{in}\ \Omega^-, u_{|\partial\Omega}=0\,, 
 }
\\ \nonumber
&&\qquad
u \ \text{satisfies\ the\ outgoing\ $k$-radiation\ condition} \}
\\ \nonumber
&&\qquad
=\{
u\in C^{\max\{1,m\},\alpha}_{{\mathrm{loc}}}(\overline{\Omega^-}):\,\Delta u+k^{2}u=0\ 
\text{in}\ \Omega^-, u_{|\partial\Omega}=0\,, 
\\ \nonumber
&&\qquad\quad
u \ \text{satisfies\ the\ outgoing\ $k$-radiation\ condition} \}
\\ \nonumber
&&\qquad
=\{
u\in C^{\max\{1,m\},\alpha}_{{\mathrm{loc}}}(\overline{\Omega^-}):\,\Delta u+k^{2}u=0\ 
\text{in}\ \Omega^-, u_{|\partial\Omega}=0\,, 
\\ \nonumber
&&\qquad\quad
u(x)=0\ \forall x\in (\Omega^-)_0 \}
\end{eqnarray}
into the space
\begin{equation}\label{mldc_corol:eipedd2}
\{
\mu\in C^{\max\{1,m\}-1,\alpha}(\partial\Omega):\,(\ref{mldc_thm:eiped4}) \ holds\}\,,
\end{equation}
where $(\Omega^-)_0$ denotes the unbounded connected component of $\Omega^-$.
Moreover,
\begin{eqnarray}\label{mldc_corol:eipedd3}
\lefteqn{
\{
\mu\in C^{\max\{1,m\}-1,\alpha}(\partial\Omega):\,(\ref{mldc_thm:eiped4}) \ holds\}
}
\\ \nonumber
&&\qquad\qquad 
=
\left\{
 \mu\in V^{m-1,\alpha}(\partial\Omega):\,(\ref{mldc_thm:eiped4}) \ holds\right\}\,,
 \end{eqnarray}
(cf.~(\ref{mldc_eq:vm-1a})).

\item[(ii)] Assume that $k^{2}$ is not a Neumann eigenvalue for $-\Delta$ in $\Omega$. The operator $\frac{\partial}{\partial\nu_{\Omega^-}}$ induces an isomorphism from the space in (\ref{mldc_corol:eipedd1}) onto the space 
\begin{equation}\label{mldc_corol:eipedd4}
\{
\mu\in C^{\max\{1,m\}-1,\alpha}(\partial\Omega):\,(\ref{mldc_thm:eiped4}) \ holds\} 
\end{equation}
 and the inverse of such an isomorphism is the map $-v_\Omega^-[\tilde{S}_{n,k;r},\cdot]$. In particular, the dimension of the space in (\ref{mldc_corol:eipedd4}) equals the geometric multiplicity of $k^{2}$ as a Dirichlet eigenvalue of $-\Delta$ in $\bigcup_{j=1}^{\kappa^-}(\Omega^-)_j$
 (or of $-\Delta$ in $\Omega^-$, with the understanding that the eigenfunctions should satisfy the outgoing $k$-radiation condition).
\end{enumerate}
\end{corollary}
{\bf Proof.} (i) By Theorem \ref{mldc_thm:eiped} (i) and by the uniqueness theorem for the Dirichlet problem  (\ref{mldc_thm:eiped1}) on the unbounded connected component $(\Omega^-)_0$ of $\Omega^-$  (cf. \cite[Thm.~1.3 (ii)]{La25}), the equalities in  (\ref{mldc_corol:eipedd1}) hold true.
Then again Theorem \ref{mldc_thm:eiped} (i) implies that  $\frac{\partial u}{\partial\nu_{\Omega^-}}$ is a linear operator from the space in (\ref{mldc_corol:eipedd1}) to the space in (\ref{mldc_corol:eipedd2}). By the uniqueness theorem for the Dirichlet problem  (\ref{mldc_thm:eiped1}) on the unbounded connected component $(\Omega^-)_0$ of $\Omega^-$ (cf. \cite[Thm.~1.3]{La25}) and by the Holmgren  Uniqueness  Corollary \ref{mldc_corol:holmgr}   of the Appendix on the bounded connected components of $\Omega^-$, $\frac{\partial u}{\partial\nu_{\Omega^-}}$ is injective in the space in (\ref{mldc_corol:eipedd1}). In order to show the equality in (\ref{mldc_corol:eipedd3}), it suffices to assume that $m=0$, that $\mu\in V^{-1,\alpha}(\partial\Omega)$ solves
equation 
\begin{equation}\label{mldc_corol:eipedd5}
-\frac{1}{2}\mu 
+
 W_\Omega^t[\tilde{S}_{n,k;r},\mu]
 =0
\end{equation}
 and to show that $\mu\in C^{0,\alpha}(\partial\Omega)$. If $\mu\in V^{-1,\alpha}(\partial\Omega)$ solves equation (\ref{mldc_corol:eipedd5}), then 
  $v_\Omega^+[\tilde{S}_{n,k;r},\mu]\in C^{0,\alpha}(\overline{\Omega})_\Delta$
 (cf.~\cite[Thm.~5.3]{La25a}) and   the jump formula  for the distributional normal derivative of the  acoustic  single layer potential of \cite[Thm.~6.3]{La25a} implies that 
 \begin{equation}\label{mldc_corol:eipedd6}
 \frac{\partial}{\partial\nu_{\Omega}}v_\Omega^+[\tilde{S}_{n,k;r},\mu]=0\,.
\end{equation}
Then  equality (\ref{mldc_corol:eipedd6}) and the regularization Theorem of \cite[Thm.~1.2 (i)]{La25}  implies that $v_\Omega^+[\tilde{S}_{n,k;r},\mu]\in C^{1,\alpha}(\overline{\Omega})$. Since
 \[
v_\Omega^-[\tilde{S}_{n,k;r},\mu]_{|\partial\Omega}= v_\Omega^+[\tilde{S}_{n,k;r},\mu]_{|\partial\Omega}\in C^{1,\alpha}(\partial\Omega)\,,
 \]
 and   $v_\Omega^-[\tilde{S}_{n,k;r},\mu]$ in $C^{0,\alpha}_{{\mathrm{loc}}}(\overline{\Omega^-})_\Delta$
 (cf.~\cite[Thm.~6.3]{La25a}), then the regularization  Theorem of \cite[Thm.~1.2 (iv)]{La25} implies that $v_\Omega^-[\tilde{S}_{n,k;r},\mu]\in C^{1,\alpha}_{{\mathrm{loc}}}(\overline{\Omega^-})$. Thus the jump formulas for the distributional normal derivative of the  acoustic  single layer potential of \cite[Thm.~6.3]{La25a} imply that 
\[
\mu=-\frac{\partial}{\partial\nu_{\Omega^-}}v_\Omega^-[\tilde{S}_{n,k;r},\mu]-
\frac{\partial}{\partial\nu_{\Omega}}v_\Omega^+[\tilde{S}_{n,k;r},\mu]\in C^{0,\alpha}(\partial\Omega)\,.
\]
 Statement (ii) is an immediate consequence of statement (i) and of Theorem \ref{mldc_thm:eiped} (ii).\hfill  $\Box$ 

\vspace{\baselineskip}

\begin{remark}\label{mldc_rem:dimker-12i+wk2t} 
Under the assumptions of Corollary \ref{mldc_corol:eipedd} (ii),  equality (\ref{mldc_corol:eipedd1})  implies that the dimension of the space in (\ref{mldc_corol:eipedd4})  equals the sum of the geometric multiplicities of $k^{2}$ as a Dirichlet eigenvalue of $-\Delta$ in $(\Omega^{-})_{j}$ for $j\in {\mathbb{N}}\setminus\{0\}$ and $j\leq \varkappa^{-}$.  (For the definition of $\varkappa^-$, see right after (\ref{mldc_eq:exto})).

In particular, if we also assume that for  each $j\in {\mathbb{N}}\setminus\{0\}$ such that $j\leq \varkappa^{-}$, $k^2$ is not a Dirichlet eigenvalue for $-\Delta$ in $(\Omega^{-})_{j}$, then the dimension of the space in (\ref{mldc_corol:eipedd4})  equals $0$. For example, this is the case if we assume  that $\Omega^-$ is connected.
\end{remark}

Next we note that under appropriate assumptions, the space of solutions of the integral equation (\ref{mldc_thm:eiped4}) is stable under conjugation. 
\begin{remark}\label{mldc_rem:-i+wtconj} 
Let $m\in{\mathbb{N}}$, $\alpha\in]0,1[$. Let $\Omega$ be a bounded open subset of ${\mathbb{R}}^{n}$ of class $C^{\max\{1,m\},\alpha}$. Let 
$k\in {\mathbb{C}}\setminus ]-\infty,0]$, ${\mathrm{Im}}\,k\geq 0$. Assume that 
$k^{2}$ is not a Neumann eigenvalue for $-\Delta$ in $\Omega$.\par

We now show that if  $\mu \in V^{m-1,\alpha}(\partial\Omega)$ satisfies equation (\ref{mldc_thm:eiped4}), then   $\mu$ belongs to $C^{\max\{1,m\}-1,\alpha}(\partial\Omega)$ and  $\overline{\mu}$ satisfies the same equation. Indeed, if $k^2$ is not real, then Corollary \ref{mldc_corol:eipedd} (ii) implies that equation (\ref{mldc_thm:eiped4}) has only the trivial solution. If instead $k^2$ is real and if $\mu $ in $V^{m-1,\alpha}(\partial\Omega)$ satisfies equation (\ref{mldc_thm:eiped4}),  then 
Corollary \ref{mldc_corol:eipedd} (i) implies that $\mu$ belongs to $C^{\max\{1,m\}-1,\alpha}(\partial\Omega)$, Corollary \ref{mldc_corol:eipedd}  (ii) implies that $v^{-}_\Omega[\tilde{S}_{n,k;r},-\mu]$ belongs to the set in 
(\ref{mldc_corol:eipedd1}) and that $\mu=\frac{\partial}{\partial\nu_{\Omega^-}}v^{-}[\tilde{S}_{n,k;r},-\mu]$. Since $k^2$ is real, $\overline{v^{-}_\Omega[\tilde{S}_{n,k;r},\mu]}$ belongs to the set in (\ref{mldc_corol:eipedd1}) too
 and thus
Corollary \ref{mldc_corol:eipedd} (i) implies that $\frac{\partial}{\partial\nu_{\Omega^-}}\overline{v^{-}[\tilde{S}_{n,k;r},-\mu]}$
satisfies equation (\ref{mldc_thm:eiped4}). Since
\[
 \frac{\partial}{\partial\nu_{\Omega^-}}\overline{v^{-}_\Omega[\tilde{S}_{n,k;r},-\mu] } 
 =
 \overline{ \frac{\partial}{\partial\nu_{\Omega^-}}v^{-}_\Omega[\tilde{S}_{n,k;r},-\mu] }=\overline{\mu}  \,,
\]
 then $\overline{\mu}$  satisfies equation (\ref{mldc_thm:eiped4}).
\end{remark}
As is well known, if    $\alpha\in]0,1[$ and $\Omega$ is a bounded open subset of ${\mathbb{R}}^{n}$ of class $C^{1,\alpha}$ such that $\Omega^-$ is not connected, then the exterior Dirichlet problem (\ref{mldc_thm:eiped1}) has nontrivial solutions precisely  when $k^2$ belongs to an increasing sequence $(\lambda_{j}^{({\mathcal{D}})}[\Omega^-])_{j\in{\mathbb{N}} }$ of   eigenvalues in $]0,+\infty[$, and we can write such a sequence so that each eigenvalue $\lambda_{j}^{({\mathcal{D}})}[\Omega^-]$ is repeated as many times as its multiplicity.

\section{A nonvariational form of the exterior Dirichlet problem for the Helmholtz equation}\label{mldc_sec:expddh}
We now consider the exterior  Dirichlet problem in the case in which the Dirichlet datum is in the space $C^{m,\alpha}(\partial\Omega)$ (nonvariational in case $m=0$). For the corresponding classical result in case $\varkappa^-=0$, we refer to  Colton and Kress \cite[Thm.~3.21]{CoKr92}. We first introduce the following classical necessary condition for the existence.

\begin{theorem}\label{mldc_thm:nexextdirhe}
  Let $m\in {\mathbb{N}}$,   $\alpha\in]0,1[$. Let $\Omega$ be a bounded open subset of ${\mathbb{R}}^{n}$ of class $C^{\max\{1,m\},\alpha}$.  Let $k\in {\mathbb{C}}\setminus ]-\infty,0]$, ${\mathrm{Im}}\,k\geq 0$.   Let $g\in C^{m,\alpha}(\partial\Omega)$. If   the exterior Dirichlet problem 
\begin{equation}
\label{mldc_thm:nexextdirhe1}
\left\{
\begin{array}{ll}
\Delta u+k^2 u=0 &{\text{in}}\ \Omega^-\,,
\\
u=g
&{\text{on}}\ \partial\Omega\,,
\\
u\ \text{satisfies\ the\ outgoing}\  k-\text{radiation\ condition}\,,
\end{array}
\right.
\end{equation}
has a solution $u\in C^{m,\alpha}_{ {\mathrm{loc}} }(\overline{\Omega^-})$, then
\begin{equation}\label{mldc_thm:nexextdirhe2}
\int_{\partial\Omega}g\frac{\partial v}{\partial\nu_{\Omega^-}}\,d\sigma=0  \,,
\end{equation}
for all $v\in C^{\max\{1,m\},\alpha}_{ {\mathrm{loc}} }(\overline{\Omega^-})$  such that
\begin{equation}\label{mldc_thm:nexextdirhe3}
\left\{
\begin{array}{ll}
 \Delta v+k^2 v=0 &\text{in}\ \Omega^-\,,
 \\
 v=0		&\text{on}\ \partial\Omega\,,
 \\
 v\ \text{satisfies\ the\ outgoing}\  k-\text{radiation\ condition}\,.
\end{array}
\right.
\end{equation}
 \end{theorem}
{\bf Proof.} Assume that a solution $u\in C^{m,\alpha}_{ {\mathrm{loc}} }(\overline{\Omega^-})$  exists and that $v\in C^{\max\{1,m\},\alpha}(\overline{\Omega})$ satisfies the homogeneous Dirichlet problem of (\ref{mldc_thm:nexextdirhe3}). Then Remark \ref{mldc_rem:helpo} and the second Green Identity for exterior domains in distributional form of (cf.~\cite[Thm.~8.18]{La25b})   imply that 
\[
\int_{\partial\Omega}g\frac{\partial v}{\partial\nu_{\Omega^-}}\,d\sigma=
\int_{\partial\Omega}u\frac{\partial v}{\partial\nu_{\Omega^-}}\,d\sigma
=\langle  \frac{\partial u}{\partial\nu_{\Omega^-}}, v_{|\partial\Omega}\rangle
=0 
\]
and thus condition (\ref{mldc_thm:nexextdirhe2}) holds true.\hfill  $\Box$ 

\vspace{\baselineskip}

By Theorem \ref{mldc_thm:nexextdirhe}, a necessary condition for the solvability of the exterior Dirichlet problem is that condition (\ref{mldc_thm:nexextdirhe2}) holds true. Thus we now introduce the following notation. If $m\in{\mathbb{N}}$, $\alpha\in]0,1[$,  $\Omega$ is a bounded open subset of ${\mathbb{R}}^{n}$ of class $C^{\max\{1,m\},\alpha}$,   $k\in {\mathbb{C}}\setminus ]-\infty,0]$, ${\mathrm{Im}}\,k\geq 0$,  then we set
\begin{eqnarray}\label{mldc_eq:xmal-}
\lefteqn{
X^{m,\alpha}_{k^2,-}\equiv\biggl\{
\frac{\partial v}{\partial\nu_{\Omega^-}}:\, v\in C^{\max\{1,m\},\alpha}_{ {\mathrm{loc}} }(\overline{\Omega^-})\,,\ 
\Delta v+k^2 v=0 \ \text{in}\ \Omega^-\,, 
}
\\ \nonumber
&&\qquad\qquad 
v=0\ \text{on}\   \partial\Omega\,,\  v\ \text{satisfies\ the\ outgoing}\  k-\text{radiation\ condition}
\biggr\}
\end{eqnarray}
and
\begin{equation}\label{mldc_eq:xmal-1}
(X^{m,\alpha}_{k^2,-})^\perp\equiv
\left\{
g\in C^{m,\alpha}(\partial\Omega):\, \int_{\partial\Omega}gf\,d\sigma=0\ \forall f\in X^{m,\alpha}_{k^2,-}
\right\}\,.
\end{equation}
By Corollary \ref{mldc_corol:eipedd} (i),    we have
 \begin{eqnarray}\label{mldc_eq:xmal-2}
\lefteqn{
X^{m,\alpha}_{k^2,-}=\biggl\{
\frac{\partial v}{\partial\nu_{\Omega^-}}:\, v\in C^{\max\{1,m\},\alpha}_{ {\mathrm{loc}} }(\overline{\Omega^-})\,,\ 
\Delta v+k^2 v=0 \ \text{in}\ \Omega^-\,, 
}
\\ \nonumber
&&\qquad\qquad\qquad\qquad 
v=0\ \text{on}\   \partial\Omega\,,\  v(x)=0\ \forall x\in (\Omega^-)_0
\biggr\}\,.
\end{eqnarray}
\begin{lemma}\label{mldc_lem:xmalco-}
 Let $m\in{\mathbb{N}}$, $\alpha\in]0,1[$. Let $\Omega$ be a bounded open subset of ${\mathbb{R}}^{n}$ of class $C^{\max\{1,m\},\alpha}$. Let $k\in {\mathbb{C}}\setminus ]-\infty,0]$, ${\mathrm{Im}}\,k\geq 0$.   If
 $\mu\in  X^{m,\alpha}_{k^2,-}$, then $\overline{\mu}\in  X^{m,\alpha}_{k^2,-}$. Moreover, if $g\in (X^{m,\alpha}_{k^2,-})^\perp$, then $\overline{g}\in  (X^{m,\alpha}_{k^2,-})^\perp$. 
\end{lemma}
{\bf Proof.} If $k^2$ is not real, then $X^{m,\alpha}_{k^2,-}=\{0\}$ and both statements are true. Let $k^2$ be real. If $\mu\in  X^{m,\alpha}_{k^2,-}$, then there exists $v$ as in (\ref{mldc_eq:xmal-2}) such that $\mu=\frac{\partial v}{\partial\nu_{\Omega^-}}$.  Since $k^2$ is real, then    $\Delta \overline{v}+k^2\overline{v}=0$ in $\Omega^-$, $\overline{v}=0$ on $\partial\Omega$ and  
$\overline{v}(x)=0$  for all $x\in (\Omega^-)_0$ and accordingly
 $\frac{\partial \overline{v}}{\partial\nu_{\Omega^-}}\in X^{m,\alpha}_{k^2,-}$ (cf.~(\ref{mldc_eq:xmal-2})). Since
\[
\overline{\mu}=\overline{\frac{\partial v}{\partial\nu_{\Omega^-}}}
=\frac{\partial \overline{v}}{\partial\nu_{\Omega^-}}\,,
\]
we have $\overline{\mu}\in X^{m,\alpha}_{k^2,-}$. Then by exploiting the definition of $(X^{m,\alpha}_{k^2,-})^\perp$, we see that if $g\in (X^{m,\alpha}_{k^2,-})^\perp$, then $\overline{g}\in  (X^{m,\alpha}_{k^2,-})^\perp$.\hfill  $\Box$ 

\vspace{\baselineskip}

\begin{proposition}\label{mldc_prop:1/2+w,vxmap-}
 Let $m\in{\mathbb{N}}$, $\alpha\in]0,1[$. Let $\Omega$ be a bounded open subset of ${\mathbb{R}}^{n}$ of class $C^{\max\{1,m\},\alpha}$. Let $k\in {\mathbb{C}}\setminus ]-\infty,0]$, ${\mathrm{Im}}\,k\geq 0$.    Then the following statements hold true. 
 \begin{enumerate}
\item[(i)] $ -\frac{1}{2}\phi+W_\Omega[\tilde{S}_{n,k;r},\phi]\in (X^{m,\alpha}_{k^2,-})^\perp$ for all $\phi\in  C^{m,\alpha}(\partial\Omega)$.
 \item[(ii)] $ V_\Omega[\tilde{S}_{n,k;r},\psi] \in (X^{m,\alpha}_{k^2,-})^\perp$ for all $\psi\in C^{\max\{1,m\}-1,\alpha}(\partial\Omega)$.
\end{enumerate}
\end{proposition}
{\bf Proof.} (i) If $\phi\in C^{m,\alpha}(\partial\Omega)$ and $f\in X^{m,\alpha}_{k^2,-}$, then 
\[
\int_{\partial\Omega}\left(- \frac{1}{2}\phi+W_\Omega[\tilde{S}_{n,k;r},\phi]\right)f\,d\sigma
=\int_{\partial\Omega}\phi \left(-\frac{1}{2}f+
W_\Omega^t[\tilde{S}_{n,k;r},f]\right)\,d\sigma
\]
and Corollary \ref{mldc_corol:eipedd} (i) ensures that the right hand side equals $0$. 

(ii) If $\psi\in C^{\max\{1,m\}-1,\alpha}(\partial\Omega)$, $f\in X^{m,\alpha}_{k^2,-}$, then there exists $v$ as in (\ref{mldc_eq:xmal-}) such that $f=\frac{\partial v}{\partial\nu_{\Omega^-}}$ and thus the second Green Identity for exterior domains implies that
\begin{eqnarray*}
\lefteqn{
\int_{\partial\Omega}v_\Omega^-[\tilde{S}_{n,k;r},\psi]f\,d\sigma
}
\\ \nonumber
&&\qquad
=
\int_{\partial\Omega}v_\Omega^-[\tilde{S}_{n,k;r},\psi]\frac{\partial v}{\partial\nu_{\Omega^-}}\,d\sigma=
\int_{\partial\Omega}\frac{\partial }{\partial\nu_{\Omega^-}}v_\Omega^-[\tilde{S}_{n,k;r},\psi]v\,d\sigma=0\,,
\end{eqnarray*}
(cf.~\textit{e.g.}, \cite[Thm.~8.18]{La25b}).\hfill  $\Box$ 

\vspace{\baselineskip}

\begin{proposition}\label{mldc_prop:silamokeri-wt-}
Let $m\in{\mathbb{N}}$, $\alpha\in]0,1[$. Let $\Omega$ be a bounded open subset of ${\mathbb{R}}^{n}$ of class $C^{\max\{1,m\},\alpha}$.  
 Then the following statements hold.
\begin{enumerate}
\item[(i)] Let $\lambda\in {\mathbb{C}}$. Let $S_{n,\lambda}$ be a fundamental solution of $\Delta+\lambda$. If 
\[
\phi\in {\mathrm{Ker}}\left(
-\frac{1}{2}I+W_\Omega^t[S_{n,\lambda},\cdot]
\right)_{|C^{\max\{1,m\}-1,\alpha}(\partial\Omega)} 
\]
 (cf.~Remark \ref{mldc_rem:wtnotation}), 
 then $v^+\equiv v_\Omega^+[S_{n,\lambda},\phi]$ belongs to  $C^{\max\{1,m\},\alpha}(\overline{\Omega})$  
 and solves the Neumann problem
 \begin{equation}\label{mldc_prop:silamokeri-wt-0}
 \Delta v^++\lambda v^+=0\qquad\text{in}\ \Omega\,,
 \qquad\frac{\partial v^+}{\partial\nu_{\Omega}}=0\qquad\text{on}\ \partial\Omega\,.
 \end{equation}
\item[(ii)] Let $k\in {\mathbb{C}}\setminus ]-\infty,0]$, ${\mathrm{Im}}\,k\geq 0$. If 
\[
\phi\in {\mathrm{Ker}}\left(
-\frac{1}{2}I+W_\Omega^t[\tilde{S}_{n,k;r},\cdot] 
\right)_{|C^{\max\{1,m\}-1,\alpha}(\partial\Omega)}
\]
 and $v_\Omega^+[\tilde{S}_{n,k;r},\phi]_{|\partial\Omega}=0$ on $\partial\Omega$, then 
$\phi\in X^{m,\alpha}_{k^2,-}$.
\item[(iii)]  Let $k\in {\mathbb{C}}\setminus ]-\infty,0]$, ${\mathrm{Im}}\,k\geq 0$. If \[
\phi\in {\mathrm{Ker}}\left(
-\frac{1}{2}I+W_\Omega^t[\tilde{S}_{n,k;r},\cdot] 
\right)_{|C^{\max\{1,m\}-1,\alpha}(\partial\Omega)}
\]
 and $\phi\in X^{m,\alpha}_{k^2,-}$, then $v_\Omega^+[\tilde{S}_{n,k;r},\phi]=0$ in $\overline{\Omega}$.
\end{enumerate}
\end{proposition} 
{\bf Proof.} (i)   Classical regularity results on the acoustic single   layer potential and the classical jump formulas for the normal derivative of the  acoustic  single   layer potential imply that 
$v_\Omega^+[S_{n,\lambda},\phi]$ belongs to $ C^{\max\{1,m\},\alpha}(\overline{\Omega}) $   
and solves the Neumann problem (\ref{mldc_prop:silamokeri-wt-0}) (cf.~\textit{e.g.}, \cite[Thm.~7.1]{DoLa17}). (ii) By classical regularity results on the  acoustic  single  layer potential,   
$v_\Omega^-[\tilde{S}_{n,k;r},\phi]$ belongs to $ C^{\max\{1,m\},\alpha}_{{\mathrm{loc}}}(\overline{\Omega^-}) $ 
(cf.~\textit{e.g.}, \cite[Thm.~7.1]{DoLa17}). 
Also, the classical jump conditions for the  acoustic  single  layer potential imply that
\[
v_\Omega^-[\tilde{S}_{n,k;r},\phi]_{|\partial\Omega}=v_\Omega^+[\tilde{S}_{n,k;r},\phi]_{|\partial\Omega}=0
\]
 (cf.~\textit{e.g.}, \cite[Lem.~4.2 (i), Lem.~6.2]{DoLa17}). Then
$v_\Omega^-[\tilde{S}_{n,k;r},\phi]$ satisfies the Dirichlet problem
\[
\Delta v_\Omega^-[\tilde{S}_{n,k;r},\phi]+k^2 v_\Omega^-[\tilde{S}_{n,k;r},\phi]=0\quad\text{in}\ \Omega^-\,,\qquad v_\Omega^-[\tilde{S}_{n,k;r},\phi]_{|\partial\Omega}=0 \,.
\]
Since $v_\Omega^-[\tilde{S}_{n,k;r},\phi]$ is known to satisfy the outgoing $k$-radiation condition (cf.~\textit{e.g.}, \cite[Thm.~6.20 (i)]{La25b}), then
the definition 
(\ref{mldc_eq:xmal-}) of $X^{m,\alpha}_{k^2,-}$  implies that  
$
\frac{\partial}{\partial\nu_{\Omega^-}}v_\Omega^-[\tilde{S}_{n,k;r},\phi]$ belongs to $ X^{m,\alpha}_{k^2,-}$. On the other hand, statement (i) implies that $\frac{\partial}{\partial\nu_{\Omega}}v_\Omega^+[\tilde{S}_{n,k;r},\phi]=0$ on $\partial\Omega$. Then the classical jump conditions for the normal derivative of the  acoustic  single  layer potential imply that
\[
\phi=-\frac{\partial}{\partial\nu_{\Omega^-}}v_\Omega^-[\tilde{S}_{n,k;r},\phi]
-\frac{\partial}{\partial\nu_{\Omega}}v_\Omega^+[\tilde{S}_{n,k;r},\phi]
=-\frac{\partial}{\partial\nu_{\Omega^-}}v_\Omega^-[\tilde{S}_{n,k;r},\phi]\in X^{m,\alpha}_{k^2,-} 
\]
(cf.~\textit{e.g.}, \cite[Thm.~7.1]{DoLa17}).   (iii) Since 
\[
\phi\in {\mathrm{Ker}}\left(
-\frac{1}{2}I+W_\Omega^t[\tilde{S}_{n,k;r},\cdot]
\right)_{|C^{\max\{1,m\}-1,\alpha}(\partial\Omega)} 
\,,
\]
classical regularity results on the  acoustic  single  layer potential and the classical jump formulas for the normal derivative of the  acoustic  single   layer potential imply that 
$v_\Omega^+[\tilde{S}_{n,k;r},\phi]$ belongs to $ C^{\max\{1,m\},\alpha}(\overline{\Omega}) $ and solves the Neumann problem (\ref{mldc_prop:silamokeri-wt-0}) with $\lambda=k^2$ (cf.~\textit{e.g.}, \cite[Thm.~7.1]{DoLa17}).  By the Holmgren uniqueness  Corollary \ref{mldc_corol:holmgr}   of the Appendix, it suffices to show that $v_\Omega^+\left[\tilde{S}_{n,k;r},\phi\right]$ equals $0$ on $\partial\Omega$.
Since  $\phi\in X^{m,\alpha}_{k^2,-}$, then there exists a function  $v$ in $C^{\max\{1,m\},\alpha}_{{\mathrm{loc}}}(\overline{\Omega^-})$ such that
\begin{eqnarray*}
\lefteqn{\Delta v+k^2 v=0 \ \text{in}\ \Omega^-\,,\  v=0\ \text{on}\ \partial\Omega\,,
}
\\ \nonumber
&&\qquad
v\ \text{satisfies\ the\ outgoing}\  k-\text{radiation\ condition}
\end{eqnarray*}
and $\phi=\frac{\partial v}{\partial\nu_{\Omega^-}}$. 
Then the classical third Green Identity for exterior domains implies that 
\[
v=-w_\Omega^-[\tilde{S}_{n,k;r},v_{|\partial\Omega}]-v_\Omega^-\left[\tilde{S}_{n,k;r},\frac{\partial v}{\partial\nu_{\Omega^-}}\right]=-v_\Omega^-\left[\tilde{S}_{n,k;r},\frac{\partial v}{\partial\nu_{\Omega^-}}\right]
\qquad\text{in}\ \Omega^-
\]
(cf.~\textit{e.g.}, \cite[Thm.~8.24]{La25b}) and the classical jump conditions for the    acoustic  single  layer potential imply that
\[
v_\Omega^+\left[\tilde{S}_{n,k;r},\phi\right]_{|\partial\Omega}=
v_\Omega^+\left[\tilde{S}_{n,k;r},\frac{\partial v}{\partial\nu_{\Omega^-}}\right]_{|\partial\Omega}=
v_\Omega^-\left[\tilde{S}_{n,k;r},\frac{\partial v}{\partial\nu_{\Omega^-}}\right]_{|\partial\Omega}
=-v_{|\partial\Omega}=0
\]
(cf.~\textit{e.g.}, \cite[Lem.~4.2 (i), Lem.~6.2]{DoLa17}).   Hence, the proof is complete. \hfill  $\Box$ 

\vspace{\baselineskip}

\begin{proposition}\label{mldc_prop:isomquotd-}
 Let $m\in{\mathbb{N}}$, $\alpha\in]0,1[$. Let $\Omega$ be a bounded open subset of ${\mathbb{R}}^{n}$ of class $C^{\max\{1,m\},\alpha}$.  
  Then the following statements hold true (cf.~Remark \ref{mldc_rem:wtnotation}). 
\begin{enumerate}
\item[(i)] Let $\lambda\in {\mathbb{C}}$. Let $S_{n,\lambda}$ be a fundamental solution of $\Delta+\lambda$.  Then 
\[
{\mathrm{Im}}\, \left(
-\frac{1}{2}I+W_\Omega[S_{n,\lambda},\cdot]
\right)_{|C^{m,\alpha}(\partial\Omega)}
\]
 is a closed subspace of
$C^{m,\alpha}(\partial\Omega)$ and 
the quotient space 
\[
\frac{
C^{m,\alpha}(\partial\Omega)
}{
{\mathrm{Im}}\, \left(
-\frac{1}{2}I+W_\Omega[S_{n,\lambda},\cdot]
\right)_{|C^{m,\alpha}(\partial\Omega)}
}
\]
 has dimension equal to the finite dimension of  
 \[
 {\mathrm{Ker}}\, \left(
-\frac{1}{2}I+W_\Omega^t[S_{n,\lambda},\cdot]
\right)_{|C^{\max\{1,m\}-1,\alpha}(\partial\Omega)}
\,.
\]
\item[(ii)] Let $k\in {\mathbb{C}}\setminus ]-\infty,0]$, ${\mathrm{Im}}\,k\geq 0$.  Then ${\mathrm{Im}}\, \left(
-\frac{1}{2}I+W_\Omega[\tilde{S}_{n,k;r},\cdot]
\right)_{|C^{m,\alpha}(\partial\Omega)}
$ is a closed subspace of
$(X^{m,\alpha}_{k^2,-})^\perp$. 
\item[(iii)] Let $k\in {\mathbb{C}}\setminus ]-\infty,0]$, ${\mathrm{Im}}\,k\geq 0$.  Then the quotient space $\frac{
C^{m,\alpha}(\partial\Omega)
}{
(X^{m,\alpha}_{k^2,-})^\perp
}$ is linearly homeomorphic to $X^{m,\alpha}_{k^2,-}$ and has finite dimension.
\item[(iv)] Let $k\in {\mathbb{C}}\setminus ]-\infty,0]$, ${\mathrm{Im}}\,k\geq 0$. Then the  following formula holds
\begin{eqnarray*}
\lefteqn{
{\mathrm{dim}}\,\left(
\frac{(X^{m,\alpha}_{k^2,-})^\perp}{
{\mathrm{Im}}\, \left(
-\frac{1}{2}I+W_\Omega[\tilde{S}_{n,k;r},\cdot]
\right)_{|C^{m,\alpha}(\partial\Omega)}
}\right)
}
\\ \nonumber
&&\qquad
={\mathrm{dim}}\,{\mathrm{Ker}}\, \left(
-\frac{1}{2}I+W_\Omega^t[\tilde{S}_{n,k;r},\cdot]
\right)_{|C^{\max\{1,m\}-1,\alpha}(\partial\Omega)}
-{\mathrm{dim}}\,X^{m,\alpha}_{k^2,-}\,.
\end{eqnarray*}
\end{enumerate}
\end{proposition}
{\bf Proof.} (i) Since $W_\Omega[S_{n,\lambda},\cdot]_{|C^{m,\alpha}(\partial\Omega)} $ is compact in $C^{m,\alpha}(\partial\Omega)$ (cf.~\textit{e.g.}, \cite[Thm.~7.4]{DoLa17} in case $m=0$ and  \cite[Cor.~9.1]{DoLa17} in case $m\geq 1$) and
$W_\Omega^t[S_{n,\lambda},\cdot]_{|C^{\max\{1,m\}-1,\alpha}(\partial\Omega)} $ is compact in $C^{\max\{1,m\}-1,\alpha}(\partial\Omega)$ (cf.~\textit{e.g.}, \cite[Cor.~10.1]{DoLa17}), then the Fredholm Alternative Theorem of Wendland \cite{We67}, \cite{We70} in the duality pairing
\[
\left(C^{\max\{1,m\}-1,\alpha} (\partial\Omega),C^{m,\alpha} (\partial\Omega)\right) 
\]
(cf.~\textit{e.g.}, Kress~\cite[Thm.~4.17]{Kr14}) implies the validity of the following equality
\begin{eqnarray*}
\lefteqn{
{\mathrm{Im}}\, \left(
-\frac{1}{2}I+W_\Omega[S_{n,\lambda},\cdot]
\right)_{|C^{m,\alpha}(\partial\Omega)}
}
\\ \nonumber
&&\qquad
=\left(
{\mathrm{Ker}}\, \left(
-\frac{1}{2}I+W_\Omega^t[S_{n,\lambda},\cdot]
\right)_{|C^{\max\{1,m\}-1,\alpha}(\partial\Omega)}
\right)^\perp\,,
\end{eqnarray*}
where the orthogonality has to be understood according to the above duality pairing
and that the finite dimension $s_{-}$ of ${\mathrm{Ker}}\, \left(
-\frac{1}{2}I+W_\Omega^t[S_{n,\lambda},\cdot]
\right)_{|C^{\max\{1,m\}-1,\alpha}(\partial\Omega)}$
equals that of ${\mathrm{Ker}}\, \left(
-\frac{1}{2}I+W_\Omega[S_{n,\lambda},\cdot]
\right)
_{|C^{m,\alpha}(\partial\Omega)}$. If $s_{-}=0$, then 
\[
\left(
{\mathrm{Ker}}\, \left(
-\frac{1}{2}I+W_\Omega^t[S_{n,\lambda},\cdot]
\right)_{|C^{\max\{1,m\}-1,\alpha}(\partial\Omega)}
\right)^\perp=C^{m,\alpha}(\partial\Omega)
\]
and statement (i) is obviously true. Thus we can assume that $s_{-}\geq 1$. 
Let $\{\mu_j\}_{j=1}^{s_{-}}$ be a basis of ${\mathrm{Ker}}\, \left(
-\frac{1}{2}I+W_\Omega^t[S_{n,\lambda},\cdot]
\right)_{|C^{\max\{1,m\}-1,\alpha}(\partial\Omega)}
$. 

Then we consider the map $\Lambda$ from $C^{m,\alpha}(\partial\Omega)$ to ${\mathbb{C}}^{s_-}$ that is defined by the equality
\[
\Lambda[g]\equiv\left(\int_{\partial\Omega}g\mu_1\,d\sigma, 
\dots,
\int_{\partial\Omega}g\mu_{s_{-}}\,d\sigma
\right)\qquad\forall g\in C^{m,\alpha}(\partial\Omega)\,.
\]
If $\Lambda$ is not surjective, then there exists $(c_1,\dots,c_{s_{-}})\in  {\mathbb{C}}^{s_{-}}$ such that
\[
0=\sum_{j=1}^{s_{-}}c_j\int_{\partial\Omega}g\mu_j\,d\sigma=\int_{\partial\Omega}g\sum_{j=1}^{s_{-}}c_j\mu_j\,d\sigma
\qquad\forall g\in C^{m,\alpha}(\partial\Omega)\,.
\]
Hence, $\sum_{j=1}^{s_{-}}c_j\mu_j=0$, a contradiction. Then $\Lambda$ is surjective. Since
\begin{eqnarray*}
\lefteqn{
{\mathrm{Ker}}\, \Lambda=\left(
{\mathrm{Ker}}\, \left(
-\frac{1}{2}I+W_\Omega^t[S_{n,\lambda},\cdot]
\right)_{|C^{\max\{1,m\}-1,\alpha}(\partial\Omega)}
\right)^\perp
}
\\ \nonumber
&&\qquad\qquad\qquad\qquad
=
{\mathrm{Im}}\, \left(
-\frac{1}{2}I+W_\Omega[S_{n,\lambda},\cdot]
\right)_{|C^{m,\alpha}(\partial\Omega)}
\,,
\end{eqnarray*}
 the Homomorphism Theorem implies that 
\[
{\mathbb{C}}^{s_{-}} \quad \text{is\ isomorphic\ to}\quad  
\frac{
C^{m,\alpha}(\partial\Omega)
}{
{\mathrm{Im}}\, \left(
-\frac{1}{2}I+W_\Omega[S_{n,\lambda},\cdot]
\right)_{|C^{m,\alpha}(\partial\Omega)}
}
\]
and accordingly, statement (i) holds true. Statement (ii) holds by statement (i) and by  Proposition \ref{mldc_prop:1/2+w,vxmap-} (i).
Since $X^{m,\alpha}_{k^2,-}$ can be identified with a linear subspace of the dual of 
$C^{m,\alpha}(\partial\Omega)$, it is known that the dual
\[
\left(\frac{
C^{m,\alpha}(\partial\Omega)
}{
(X^{m,\alpha}_{k^2,-})^\perp
}\right)'\quad\text{is\ isomorphic\ to}\ (X^{m,\alpha}_{k^2,-})^{\perp\perp}
\]
(cf.~\textit{e.g.}, Rudin, \cite[Thm.~4.9, Chap. IV]{Ru91}). On the other hand, $X^{m,\alpha}_{k^2,-}$ is finite dimensional and accordingly, $X^{m,\alpha}_{k^2,-}$  is weak$^\ast$ closed and
\[
(X^{m,\alpha}_{k^2,-})^{\perp\perp}=X^{m,\alpha}_{k^2,-}
\]
(cf.~\textit{e.g.}, Rudin, \cite[Thm.~4.7, Chap. IV]{Ru91}). Hence, $\left(\frac{
C^{m,\alpha}(\partial\Omega)
}{
(X^{m,\alpha}_{k^2,-})^\perp
}\right)'$ is linearly  isomorphic to $X^{m,\alpha}_{k^2,-}$ and has finite dimension. Hence, 
$\frac{
C^{m,\alpha}(\partial\Omega)
}{
(X^{m,\alpha}_{k^2,-})^\perp
}$ is also linearly isomorphic to $X^{m,\alpha}_{k^2,-}$. On the other hand linear isomorphisms between normed spaces of finite dimension are also homeomorphisms and thus statement (iii) holds true. 

 By the Homomorphism Theorem (cf.~\textit{e.g.}, Greub \cite[7. b), p.~49]{Gr75}), the quotient $\frac{
C^{m,\alpha}(\partial\Omega)
}{
(X^{m,\alpha}_{k^2,-})^\perp
}$ is linearly isomorphic to the quotient
\[
\left(\frac{
C^{m,\alpha}(\partial\Omega)
}{
{\mathrm{Im}}\, \left(
-\frac{1}{2}I+W_\Omega[\tilde{S}_{n,k;r},\cdot]
\right)_{|C^{m,\alpha}(\partial\Omega)}
}\right)
{\displaystyle\diagup}
\left(\frac{
(X^{m,\alpha}_{k^2,-})^\perp
}{
{\mathrm{Im}}\, \left(
-\frac{1}{2}I+W_\Omega[\tilde{S}_{n,k;r},\cdot]
\right)_{|C^{m,\alpha}(\partial\Omega)}
}\right) 
\]
that has dimension equal to 
\begin{eqnarray*}
\lefteqn{
{\mathrm{dim}}\,{\mathrm{Ker}}\, \left(
-\frac{1}{2}I+W_\Omega^t[\tilde{S}_{n,k;r},\cdot]
\right)_{|C^{\max\{1,m\}-1,\alpha}(\partial\Omega)}
}
\\ \nonumber
&&\qquad
-{\mathrm{dim}}\,
\left(\frac{
(X^{m,\alpha}_{k^2,-})^\perp
}{
{\mathrm{Im}}\, \left(
-\frac{1}{2}I+W_\Omega[\tilde{S}_{n,k;r},\cdot]
\right)_{|C^{m,\alpha}(\partial\Omega)}
}\right)
\end{eqnarray*}
by statement (i). Hence, statement (iii) implies the validity of the formula in statement (iv)  and the proof is complete.\hfill  $\Box$ 

\vspace{\baselineskip}

Then we can prove the following statement that says that the space of compatible data $X^{m,\alpha}_{k^2,-}$ can be written as an algebraic sum of data that equal the boundary values of an  acoustic  double   layer potential and of an  acoustic  single  layer potential. 
\begin{theorem}\label{mldc_thm:codasumwv-}
Let $m\in{\mathbb{N}}$, $\alpha\in]0,1[$. Let $\Omega$ be a bounded open subset of ${\mathbb{R}}^{n}$ of class $C^{\max\{1,m\},\alpha}$.  
 Let $k\in {\mathbb{C}}\setminus ]-\infty,0]$, ${\mathrm{Im}}\,k\geq 0$.  Then
\begin{equation}\label{mldc_thm:codasumwv-1}
(X^{m,\alpha}_{k^2,-})^\perp={\mathrm{Im}}\, \left(
-\frac{1}{2}I+W_\Omega[\tilde{S}_{n,k;r},\cdot]
\right)_{|C^{m,\alpha}(\partial\Omega)}
+
V_\Omega[\tilde{S}_{n,k;r},C^{\max\{1,m\}-1,\alpha}(\partial\Omega)]\,.
\end{equation}
\end{theorem}
{\bf Proof.} By Proposition \ref{mldc_prop:isomquotd-} (i), (ii), (iii), ${\mathrm{Im}}\, \left(
-\frac{1}{2}I+W_\Omega[\tilde{S}_{n,k;r},\cdot]
\right)_{|C^{m,\alpha}(\partial\Omega)}
$ is a closed subspace of
$(X^{m,\alpha}_{k^2,-})^\perp$ of finite codimension. Let $\pi$ be the canonical projection of 
$(X^{m,\alpha}_{k^2,-})^\perp$ onto the finite dimensional quotient space 
\[
Q^{m,\alpha}_{k^2,-}\equiv 
\frac{(X^{m,\alpha}_{k^2,-})^\perp}{
{\mathrm{Im}}\, \left(
-\frac{1}{2}I+W_\Omega[\tilde{S}_{n,k;r},\cdot]\right)_{|C^{m,\alpha}(\partial\Omega)} 
} \,.
\]
 By Proposition \ref{mldc_prop:1/2+w,vxmap-}  (ii), the image $V_\Omega
 [\tilde{S}_{n,k;r},C^{\max\{1,m\}-1,\alpha}(\partial\Omega)]$ is a subspace of $(X^{m,\alpha}_{k^2,-})^\perp$ and thus  the $\pi$-image $\pi\left[Y^{m,\alpha}_{k^2,-}\right]$
of
\[
Y^{m,\alpha}_{k^2,-}\equiv  {\mathrm{Im}}\, \left(
-\frac{1}{2}I+W_\Omega[\tilde{S}_{n,k;r},\cdot]\right)
_{|C^{m,\alpha}(\partial\Omega)}
+
V_\Omega[\tilde{S}_{n,k;r},C^{\max\{1,m\}-1,\alpha}(\partial\Omega)] 
\]
is a   subspace of the finite dimensional Banach quotient space $ 
Q^{m,\alpha}_{k^2,-}$ and is accordingly closed in $Q^{m,\alpha}_{k^2,-}$. Since the canonical projection $\pi$ is continuous, then the preimage
\[
\pi^{\leftarrow}\left[\pi\left[Y^{m,\alpha}_{k^2,-} \right]\right]
\]
is closed in $(X^{m,\alpha}_{k^2,-})^\perp$. Now a simple computation based on the definition of $\pi$  shows that
\[
\pi^{\leftarrow}\left[\pi\left[Y^{m,\alpha}_{k^2,-} \right]\right]=Y^{m,\alpha}_{k^2,-}
\]
and accordingly $Y^{m,\alpha}_{k^2,-}$ is a closed subspace of $(X^{m,\alpha}_{k^2,-})^\perp$. Assume by contradiction that there exists
\[
\tau_-\in (X^{m,\alpha}_{k^2,-})^\perp\setminus Y^{m,\alpha}_{k^2,-}\,.
\]
Then the Hahn-Banach Theorem implies the existence of  $\tilde{\mu}_-\in \left(C^{m,\alpha}(\partial\Omega)\right)'$ such that
\begin{equation}\label{mldc_thm:codasumwv-2}
\langle \tilde{\mu}_-,\tau_-\rangle=1\,,\qquad \langle \tilde{\mu}_-,\eta\rangle=0\qquad
\forall \eta \in Y^{m,\alpha}_{k^2,-}\,.
\end{equation}
Since $\tilde{\mu}_-$ annihilates ${\mathrm{Im}}\, \left(
-\frac{1}{2}I+W_\Omega[\tilde{S}_{n,k;r},\cdot]
\right)_{|C^{m,\alpha}(\partial\Omega)}$, then
\[
\tilde{\mu}_-\in {\mathrm{Ker}}\, \left(
-\frac{1}{2}I+W_\Omega^t[\tilde{S}_{n,k;r},\cdot]\right)\,,
\]
where $W_\Omega^t[\tilde{S}_{n,k;r},\cdot]$ is the transpose operator of $W_\Omega[\tilde{S}_{n,k;r},\cdot]$ with respect to the natural duality pairing
\[
\left(\left(C^{m,\alpha}(\partial\Omega)\right)',C^{m,\alpha}(\partial\Omega)\right)
\]
Since $W_\Omega[\tilde{S}_{n,k;r},\cdot]$ is compact in $C^{m,\alpha}(\partial\Omega)$ (cf.~\textit{e.g.}, \cite[Thm.~7.4]{DoLa17} in case $m=0$ and  \cite[Cor.~9.1]{DoLa17} in case $m\geq 1$), the Fredholm Alternative Theorem implies that
\[
{\mathrm{dim}}\,{\mathrm{Ker}}\, \left(
-\frac{1}{2}I+W_\Omega^t[\tilde{S}_{n,k;r},\cdot]\right)
={\mathrm{dim}}\,{\mathrm{Ker}}\,
\left(
-\frac{1}{2}I+W_\Omega[\tilde{S}_{n,k;r},\cdot]\right)
_{|C^{m,\alpha}(\partial\Omega)}
 \,.
\]
On the other hand, also the operator $W_\Omega^t[\tilde{S}_{n,k;r},\cdot]_{|C^{\max\{1,m\}-1,\alpha}(\partial\Omega)} $ is compact in $C^{\max\{1,m\}-1,\alpha}(\partial\Omega)$ (cf.~\textit{e.g.}, \cite[Cor.~10.1]{DoLa17}) and thus the  Fredholm Alternative Theorem of Wendland \cite{We67}, \cite{We70} in the duality pairing
\begin{eqnarray}\label{mldc_thm:codasumwv-3}
\left(C^{\max\{1,m\}-1,\alpha} (\partial\Omega),C^{m,\alpha} (\partial\Omega)\right) 
\end{eqnarray}
(cf.~\textit{e.g.}, Kress~\cite[Thm.~4.17]{Kr14}) implies that
\begin{eqnarray*}
\lefteqn{
{\mathrm{dim}}\,{\mathrm{Ker}}\, \left(
-\frac{1}{2}I+W_\Omega^t[\tilde{S}_{n,k;r},\cdot]
\right)_{|C^{\max\{1,m\}-1,\alpha}(\partial\Omega)}
}
\\ \nonumber
&&\qquad\qquad\qquad
={\mathrm{dim}}\,{\mathrm{Ker}}\,
\left(
-\frac{1}{2}I+W_\Omega[\tilde{S}_{n,k;r},\cdot]\right)_{|C^{m,\alpha}(\partial\Omega)}
\,,
\end{eqnarray*}
where now $W_\Omega^t[\tilde{S}_{n,k;r},\cdot]_{|C^{\max\{1,m\}-1,\alpha}(\partial\Omega)}$ is the transpose with respect to the duality pairing in (\ref{mldc_thm:codasumwv-3}). Since the canonical embedding ${\mathcal{J}}$ of Lemma \ref{mldc_lem:caincl} maps 
\begin{equation}\label{mldc_thm:codasumwv-4}
{\mathrm{Ker}}\,\left(
-\frac{1}{2}I+W_\Omega^t[\tilde{S}_{n,k;r},\cdot]
\right)
_{|C^{\max\{1,m\}-1,\alpha}(\partial\Omega)}\ \ \text{into}\ \  {\mathrm{Ker}}\, \left(
-\frac{1}{2}I+W_\Omega^t[\tilde{S}_{n,k;r},\cdot]\right)\,,
\end{equation}
the above finite dimensional equalities imply that ${\mathcal{J}}$ restricts a bijection between the spaces in (\ref{mldc_thm:codasumwv-4}). Hence, there exists 
\[
\mu\in {\mathrm{Ker}}\, \left(
-\frac{1}{2}I+W_\Omega^t[\tilde{S}_{n,k;r},\cdot]
\right)
_{|C^{\max\{1,m\}-1,\alpha}(\partial\Omega)}
\]
 such that $\tilde{\mu}_-={\mathcal{J}}[\mu_-]$. Since $\tilde{\mu}_-$ annihilates $Y^{m,\alpha}_{k^2,-}$, then $\tilde{\mu}_-$ annihilates its subspace $V_\Omega[\tilde{S}_{n,k;r},C^{\max\{1,m\}-1,\alpha}(\partial\Omega)] $ and accordingly,
\begin{eqnarray*}
\lefteqn{
0=\langle \tilde{\mu},V_\Omega[\tilde{S}_{n,k;r},\eta]\rangle=\langle {\mathcal{J}}[\mu_-], V_\Omega[\tilde{S}_{n,k;r},\eta]\rangle
}
\\ \nonumber
&& 
=\int_{\partial\Omega}\mu_-  V_\Omega[\tilde{S}_{n,k;r},\eta]\,d\sigma
=\int_{\partial\Omega} V_\Omega[\tilde{S}_{n,k;r},\mu_-]\eta\,d\sigma
\end{eqnarray*}
for all $\eta\in C^{\max\{1,m\}-1,\alpha}(\partial\Omega)$. Hence, $V_\Omega[\tilde{S}_{n,k;r},\mu_-]=0$. Then Proposition \ref{mldc_prop:silamokeri-wt-} (ii) implies that $\mu_-\in X^{m,\alpha}_{k^2,-}$. Since $\tau\in (X^{m,\alpha}_{k^2,-})^\perp$, we must have 
\[
\langle \tilde{\mu}_-,\tau_-\rangle=\int_{\partial\Omega}\mu_-\tau_-\,d\sigma=0\,,
\]
a contradiction (see (\ref{mldc_thm:codasumwv-2})). Hence, $(X^{m,\alpha}_{k^2,-})^\perp=Y^{m,\alpha}_{k^2,-}$ and the proof is complete.\hfill  $\Box$ 

\vspace{\baselineskip}

Then by combining the necessary condition of Theorem \ref{mldc_thm:nexextdirhe}, Theorem \ref{mldc_thm:codasumwv-} and the jump formulas for the  acoustic  single and the double   layer potential, we deduce the validity of the following statement.

\begin{theorem}[of existence for the exterior Dirichlet problem]\label{mldc_thm:exextdirhe}
  Let $m\in {\mathbb{N}}$,   $\alpha\in]0,1[$. Let $\Omega$ be a bounded open subset of ${\mathbb{R}}^{n}$ of class $C^{\max\{1,m\},\alpha}$. 
  
  Let $k\in {\mathbb{C}}\setminus ]-\infty,0]$, ${\mathrm{Im}}\,k\geq 0$.   Let $g\in C^{m,\alpha}(\partial\Omega)$. Then  the exterior Dirichlet problem 
\begin{equation}
\label{mldc_thm:exextdirhe1}
\left\{
\begin{array}{ll}
\Delta u+k^2 u=0 &{\text{in}}\ \Omega^-\,,
\\
u=g
&{\text{on}}\ \partial\Omega\,,
\\
u\ \text{satisfies\ the\ outgoing}\  k-\text{radiation\ condition}\,,
\end{array}
\right.
\end{equation}
has a solution $u\in C^{m,\alpha}_{ {\mathrm{loc}} }(\overline{\Omega^-})$ if and only if the compatibility condition
\begin{equation}\label{mldc_thm:exextdirhe2}
\int_{\partial\Omega}g\frac{\partial v}{\partial\nu_\Omega^-}\,d\sigma=0  \,,
\end{equation}
is satisfied for all $v\in C^{\max\{1,m\},\alpha}_{ {\mathrm{loc}} }(\overline{\Omega^-})$  such that
\begin{equation}\label{mldc_thm:exextdirhe3}
\left\{
\begin{array}{ll}
 \Delta v+k^2 v=0 &\text{in}\ \Omega^-\,,
 \\
 v=0		&\text{on}\ \partial\Omega\,,
 \\
 v\ \text{satisfies\ the\ outgoing}\  k-\text{radiation\ condition}\,.
\end{array}
\right.
\end{equation}
\end{theorem}
{\bf Proof.} The necessity follows by Theorem \ref{mldc_thm:nexextdirhe}. If $g$ satisfies condition (\ref{mldc_thm:exextdirhe2}), then $g\in (X^{m,\alpha}_{k^2,-})^\perp$ and Theorem \ref{mldc_thm:codasumwv-}
implies the existence of $\phi\in C^{m,\alpha}(\partial\Omega)$ and of $\psi\in C^{\max\{1,m\}-1,\alpha}(\partial\Omega)$ such that
\[
g=-\frac{1}{2}\phi+W_\Omega[\tilde{S}_{n,k;r},\phi]+V_\Omega[\tilde{S}_{n,k;r},\psi]\,.
\]
Thus if we set
\[
u=w_\Omega^-[\tilde{S}_{n,k;r},\phi] +v_\Omega^-[\tilde{S}_{n,k;r},\psi]\,,
\]
then $u\in  C^{m,\alpha}_{{\mathrm{loc}} }(\overline{\Omega^-})$ by known  regularity results for the  acoustic  single layer potential (cf.~\textit{e.g.}, \cite[Thm.~7.1]{DoLa17}) and by known regularity results for the  acoustic  double   layer   potential (cf.~\textit{e.g.}, 
 \cite[Thm.~7.3]{DoLa17} in case $m\geq 1$, 
 \cite[Thm.~5.1]{La25a} in case $m=0$). Also, the known  jump formulas for the acoustic single and double  layer potential (cf.~\textit{e.g.},  \cite[Thms.~7.1, 7.3]{DoLa17})  imply that $u_{|\partial\Omega}=g$.\hfill  $\Box$ 

\vspace{\baselineskip}

\section{A nonvariational form of the interior Neumann problem for the Helmholtz equation}\label{mldc_sec:indneh}

We now consider the interior  Neumann problem for the Helmholtz equation in the case in which the Neumann datum is in the space $V^{m-1,\alpha}(\partial\Omega)$  (see (\ref{mldc_eq:vm-1a}). Nonvariational in case $m=0$).  For the corresponding classical result in case $\varkappa^-=0$, we refer to  Colton and Kress \cite[Thm.~3.20]{CoKr92}.  (For the definition of $\varkappa^-$, see right after (\ref{mldc_eq:exto})).  For case $m=0$ and in the specific case in which solutions can be represented in terms of a single acoustic layer, with some extra assumptions on the bounded connected components of the exterior $\Omega^-$, we refer to
\cite[Thm.~10.1]{La25b}.  
  We first  introduce the following necessary condition for the existence (see also \cite[Thm.~10.1]{La25b}, which we follow).
\begin{theorem} \label{mldc_thm:nexintneuhe}
Let $m\in {\mathbb{N}}$,   $\alpha\in]0,1[$. Let $\Omega$ be a bounded open subset of ${\mathbb{R}}^{n}$ of class $C^{\max\{1,m\},\alpha}$. Let  $\lambda\in {\mathbb{C}}$.  Let $g\in V^{m-1,\alpha}(\partial\Omega)$. If the interior Neumann problem
 \begin{equation}
\label{mldc_thm:nexintneuhe1}
\left\{
\begin{array}{ll}
\Delta u+\lambda u=0 &{\text{in}}\ \Omega\,,
\\
\frac{\partial u}{\partial\nu_\Omega}=g
&{\text{on}}\ \partial\Omega\,
 \end{array}
\right.
\end{equation}
has   a  solution $u\in C^{0,\alpha}(\overline{\Omega})_\Delta$, then
\begin{equation}
\label{mldc_thm:nexintneuhe2}
\langle g,v_{|\partial\Omega}\rangle =0
\end{equation}
for all $v\in C^{\max\{1,m\},\alpha}(\overline{\Omega})$ such that
\begin{equation}\label{mldc_thm:nexintneuhe3}
\Delta v+\lambda v=0 \qquad\text{in}\ \Omega\,,\qquad 
\frac{\partial v}{\partial\nu_\Omega}=0\qquad\text{on}\ \partial\Omega\,.
\end{equation}
\end{theorem}
{\bf Proof.}  Assume that a solution $u\in C^{0,\alpha}(\overline{\Omega})_\Delta$ exists and that $v\in C^{\max\{1,m\},\alpha}(\overline{\Omega})$ satisfies (\ref{mldc_thm:nexintneuhe3}). 
Then Remark \ref{mldc_rem:helpo} implies that $v\in C^{1,\alpha}(\overline{\Omega})_\Delta$  and  
 the second Green Identity in distributional form of \cite[Thm.~8.1]{La25b}
 and equality (\ref{mldc_prop:nschext3}) imply that
 \begin{eqnarray*}
\lefteqn{
\langle g,v\rangle =\langle  \frac{\partial u}{\partial\nu_\Omega}, v_{|\partial\Omega}\rangle 
}
\\ \nonumber
&&\qquad\qquad
=\int_{\partial\Omega}u \frac{\partial v}{\partial\nu_\Omega}  \,d\sigma
-
\int_\Omega u(\Delta v+\lambda v)\,dx+\langle E^\sharp[\Delta u+\lambda u],v\rangle =0\,.
\end{eqnarray*}\hfill  $\Box$ 

\vspace{\baselineskip}

By Theorem \ref{mldc_thm:nexintneuhe}, a necessary condition for the solvability of the Neumann problem is that condition (\ref{mldc_thm:nexintneuhe2}) holds true. Thus we now introduce the following notation. If $m\in{\mathbb{N}}$, $\alpha\in]0,1[$,  $\Omega$ is a bounded open subset of ${\mathbb{R}}^{n}$ of class $C^{\max\{1,m\},\alpha}$, $\lambda\in {\mathbb{C}}$, then we set
\begin{equation}\label{mldc_eq:smal}
S^{m,\alpha}_\lambda\equiv\left\{
 v_{|\partial\Omega}:\, v\in C^{\max\{1,m\},\alpha}(\overline{\Omega})\,,\ 
\Delta v+\lambda v=0 \ \text{in}\ \Omega\,,\  \frac{\partial v}{\partial\nu_\Omega}=0\ \text{on}\ \partial\Omega\right\}
\end{equation}
and
\begin{equation}\label{mldc_eq:smal1}
(S^{m,\alpha}_\lambda)^\perp\equiv
\left\{
g\in \left(C^{\max\{1,m\},\alpha}(\partial\Omega)\right)':\, \langle g, f \rangle=0\ \forall f \in S^{m,\alpha}_\lambda
\right\}\,.
\end{equation}

\begin{proposition}\label{mldc_prop:-1/2+w,dnw+}
 Let $m\in{\mathbb{N}}$, $\alpha\in]0,1[$. Let $\Omega$ be a bounded open subset of ${\mathbb{R}}^{n}$ of class $C^{\max\{1,m\},\alpha}$. Let $\lambda\in {\mathbb{C}}$. Let $S_{n,\lambda}$ be a fundamental solution of $\Delta+\lambda$. Then the following statements hold true. 
 \begin{enumerate}
\item[(i)] $ -\frac{1}{2}\phi+W_\Omega^t[S_{n,\lambda},\phi]\in (S^{m,\alpha}_\lambda)^\perp$ for all $\phi\in  \left(C^{\max\{1,m\},\alpha}(\partial\Omega)\right)'$.
 \item[(ii)]  $\frac{\partial}{\partial\nu_\Omega}w_\Omega^+[S_{n,\lambda},\psi]\in 
   (S^{m,\alpha}_\lambda)^\perp\cap V^{m-1,\alpha}(\partial\Omega)$ for all $\psi\in C^{m,\alpha}(\partial\Omega)$.
\end{enumerate}
\end{proposition}
{\bf Proof.} (i) If $\phi\in  \left(C^{\max\{1,m\},\alpha}(\partial\Omega)\right)'$ and $f\in S^{m,\alpha}_\lambda$, then
\[
\langle -\frac{1}{2}\phi+W_\Omega^t[S_{n,\lambda},\phi],f \rangle=
\langle\phi,-\frac{1}{2}f +W_\Omega[S_{n,\lambda},f]\rangle\,.
\]
Since Corollary 9.9 (i) of \cite{La25b} ensures that $-\frac{1}{2}f +W_\Omega[S_{n,\lambda},f]=0$, we conclude that the right hand side equals $0$. 

(ii) By known regularity results for the  acoustic  double layer   potential (cf.~\textit{e.g.}, 
 \cite[Thm.~7.3]{DoLa17} in case $m\geq 1$, 
 \cite[Thm.~5.1]{La25a} in case $m=0$), we have $w_\Omega^+[S_{n,\lambda},\psi]\in 
 C^{m,\alpha}(\overline{\Omega})$ and
  accordingly $\frac{\partial}{\partial\nu_\Omega}w_\Omega^+[S_{n,\lambda},\psi]$ belongs to $ V^{m-1,\alpha}(\partial\Omega)$ (in case $m=0$, see also Remark \ref{mldc_rem:helpo}
 and Proposition \ref{mldc_prop:ricodnu}). Moreover, 
the membership of $\frac{\partial}{\partial\nu_\Omega}w_\Omega^+[S_{n,\lambda},\psi]$ in $ 
   (S^{m,\alpha}_\lambda)^\perp$ follows  by  Remark \ref{mldc_rem:helpo} and Theorem \ref{mldc_thm:nexintneuhe} with $g=\frac{\partial}{\partial\nu_\Omega}w_\Omega^+[S_{n,\lambda},\psi]$.\hfill  $\Box$ 

\vspace{\baselineskip}

\begin{proposition}\label{mldc_prop:isomquotn}
Let $m\in{\mathbb{N}}$, $\alpha\in]0,1[$. Let $\Omega$ be a bounded open subset of ${\mathbb{R}}^{n}$ of class $C^{\max\{1,m\},\alpha}$.  
 Let $\lambda\in {\mathbb{C}}$. Let $S_{n,\lambda}$ be a fundamental solution of $\Delta+\lambda$. Then 
 \begin{equation}\label{mldc_prop:isomquotn1}
{\mathrm{Im}}\, \left(
-\frac{1}{2}I+W_\Omega^t[S_{n,\lambda},\cdot] \right)
=\left({\mathrm{Ker}}\,\left(-\frac{1}{2}I+W_\Omega[S_{n,\lambda},\cdot]\right)_{|C^{\max\{1,m\},\alpha}(\partial\Omega)}\right)^\perp
\end{equation}
 and ${\mathrm{Im}}\, \left(
-\frac{1}{2}I+W_\Omega^t[S_{n,\lambda},\cdot] \right)$ is a weakly$^\ast$ closed subspace of
$(S^{m,\alpha}_\lambda)^\perp$ of finite codimension. 
Here,
$W_\Omega^t[\tilde{S}_{n,k;r},\cdot]$ is the transpose to the operator $W_\Omega[\tilde{S}_{n,k;r},\cdot]$ in $C^{\max\{1,m\},\alpha}(\partial\Omega)$ and the annihilator is taken with respect to the natural duality pairing (\ref{mldc_rem:wtnotation1}).
\end{proposition}
{\bf Proof.} Since $W_\Omega[S_{n,\lambda},\cdot]$ is compact in  $C^{\max\{1,m\},\alpha}(\partial\Omega)$  (cf.~\textit{e.g.},   \cite[Cor.~9.1]{DoLa17}), then the Fredholm Alternative Theorem  in the natural duality pairing
(\ref{mldc_rem:wtnotation1}) implies the validity of   equality (\ref{mldc_prop:isomquotn1}) and that 
${\mathrm{Im}}\, \left(
-\frac{1}{2}I+W_\Omega^t[S_{n,\lambda},\cdot] \right)$ has finite codimension in $\left(C^{\max\{1,m\},\alpha}(\overline{\Omega})\right)'$.  Since 
\[
\left({\mathrm{Ker}}\,\left(-\frac{1}{2}I+W_\Omega[S_{n,\lambda},\cdot]\right)_{|C^{\max\{1,m\},\alpha}(\partial\Omega)}\right)^\perp
\]
 is weakly$^\ast$ closed in $\left(C^{\max\{1,m\},\alpha}(\partial\Omega)\right)'$, then the space  ${\mathrm{Im}}\, \left(
-\frac{1}{2}I+W_\Omega^t[S_{n,\lambda},\cdot] \right)$ is also weakly$^\ast$ closed in    $\left(C^{\max\{1,m\},\alpha}(\partial\Omega)\right)'$. 

By Proposition \ref{mldc_prop:-1/2+w,dnw+} (i), we know that the space ${\mathrm{Im}}\, \left(
-\frac{1}{2}I+W_\Omega^t[S_{n,\lambda},\cdot] \right)$ is contained in $(S^{m,\alpha}_\lambda)^\perp$, which is also weakly$^\ast$ closed in $\left(C^{\max\{1,m\},\alpha}(\partial\Omega)\right)'$. Hence, the proof is complete.\hfill  $\Box$ 

\vspace{\baselineskip}

Then we can prove the following statement that says that the space of compatible data can be written as an algebraic sum of normal derivatives of an  acoustic  single   layer potential and of normal derivatives of an  acoustic  double   layer potential. 
\begin{theorem}\label{mldc_thm:codasumnwv}
Let $m\in{\mathbb{N}}$, $\alpha\in]0,1[$. Let $\Omega$ be a bounded open subset of ${\mathbb{R}}^{n}$ of class $C^{\max\{1,m\},\alpha}$.  Let $k\in {\mathbb{C}}\setminus ]-\infty,0]$, ${\mathrm{Im}}\,k\geq 0$. Then
\begin{equation}\label{mldc_thm:codasumnwv1}
(S^{m,\alpha}_{k^2})^\perp=
{\mathrm{Im}}\, \left(
-\frac{1}{2}I+W_\Omega^t[\tilde{S}_{n,k;r},\cdot] \right)
+\frac{\partial}{\partial\nu_\Omega}w_\Omega^+\left[\tilde{S}_{n,k;r},
C^{m,\alpha}(\partial\Omega)
\right]\,,
\end{equation}
where  $W_\Omega^t[\tilde{S}_{n,k;r},\cdot]$ is the transpose to the operator $W_\Omega[\tilde{S}_{n,k;r},\cdot]$ in $C^{\max\{1,m\},\alpha}(\partial\Omega)$.
\end{theorem}
{\bf Proof.} By Proposition \ref{mldc_prop:isomquotn}, ${\mathrm{Im}}\, \left(
-\frac{1}{2}I+W_\Omega^t[\tilde{S}_{n,k;r},\cdot] \right)$ is a weakly$^\ast$ closed subspace of
$(S^{m,\alpha}_{k^2})^\perp$. Let $\pi$ be the canonical projection of  the topological vector space 
$(S^{m,\alpha}_{k^2})^\perp$ with the weak$^\ast$ topology onto the finite dimensional quotient space
\[
Q^{m,\alpha}_{k^2}\equiv 
\frac{	
(S^{m,\alpha}_{k^2})^\perp
}{
{\mathrm{Im}}\, \left(
-\frac{1}{2}I+W_\Omega^t[\tilde{S}_{n,k;r},\cdot]
\right)
} \,.
\]
By Proposition \ref{mldc_prop:-1/2+w,dnw+} (ii),
\[
\frac{\partial}{\partial\nu_\Omega}w_\Omega^+\left[\tilde{S}_{n,k;r},
C^{m,\alpha}(\partial\Omega)
\right]
\leq (S^{m,\alpha}_{k^2})^\perp
\]
 and thus  the $\pi$-image $\pi\left[T^{m,\alpha}_{k^2}\right]$
of
\[
T^{m,\alpha}_{k^2}
\equiv
{\mathrm{Im}}\, \left(
-\frac{1}{2}I+W_\Omega^t[\tilde{S}_{n,k;r},\cdot]
\right)
+\frac{\partial}{\partial\nu_\Omega}w_\Omega^+\left[\tilde{S}_{n,k;r},
C^{m,\alpha}(\partial\Omega)
\right]
\]
is a   subspace of the finite dimensional Hausdorff topological     quotient vector space $ 
Q^{m,\alpha}_{k^2}$ and is accordingly closed in $Q^{m,\alpha}_{k^2}$. Since the canonical projection $\pi$ is continuous, then the preimage
\[
\pi^{\leftarrow}\left[\pi\left[T^{m,\alpha}_{k^2} \right]\right]
\]
is weakly$^\ast$ closed in $(S^{m,\alpha}_{k^2} )^\perp$.  Now a simple computation based on the definition of $\pi$  shows that
\[
\pi^{\leftarrow}\left[\pi\left[T^{m,\alpha}_{k^2} \right]\right]=T^{m,\alpha}_{k^2}
\]
and accordingly $T^{m,\alpha}_{k^2}$ is a weakly$^\ast$ closed subspace of $(S^{m,\alpha}_{k^2})^\perp$. Assume by contradiction that $T^{m,\alpha}_{k^2}\subsetneq (S^{m,\alpha}_{k^2})^\perp$. Then   there exists
\[
\tau\in (S^{m,\alpha}_{k^2})^\perp\setminus T^{m,\alpha}_{k^2}\,.
\]
By the Hahn-Banach Theorem in the dual space
\[
\left(C^{\max\{1,m\},\alpha}(\partial\Omega)\right)'_s
\]
with the weak$^\ast$ topology, there exists
\[
\tilde{\mu}\in\left(\left(C^{\max\{1,m\},\alpha}(\partial\Omega)\right)'_s\right)'
\]
such that
\begin{equation}\label{mldc_thm:codasumnwv2}
\langle \tilde{\mu},\tau\rangle=1\,,\qquad \langle \tilde{\mu},x'\rangle=0\qquad
\forall x' \in T^{m,\alpha}_{k^2}\,.
\end{equation}
Now the canonical isometry of $C^{\max\{1,m\},\alpha}(\partial\Omega)$ is an isomorphism of the space $C^{\max\{1,m\},\alpha}(\partial\Omega)$ onto the bi-dual
\[
\left(\left(C^{\max\{1,m\},\alpha}(\partial\Omega)\right)'_s\right)'\,,
\]
which is a  subspace of the strong bidual $\left(C^{\max\{1,m\},\alpha}(\partial\Omega)\right)''$
(cf.~\textit{e.g.}, Brezis \cite[Prop.~3.14, p.~64]{Br11}). Hence, there exists $\mu\in C^{\max\{1,m\},\alpha}(\partial\Omega)$ such that
\[
\langle \tilde{\mu},x'\rangle=\langle x',\mu\rangle\qquad\forall x'\in 
\left(C^{\max\{1,m\},\alpha}(\partial\Omega)\right)'
\]
and we must have
\begin{equation}\label{mldc_thm:codasumnwv3}
\langle\tau,\mu \rangle=1\,,\qquad \langle x',\mu\rangle=0\qquad
\forall x' \in T^{m,\alpha}_{k^2}\,.
\end{equation}
In particular, $\mu$ annihilates the subspace ${\mathrm{Im}}\, \left(
-\frac{1}{2}I+W_\Omega^t[\tilde{S}_{n,k;r},\cdot]
\right)$ of $T^{m,\alpha}_{k^2}$ and accordingly Proposition  \ref{mldc_prop:isomquotn}    implies that 
\begin{eqnarray}\label{mldc_thm:codasumnwv4}
\lefteqn{\mu\in
\left(
{\mathrm{Ker}}\,\left(-\frac{1}{2}I+W_\Omega[\tilde{S}_{n,k;r},\cdot]\right)_{|C^{\max\{1,m\},\alpha}(\partial\Omega)}\right)^{\perp\perp}
}
\\ \nonumber
&&\qquad\qquad\qquad\qquad
=
{\mathrm{Ker}}\,\left(-\frac{1}{2}I+W_\Omega[\tilde{S}_{n,k;r},\cdot]\right)_{|C^{\max\{1,m\},\alpha}(\partial\Omega)}\,,
\end{eqnarray}
(cf.~\textit{e.g.}, Rudin, \cite[Thm.~4.7, Chap. IV]{Ru91}).
Hence,  the jump formulas for the acoustic double  layer  potential imply that 
\begin{equation}\label{mldc_thm:codasumnwv4a}
w_\Omega^-\left[\tilde{S}_{n,k;r},\mu\right]=0\qquad\text{on}\ \partial\Omega
\end{equation}
and accordingly that 
$w_\Omega^-\left[\tilde{S}_{n,k;r},\mu\right]$ solves the Helmholtz equation in $\Omega^-$ that corresponds to $k^2$ with zero Dirichlet boundary condition (cf.~\textit{e.g.},  \cite[Thm.~7.3 (i)]{DoLa17}). By known regularity results for the  acoustic  double layer   potential (cf.~\textit{e.g.}, 
 \cite[Thm.~7.3]{DoLa17}), we have $w_\Omega^-[\tilde{S}_{n,k;r},\mu]\in 
 C^{\max\{1,m\},\alpha}_{{\mathrm{loc}}}(\overline{\Omega^-})$. Since $w_\Omega^-\left[\tilde{S}_{n,k;r},\mu\right]$ solves the Helmholtz equation in $\Omega^-$ that corresponds to $k^2$ with zero Dirichlet boundary condition and satisfies the outgoing $k$-radiation condition (cf.~\textit{e.g.}, \cite[Thm.~6.20 (ii)]{La25b}), Corollary \ref{mldc_corol:eipedd} (i) implies that
\begin{equation}\label{mldc_thm:codasumnwv4b}
\frac{\partial}{\partial\nu_{\Omega^-}}w_\Omega^-\left[\tilde{S}_{n,k;r},
\mu
\right]\in {\mathrm{Ker}}\,\left(-\frac{1}{2}I+W_\Omega^t[\tilde{S}_{n,k;r},\cdot]\right)_{|C^{\max\{1,m\}-1,\alpha}(\partial\Omega)}\,.
\end{equation}
Then the jump formula for the   normal derivative of the   acoustic  double  layer potential implies that
\begin{eqnarray}\label{mldc_thm:codasumnwv5}
\lefteqn{-\frac{\partial}{\partial\nu_{\Omega}}w_\Omega^+\left[\tilde{S}_{n,k;r},
\mu
\right]
}
\\ \nonumber
&&\qquad 
=
\frac{\partial}{\partial\nu_{\Omega^-}}w_\Omega^-\left[\tilde{S}_{n,k;r},
\mu
\right]\in {\mathrm{Ker}}\,\left(-\frac{1}{2}I+W_\Omega^t[\tilde{S}_{n,k;r},\cdot]\right)_{|C^{\max\{1,m\}-1,\alpha}(\partial\Omega)}
\end{eqnarray}
(cf.~\textit{e.g.}, \cite[Thm.~6.5]{La25a}). Moreover, the jump formula for the  acoustic  double layer potential implies that
\begin{equation}\label{mldc_thm:codasumnwv6}
\mu=w_\Omega^+\left[\tilde{S}_{n,k;r},\mu\right]-w_\Omega^-\left[\tilde{S}_{n,k;r},\mu\right]
=w_\Omega^+\left[\tilde{S}_{n,k;r},\mu\right]\qquad\text{on}\ \partial\Omega\,,
\end{equation}
(cf.~(\ref{mldc_thm:dlay1a}), (\ref{mldc_thm:codasumnwv4a})). Since $\mu$ annihilates the subspace $\frac{\partial}{\partial\nu_\Omega}w_\Omega^+\left[\tilde{S}_{n,k;r},
C^{m,\alpha}(\partial\Omega)
\right]$ of $T^{m,\alpha}_{k^2}$, we have
\[
\langle
\frac{\partial}{\partial\nu_\Omega}w_\Omega^+\left[\tilde{S}_{n,k;r},\eta\right]
,\mu \rangle=0
\qquad\forall \eta\in C^{m,\alpha}(\partial\Omega)\,.
\]
Hence, Remark \ref{mldc_rem:helpo} (i) and the second Green Identity in distributional form of \cite[Thm.~8.1]{La25b} imply that
\begin{eqnarray*}
\lefteqn{
0=\langle
\frac{\partial}{\partial\nu_\Omega}w_\Omega^+\left[\tilde{S}_{n,k;r},\eta\right]
,\mu \rangle
}
\\ \nonumber
&&\qquad
=\langle
\frac{\partial}{\partial\nu_\Omega}w_\Omega^+\left[\tilde{S}_{n,k;r},\eta\right]
,w_\Omega^+\left[\tilde{S}_{n,k;r},\mu\right]_{|\partial\Omega} \rangle
\\ \nonumber
&&\qquad
=\int_{\partial\Omega}w_\Omega^+\left[\tilde{S}_{n,k;r},\eta\right]
\frac{\partial}{\partial\nu_\Omega}w_\Omega^+\left[\tilde{S}_{n,k;r},\mu\right]\,d\sigma
\\ \nonumber
&&\qquad\quad
+\langle E^\sharp\left[\Delta \left(w_\Omega^+\left[\tilde{S}_{n,k;r},\eta\right]\right)+k^2
w_\Omega^+\left[\tilde{S}_{n,k;r},\eta\right]
\right],
w_\Omega^+\left[\tilde{S}_{n,k;r},\mu\right]\rangle
\\ \nonumber
&&\qquad\quad
-\int_\Omega w_\Omega^+\left[\tilde{S}_{n,k;r},\eta\right]
\left\{\Delta \left(w_\Omega^+\left[\tilde{S}_{n,k;r},\mu\right]\right)
+k^2w_\Omega^+\left[\tilde{S}_{n,k;r},\mu\right]
\right\}\,dx
\\ \nonumber
&&\qquad
=\int_{\partial\Omega}w_\Omega^+\left[\tilde{S}_{n,k;r},\eta\right]
\frac{\partial}{\partial\nu_\Omega}w_\Omega^+\left[\tilde{S}_{n,k;r},\mu\right]\,d\sigma
\qquad\forall \eta\in C^{m,\alpha}(\partial\Omega)\,.
\end{eqnarray*}
In particular, the same equality holds for all $\eta\in  C^{\max\{m,1\},\alpha}(\partial\Omega)$. Hence, the jump formula (\ref{mldc_thm:dlay1a}) for the acoustic double layer potential   implies that
\begin{eqnarray}\label{mldc_thm:codasumnwv7}
\lefteqn{\frac{\partial}{\partial\nu_\Omega}w_\Omega^+\left[\tilde{S}_{n,k;r},\mu\right]
\in 
\left({\mathrm{Im}}\,
\left(\frac{1}{2}I+W_\Omega[\tilde{S}_{n,k;r},\cdot]\right)_{|C^{\max\{m,1\},\alpha}(\partial\Omega)}\right)^\perp
}
\\ \nonumber
&&\qquad\qquad\qquad\qquad\qquad\qquad\qquad\qquad\qquad\qquad
={\mathrm{Ker}}\,\left(\frac{1}{2}I+W_\Omega^t[\tilde{S}_{n,k;r},\cdot]\right)\,,
 \end{eqnarray}
 where the annihilator is taken with respect to the natural duality 
 (\ref{mldc_rem:wtnotation1}) of the space $C^{\max\{m,1\},\alpha}(\partial\Omega)$ (cf.~\textit{e.g.}, Rudin, \cite[Thm.~4.12, Chap. IV]{Ru91}). 
Then by combining  (\ref{mldc_thm:codasumnwv5}) and the membership of (\ref{mldc_thm:codasumnwv7}), we deduce that
\[
\frac{\partial}{\partial\nu_\Omega}w_\Omega^+\left[\tilde{S}_{n,k;r},\mu\right]
=0\,.
\]
 Since $w_\Omega^+\left[\tilde{S}_{n,k;r},\mu\right]$ solves the Helmholtz equation in $\Omega$ with zero
Neumann boundary conditions, then we have
\[
\mu=w_\Omega^+\left[\tilde{S}_{n,k;r},\mu\right]_{|\partial\Omega}\in S^{m,\alpha}_{k^2}\,.
\]
Since $\tau\in (S^{m,\alpha}_{k^2})^\perp$, we have $\langle\tau,\mu\rangle=0$, a contradiction (see (\ref{mldc_thm:codasumnwv3})). Hence, $T^{m,\alpha}_{k^2} = (S^{m,\alpha}_{k^2})^\perp$ and thus the proof is complete.\hfill  $\Box$ 

\vspace{\baselineskip}

By exploiting the direct sum of Theorem \ref{mldc_thm:codasumnwv}, one could solve the Neumann problem with a distributional datum in $(S^{m,\alpha}_{k^2})^\perp$. We now show that if the datum is also in $V^{m-1,\alpha}(\partial\Omega)$, then the corresponding solutions belong  to $C^{m,\alpha}(\overline{\Omega})$.
\begin{theorem}\label{mldc_thm:pm1/2wtreg}
 Let  $m\in{\mathbb{N}}$, $\alpha\in]0,1[$. Let $\Omega$ be a bounded open subset of ${\mathbb{R}}^{n}$ of class $C^{\max\{1,m\},\alpha}$.  Let $\lambda\in {\mathbb{C}}$. Let $S_{n,\lambda}$ be a fundamental solution of $\Delta+\lambda$.  Let $\tau\in V^{m-1,\alpha}(\partial\Omega)$. If $\mu_\pm\in 
\left(C^{\max\{1,m\},\alpha}(\partial\Omega)\right)'$ and
\[
\pm\frac{1}{2}\mu_\pm +W_\Omega^t[S_{n,\lambda},\mu_\pm]=\tau\,,
\] 
then $\mu_\pm\in V^{m-1,\alpha}(\partial\Omega)$,  where we understand that in the symbol $\pm$ we can take either only the `plus' sign $+$  or only the `minus' sign  $-$ (cf.~(\ref{mldc_eq:vm-1a}), Remark \ref{mldc_rem:wtnotation}). 
\end{theorem}
{\bf Proof.} We first note that 
\begin{eqnarray*} 
\lefteqn{
\tau\in {\mathrm{Im}}\,\left(\pm\frac{1}{2}I +W_\Omega^t[S_{n,\lambda},\cdot]\right)
\subseteq \left(
{\mathrm{Im}}\,\left(\pm\frac{1}{2}I +W_\Omega^t[S_{n,\lambda},\cdot]\right)\right)^{\perp\perp}
}
\\ \nonumber
&&\qquad\qquad\qquad\qquad 
=\left({\mathrm{Ker}}\,\left(\pm\frac{1}{2}I +W_\Omega[S_{n,\lambda},\cdot]\right)_{|C^{\max\{1,m\},\alpha}(\partial\Omega)}\right)^\perp\,,
\end{eqnarray*}
where the annihilators are taken with respect to the natural duality 
 (\ref{mldc_rem:wtnotation1}) of $C^{\max\{m,1\},\alpha}(\partial\Omega)$ 
(cf.~\textit{e.g.}, Rudin, \cite[Thms.~4.7, 4.12, Chap. IV]{Ru91}). Hence, \[
\langle\tau,f\rangle=0\qquad\forall f \in{\mathrm{Ker}}\,\left(\pm\frac{1}{2}I +W_\Omega[S_{n,\lambda},\cdot]\right)_{|C^{\max\{1,m\},\alpha}(\partial\Omega)}\,.
\]
Next we note that $W_\Omega[S_{n,\lambda},\cdot]$ is compact in $C^{\max\{1,m\},\alpha}(\partial\Omega)$ (cf.~\textit{e.g.}, \cite[Cor.~9.1]{DoLa17}) and that 
$W_\Omega^t[S_{n,\lambda},\cdot]$ is compact in $V^{m-1,\alpha}(\partial\Omega)$ (cf.~\textit{e.g.}, \cite[Cor.~10.1]{DoLa17} in case $m\geq 1$ and \cite[Cor.~8.5]{La25a} in case $m=0$). Then the  Fredholm Alternative Theorem of Wendland \cite{We67}, \cite{We70}
in the duality pairing of (\ref{mldc_rem:wtnotation2}) implies that there exists
\[
\mu^\sharp_\pm\in V^{m-1,\alpha}(\partial\Omega)
\]
such that
\[
\pm\frac{1}{2}\mu^\sharp_\pm +W_\Omega^t[S_{n,\lambda},\mu^\sharp_\pm]=\tau
\]
and that
\begin{eqnarray*}
\lefteqn{{\mathrm{dim}}\,
{\mathrm{Ker}}\,\left(\pm\frac{1}{2}I +W_\Omega[S_{n,\lambda},\cdot]\right)_{|C^{\max\{1,m\},\alpha}(\partial\Omega)}
}
\\ \nonumber
&&\qquad\qquad\qquad\qquad={\mathrm{dim}}\,
{\mathrm{Ker}}\,\left(\pm\frac{1}{2}I +W_\Omega^t[S_{n,\lambda},\cdot]\right)_{|V^{m-1,\alpha}(\partial\Omega)}\,,
\end{eqnarray*}
 (cf.~\textit{e.g.}, Kress~\cite[Thm.~4.17]{Kr14}).
By  the Fredholm Alternative in the natural duality pairing (\ref{mldc_rem:wtnotation1}), we have
\begin{eqnarray*} 
\lefteqn{
{\mathrm{dim}}\,
{\mathrm{Ker}}\,\left(\pm\frac{1}{2}I +W_\Omega[S_{n,\lambda},\cdot]\right)_{|C^{\max\{1,m\},\alpha}(\partial\Omega)}
}
\\ \nonumber
&&\qquad\qquad\qquad 
={\mathrm{dim}}\,
{\mathrm{Ker}}\,\left(\pm\frac{1}{2}I +W_\Omega^t[S_{n,\lambda},\cdot]\right)  \,.
\end{eqnarray*}
Hence,
\[
{\mathrm{dim}}\,
{\mathrm{Ker}}\,\left(\pm\frac{1}{2}I +W_\Omega^t[S_{n,\lambda},\cdot]\right)_{|V^{m-1,\alpha}(\partial\Omega)}=
{\mathrm{dim}}\,
{\mathrm{Ker}}\,\left(\pm\frac{1}{2}I +W_\Omega^t[S_{n,\lambda},\cdot]\right)\,.
\]
Since $V^{m-1,\alpha}(\partial\Omega)\leq \left(C^{\max\{1,m\},\alpha}(\partial\Omega)\right)'$, we conclude that
\[
{\mathrm{Ker}}\,\left(\pm\frac{1}{2}I +W_\Omega^t[S_{n,\lambda},\cdot]\right)_{|V^{m-1,\alpha}(\partial\Omega)}={\mathrm{Ker}}\,\left(\pm\frac{1}{2}I +W_\Omega^t[S_{n,\lambda},\cdot]\right)\,.
\]
Since $\mu_\pm-\mu^\sharp_\pm\in {\mathrm{Ker}}\,\left(\pm\frac{1}{2}I +W_\Omega^t[S_{n,\lambda},\cdot]\right)$, we have 
\[
\mu_\pm-\mu^\sharp_\pm\in  {\mathrm{Ker}}\,\left(\pm\frac{1}{2}I +W_\Omega^t[S_{n,\lambda},\cdot]\right)_{|V^{m-1,\alpha}(\partial\Omega)}\,.
\]
Then the membership of $\mu^\sharp_\pm$ in $V^{m-1,\alpha}(\partial\Omega)$ implies that
$\mu_\pm$ belongs to $V^{m-1,\alpha}(\partial\Omega)$ and thus the proof is complete.\hfill  $\Box$ 

\vspace{\baselineskip}

Then by combining Proposition \ref{mldc_prop:-1/2+w,dnw+} (ii), Theorem \ref{mldc_thm:codasumnwv} and the regularity Theorem \ref{mldc_thm:pm1/2wtreg}, we deduce the validity of the following statement.
\begin{theorem}\label{mldc_thm:codamasumnwv}
Let $m\in{\mathbb{N}}$, $\alpha\in]0,1[$. Let $\Omega$ be a bounded open subset of ${\mathbb{R}}^{n}$ of class $C^{\max\{1,m\},\alpha}$.  Let $k\in {\mathbb{C}}\setminus ]-\infty,0]$, ${\mathrm{Im}}\,k\geq 0$. Then
\begin{eqnarray}\label{mldc_thm:codamasumnwv1}
\lefteqn{(S^{m,\alpha}_{k^2})^\perp
\cap V^{m-1,\alpha}(\partial\Omega)
}
\\ \nonumber
&&\qquad=
{\mathrm{Im}}\, \left(
-\frac{1}{2}I+W_\Omega^t[\tilde{S}_{n,k;r},\cdot] \right)_{|V^{m-1,\alpha}(\partial\Omega)}
+\frac{\partial}{\partial\nu_\Omega}w_\Omega^+\left[\tilde{S}_{n,k;r},
C^{m,\alpha}(\partial\Omega)
\right]\,.
\end{eqnarray}
\end{theorem}
{\bf Proof.} Since $W_\Omega^t[\tilde{S}_{n,k;r},\cdot]$ maps  $V^{m-1,\alpha}(\partial\Omega)$ to itself (cf.~\textit{e.g.}, \cite[Cor.~10.1]{DoLa17} in case $m\geq 1$ and \cite[Cor.~8.5]{La25a} in case $m=0$), Proposition \ref{mldc_prop:-1/2+w,dnw+} (ii) implies that the right hand side of equality (\ref{mldc_thm:codamasumnwv1}) is contained in the left hand side of  equality (\ref{mldc_thm:codamasumnwv1}). On the other hand, if $\tau\in (S^{m,\alpha}_{k^2})^\perp
\cap V^{m-1,\alpha}(\partial\Omega)$, then Theorem \ref{mldc_thm:codasumnwv} implies that there exist $\phi\in  \left(C^{\max\{1,m\},\alpha}(\partial\Omega)\right)'$ and 
$\psi\in C^{m,\alpha}(\partial\Omega)$ such that
\begin{equation}\label{mldc_thm:codamasumnwv2}
\tau=-\frac{1}{2}\phi+W_\Omega^t[\tilde{S}_{n,k;r},\phi] 
+\frac{\partial}{\partial\nu_\Omega}w_\Omega^+\left[\tilde{S}_{n,k;r},\psi\right]\,.
\end{equation}
By Proposition \ref{mldc_prop:-1/2+w,dnw+} (ii), we have $\frac{\partial}{\partial\nu_\Omega}w_\Omega^+\left[\tilde{S}_{n,k;r},\psi\right]\in (S^{m,\alpha}_\lambda)^\perp\cap V^{m-1,\alpha}(\partial\Omega)$ and accordingly the above equality  (\ref{mldc_thm:codamasumnwv2}) implies that $-\frac{1}{2}\phi+W_\Omega^t[\tilde{S}_{n,k;r},\phi] $ belongs to
$(S^{m,\alpha}_\lambda)^\perp\cap V^{m-1,\alpha}(\partial\Omega)$. Hence, the regularity Theorem \ref{mldc_thm:pm1/2wtreg} implies that $\phi\in V^{m-1,\alpha}(\partial\Omega)$. Then equality (\ref{mldc_thm:codamasumnwv2}) implies that $\tau$ belongs to the right hand side of  equality (\ref{mldc_thm:codamasumnwv1}).\hfill  $\Box$ 

\vspace{\baselineskip}

Then by combining the necessary condition of Theorem \ref{mldc_thm:nexintneuhe}, Theorem \ref{mldc_thm:codamasumnwv} and the jump formulas for the normal derivative of the acoustic single and  double  layer potential, we deduce the validity of the following statement.
\begin{theorem}[of existence for the interior Neumann problem]\label{mldc_thm:exintneuhe}
  Let  \,  $m$ belong to ${\mathbb{N}}$.  Let  $\alpha\in]0,1[$. Let $\Omega$ be a bounded open subset of ${\mathbb{R}}^{n}$ of class $C^{\max\{1,m\},\alpha}$. 
  
  Let $k\in {\mathbb{C}}\setminus ]-\infty,0]$, ${\mathrm{Im}}\,k\geq 0$.   If  $g\in V^{m-1,\alpha}(\partial\Omega)$, then the interior Neumann problem
  \begin{equation}
\label{mldc_thm:exintneuhe1}
\left\{
\begin{array}{ll}
\Delta u+k^2 u=0 &{\text{in}}\ \Omega\,,
\\
\frac{\partial u}{\partial\nu_\Omega}=g
&{\text{on}}\ \partial\Omega\,
 \end{array}
\right.
\end{equation}
has   a  solution 
\begin{equation}\label{mldc_thm:exintneuhe1a}
u\in\left\{
\begin{array}{ll}
C^{0,\alpha}(\overline{\Omega})_\Delta &\text{if}\  m=0\,,
 \\
C^{m,\alpha}(\overline{\Omega})& \text{if}\ m\geq 1\,,
\end{array}
\right.  
\end{equation}
 if and only if 
\begin{equation}
\label{mldc_thm:exintneuhe2}
\langle g,v_{|\partial\Omega}\rangle =0
\end{equation}
for all $v\in C^{\max\{1,m\},\alpha}(\overline{\Omega})$ such that
\begin{equation}\label{mldc_thm:exintneuhe3}
\Delta v+k^2 v=0 \qquad\text{in}\ \Omega\,,\qquad \frac{\partial v}{\partial\nu_\Omega}=0\qquad\text{on}\ \partial\Omega\,.
\end{equation}
 \end{theorem}
{\bf Proof.} The necessity follows by Theorem \ref{mldc_thm:nexintneuhe}. If $g$ satisfies condition (\ref{mldc_thm:exintneuhe2}), then $g\in (S^{m,\alpha}_{k^2})^\perp\cap V^{m-1,\alpha}(\partial\Omega)$ and Theorem \ref{mldc_thm:codamasumnwv} implies the existence of
$\phi\in V^{m-1,\alpha}(\partial\Omega)$ and $\mu\in C^{m,\alpha}(\partial\Omega)$ such that
\[
g=-\frac{1}{2}\phi+W_\Omega^t[\tilde{S}_{n,k;r},\phi] +\frac{\partial}{\partial\nu_\Omega}
w_\Omega^+[\tilde{S}_{n,k;r},\mu]\,.
\]
Thus if we set
\[
u=v_\Omega^+[\tilde{S}_{n,k;r},\phi] +w_\Omega^+[\tilde{S}_{n,k;r},\mu]\,,
\]
then $u$ satisfies the memberships of (\ref{mldc_thm:exintneuhe1a}) by known  regularity results for the acoustic single layer potential (cf.~\textit{e.g.}, \cite[Thm.~7.1]{DoLa17} in case $m\geq 1$, \cite[Thm.~5.3]{La25a} in case $m=0$) and by known regularity results for the acoustic double layer   potential (cf.~\textit{e.g.}, 
 \cite[Thm.~7.3]{DoLa17} in case $m\geq 1$, 
 \cite[Thm.~5.1]{La25a} in case $m=0$). Also, the known  jump formulas for the normal derivatve of the acoustic simple layer potential of \cite[Thm.~6.3]{La25a} imply that $\frac{\partial u}{\partial\nu_\Omega}=g$. Hence, $u$ satisfies the Neumann problem (\ref{mldc_thm:exintneuhe1}).\hfill  $\Box$ 

\vspace{\baselineskip}
 
\section{A nonvariational form of the exterior Neumann problem for the Helmholtz equation}\label{mldc_sec:enpddh}
We now consider the exterior   Neumann problem in the case in which the Neumann datum is in the space $V^{m-1,\alpha}(\partial\Omega)$ (see (\ref{mldc_eq:vm-1a}). Nonvariational in case $m=0$). For the corresponding classical result in case $\Omega^-$ is connected and the Neumann data satisfy some extra regularity  assumption,  we refer to  Colton and Kress \cite[Thm.~3.25]{CoKr92}.  For case $m=0$ and in the specific case in which solutions can be represented in terms of an acoustic single  layer potential, with some extra assumptions on   $\Omega$, we refer to \cite[Thm.~12.1]{La25b}.  Next, we  introduce the following classical necessary condition for the existence (see also \cite[Thm.~12.1]{La25b}, which we follow).
\begin{theorem}\label{mldc_thm:nexextneuhe}
 Let $m\in {\mathbb{N}}$,   $\alpha\in]0,1[$.  Let $\Omega$ be a bounded open subset of ${\mathbb{R}}^{n}$ of class $C^{\max\{1,m\},\alpha}$. 
  
  Let $k\in {\mathbb{C}}\setminus ]-\infty,0]$, ${\mathrm{Im}}\,k\geq 0$.   Let  $g\in V^{m-1,\alpha}(\partial\Omega)$. If the exterior Neumann problem
  \begin{equation}
\label{mldc_thm:nexextneuhe1}
\left\{
\begin{array}{ll}
\Delta u+k^2 u=0 &{\text{in}}\ \Omega^-\,,
\\
-\frac{\partial u}{\partial\nu_{\Omega^-}}=g
&{\text{on}}\ \partial\Omega\,,
\\
u\ \text{satisfies\ the\ outgoing}\  k-\text{radiation\ condition}\,,
 \end{array}
\right.
\end{equation}
has   a  solution $u\in C^{0,\alpha}_{{\mathrm{loc}}}(\overline{\Omega^-})_\Delta$, then
\begin{equation}
\label{mldc_thm:nexextneuhe2}
\langle g,v_{|\partial\Omega}\rangle =0
\end{equation}
for all $v\in C^{\max\{1,m\},\alpha}_{{\mathrm{loc}}}(\overline{\Omega^-})$ such that
\begin{eqnarray}\label{mldc_thm:nexextneuhe3}
&&\Delta v+k^2 v=0 \qquad\text{in}\ \Omega^-\,,\qquad \frac{\partial v}{\partial\nu_{\Omega^-}}=0\qquad\text{on}\ \partial\Omega\,.
\\\nonumber
&& v \ \text{satisfies\ the\ outgoing\ $k$-radiation\ condition.	}
\end{eqnarray}
\end{theorem}
{\bf Proof.} Assume that  $u\in C^{0,\alpha}_{{\mathrm{loc}}}(\overline{\Omega^-})_\Delta$  solves the exterior Neumann problem 
(\ref{mldc_thm:nexextneuhe2}) and that $v\in C^{\max\{1,m\},\alpha}_{{\mathrm{loc}}}(\overline{\Omega^-})$ satisfies  (\ref{mldc_thm:nexextneuhe3}).  Then Remark \ref{mldc_rem:helpo} implies that $v\in C^{1,\alpha}_{{\mathrm{loc}}}(\overline{\Omega^-})_\Delta$ and 	  the second Green Identity for exterior domains of \cite[Thm.~8.18]{La25b}  implies that
\[
\langle g,v_{|\partial\Omega}\rangle=\langle-\frac{\partial u}{\partial\nu_{\Omega^-}},v_{|\partial\Omega}\rangle
=-\int_{\partial\Omega}u\frac{\partial v}{\partial\nu_{\Omega^-}}\,d\sigma=0\,.
\]
\hfill  $\Box$ 

\vspace{\baselineskip}

By Theorem \ref{mldc_thm:nexextneuhe}, a necessary condition for the solvability of the Neumann problem is that condition (\ref{mldc_thm:nexextneuhe2}) holds true. 
Thus we now introduce the following notation. If $m\in{\mathbb{N}}$, $\alpha\in]0,1[$,  $\Omega$ is a bounded open subset of ${\mathbb{R}}^{n}$ of class $C^{\max\{1,m\},\alpha}$,  $k\in {\mathbb{C}}\setminus ]-\infty,0]$, ${\mathrm{Im}}\,k\geq 0$,  then we set
\begin{eqnarray}\label{mldc_eq:smal-}
\lefteqn{S^{m,\alpha}_{k^2,-}\equiv\biggl\{
 v_{|\partial\Omega}:\, v\in C^{\max\{1,m\},\alpha}_{{\mathrm{loc}}	}(\overline{\Omega^-})\,,\ 
 }
 \\ \nonumber
 &&\qquad\qquad\quad
 \Delta v+k^2 v=0 \ \text{in}\ \Omega^-\,,\  \frac{\partial v}{\partial\nu_{\Omega^-}}=0\ \text{on}\ \partial\Omega\,,
 \\ \nonumber
 &&\qquad\qquad\quad
 v \ \text{satisfies\ the\ outgoing\ $k$-radiation\ condition}
\biggr\}
\end{eqnarray}
and
\begin{equation}\label{mldc_eq:smal-1}
(S^{m,\alpha}_{k^2,-})^\perp\equiv
\left\{
g\in \left(C^{\max\{1,m\},\alpha}(\partial\Omega)\right)':\, \langle g, f \rangle=0\ \forall f \in S^{m,\alpha}_{k^2,-}
\right\}\,.
\end{equation}
By the uniqueness Theorem of the solutions of  Helmholtz equation with Neumann boundary conditions  on the unbounded connected component $(\Omega^-)_0$ of $\Omega^-$ (cf.~\cite[Thm.~1.3 (i)]{La25}), we have
\begin{eqnarray}\label{mldc_eq:smal-2}
\lefteqn{S^{m,\alpha}_{k^2,-}=\biggl\{
 v_{|\partial\Omega}:\, v\in C^{\max\{1,m\},\alpha}_{{\mathrm{loc}}	}(\overline{\Omega^-})\,,\ 
 }
 \\ \nonumber
 &&\qquad\qquad\quad
 \Delta v+k^2 v=0 \ \text{in}\ \Omega^-\,,\  \frac{\partial v}{\partial\nu_{\Omega^-}}=0\ \text{on}\ \partial\Omega\,,
 \\ \nonumber
 &&\qquad\qquad\quad
 \  v(x)=0\ \forall x\in (\Omega^-)_0
\biggr\}\,.
\end{eqnarray}
\begin{proposition}\label{mldc_prop:-1/2+w,dnw-}
 Let $m\in{\mathbb{N}}$, $\alpha\in]0,1[$. Let $\Omega$ be a bounded open subset of ${\mathbb{R}}^{n}$ of class $C^{\max\{1,m\},\alpha}$. Let $k\in {\mathbb{C}}\setminus ]-\infty,0]$, ${\mathrm{Im}}\,k\geq 0$.  Then the following statements hold true. 
 \begin{enumerate}
\item[(i)] $  \frac{1}{2}\phi+W_\Omega^t[\tilde{S}_{n,k;r},\phi]\in (S^{m,\alpha}_{k^2,-})^\perp$ for all $\phi\in  \left(C^{\max\{1,m\},\alpha}(\partial\Omega)\right)'$.
 \item[(ii)]  $\frac{\partial}{\partial\nu_{\Omega^-}}w_\Omega^-[\tilde{S}_{n,k;r},\psi]\in 
   (S^{m,\alpha}_{k^2,-})^\perp\cap V^{m-1,\alpha}(\partial\Omega)$ for all $\psi\in C^{m,\alpha}(\partial\Omega)$.
\end{enumerate}
\end{proposition}
{\bf Proof.} (i) If $\phi\in  \left(C^{\max\{1,m\},\alpha}(\partial\Omega)\right)'$ and $f\in S^{m,\alpha}_{k^2,-}$, then
\[
\langle \frac{1}{2}\phi+W_\Omega^t[\tilde{S}_{n,k;r},\phi],f \rangle=
\langle\phi,\frac{1}{2}f +W_\Omega[\tilde{S}_{n,k;r},f]\rangle\,.
\]
Since Corollary 11.8 (i) of \cite{La25b} ensures that $\frac{1}{2}f +W_\Omega[\tilde{S}_{n,k;r},f]=0$, we conclude that the right hand side equals $0$. 

(ii) By known regularity results for the acoustic double layer   potential (cf.~\textit{e.g.}, 
 \cite[Thm.~7.3]{DoLa17} in case $m\geq 1$, 
 \cite[Thm.~5.1]{La25a} in case $m=0$), we have $w_\Omega^-[\tilde{S}_{n,k;r},\psi]\in 
 C^{m,\alpha}_{{\mathrm{loc}}}(\overline{\Omega^-})$
  and
  accordingly $\frac{\partial}{\partial\nu_{\Omega^-}}w_\Omega^-[\tilde{S}_{n,k;r},\psi]$ belongs to $ V^{m-1,\alpha}(\partial\Omega)$ (in case $m=0$, see also Remark \ref{mldc_rem:helpo}
 and Proposition \ref{mldc_prop:recodnu}). Moreover, $w_\Omega^-[\tilde{S}_{n,k;r},\psi]$ satisfies the outgoing $k$-radiation condition (cf.~\textit{e.g.}, \cite[Thm.~6.20 (ii)]{La25b}). Then 
statement (ii) follows  by  Remark \ref{mldc_rem:helpo} and Theorem \ref{mldc_thm:nexextneuhe} with $g=\frac{\partial}{\partial\nu_{\Omega^-}}w_\Omega^-[\tilde{S}_{n,k;r},\psi]$.\hfill  $\Box$ 

\vspace{\baselineskip}

\begin{proposition}\label{mldc_prop:isomquotn-}
Let $m\in{\mathbb{N}}$, $\alpha\in]0,1[$. Let $\Omega$ be a bounded open subset of ${\mathbb{R}}^{n}$ of class $C^{\max\{1,m\},\alpha}$.  
 Let $k\in {\mathbb{C}}\setminus ]-\infty,0]$, ${\mathrm{Im}}\,k\geq 0$. Then 
 \begin{equation}\label{mldc_prop:isomquotn-1}
{\mathrm{Im}}\, \left(
\frac{1}{2}I+W_\Omega^t[\tilde{S}_{n,k;r},\cdot] \right)
=\left({\mathrm{Ker}}\,\left(\frac{1}{2}I+W_\Omega[\tilde{S}_{n,k;r},\cdot]\right)_{|C^{\max\{1,m\},\alpha}(\partial\Omega)}\right)^\perp
\end{equation}
 and ${\mathrm{Im}}\, \left(
\frac{1}{2}I+W_\Omega^t[\tilde{S}_{n,k;r},\cdot] \right)$ is a weakly$^\ast$ closed subspace of
$(S^{m,\alpha}_{k^2,-})^\perp$ of finite codimension. Here,
$W_\Omega^t[\tilde{S}_{n,k;r},\cdot]$ is the transpose to the operator $W_\Omega[\tilde{S}_{n,k;r},\cdot]$ in $C^{\max\{1,m\},\alpha}(\partial\Omega)$ and the annihilator is taken with respect to the natural duality pairing (\ref{mldc_rem:wtnotation1}).
\end{proposition}
{\bf Proof.} Since $W_\Omega[\tilde{S}_{n,k;r},\cdot]$ is compact in  $C^{\max\{1,m\},\alpha}(\partial\Omega)$  (cf.~\textit{e.g.},   \cite[Cor.~9.1]{DoLa17}), then the Fredholm Alternative Theorem  in the natural duality pairing (\ref{mldc_rem:wtnotation1}) implies the validity of   equality (\ref{mldc_prop:isomquotn-1}) and that 
${\mathrm{Im}}\, \left(
-\frac{1}{2}I+W_\Omega^t[S_{n,\lambda},\cdot] \right)$ has finite codimension in $\left(C^{\max\{1,m\},\alpha}(\overline{\Omega})\right)'$.  Since 
\[
\left({\mathrm{Ker}}\,\left(\frac{1}{2}I+W_\Omega[
\tilde{S}_{n,k;r},\cdot]\right)_{|C^{\max\{1,m\},\alpha}(\partial\Omega)}\right)^\perp
\]
 is weakly$^\ast$ closed in $\left(C^{\max\{1,m\},\alpha}(\partial\Omega)\right)'$, then the space  ${\mathrm{Im}}\, \left(
\frac{1}{2}I+W_\Omega^t[\tilde{S}_{n,k;r},\cdot] \right)$ is also weakly$^\ast$ closed in    $\left(C^{\max\{1,m\},\alpha}(\partial\Omega)\right)'$. 

By Proposition \ref{mldc_prop:-1/2+w,dnw-} (i), we know that the space ${\mathrm{Im}}\, \left(
\frac{1}{2}I+W_\Omega^t[\tilde{S}_{n,k;r},\cdot] \right)$ is contained in $(S^{m,\alpha}_{k^2,-})^\perp$, which is also weakly$^\ast$ closed in $\left(C^{\max\{1,m\},\alpha}(\partial\Omega)\right)'$. Hence, the proof is complete.\hfill  $\Box$ 

\vspace{\baselineskip}

Then we can prove the following statement that says that the space of compatible data $(S^{m,\alpha}_{k^2,-})^\perp$ can be written as an algebraic sum of normal derivatives of an acoustic single layer potential and of normal derivatives of an acoustic double layer potential. 
\begin{theorem}\label{mldc_thm:codasumnwv-}
Let $m\in{\mathbb{N}}$, $\alpha\in]0,1[$. Let $\Omega$ be a bounded open subset of ${\mathbb{R}}^{n}$ of class $C^{\max\{1,m\},\alpha}$.  Let $k\in {\mathbb{C}}\setminus ]-\infty,0]$, ${\mathrm{Im}}\,k\geq 0$. Then
\begin{equation}\label{mldc_thm:codasumnwv-1}
(S^{m,\alpha}_{k^2,-})^\perp=
{\mathrm{Im}}\, \left(
\frac{1}{2}I+W_\Omega^t[\tilde{S}_{n,k;r},\cdot] \right)
+\frac{\partial}{\partial\nu_{\Omega^-}}w_\Omega^-\left[\tilde{S}_{n,k;r},
C^{m,\alpha}(\partial\Omega)
\right]\,,
\end{equation}
where  $W_\Omega^t[\tilde{S}_{n,k;r},\cdot]$ is the transpose to the operator $W_\Omega[\tilde{S}_{n,k;r},\cdot]$ in $C^{\max\{1,m\},\alpha}(\partial\Omega)$.\end{theorem}
{\bf Proof.} By Proposition \ref{mldc_prop:isomquotn-}, ${\mathrm{Im}}\, \left(
\frac{1}{2}I+W_\Omega^t[\tilde{S}_{n,k;r},\cdot] \right)$ is a weakly$^\ast$ closed subspace of
$(S^{m,\alpha}_{k^2,-})^\perp$. Let $\pi$ be the canonical projection  of  the topological vector space 
$(S^{m,\alpha}_{k^2,-})^\perp$ with the weak$^\ast$ topology onto the finite dimensional quotient space\[
Q^{m,\alpha}_{k^2,-}\equiv 
\frac{	
(S^{m,\alpha}_{k^2,-})^\perp
}{
{\mathrm{Im}}\, \left(
\frac{1}{2}I+W_\Omega^t[\tilde{S}_{n,k;r},\cdot]
\right)
} \,.
\]
By Proposition \ref{mldc_prop:-1/2+w,dnw-} (ii),
\[
\frac{\partial}{\partial\nu_{\Omega^-}}w_\Omega^-\left[\tilde{S}_{n,k;r},
C^{m,\alpha}(\partial\Omega)
\right]
\leq (S^{m,\alpha}_{k^2,-})^\perp
\]
 and thus  the $\pi$-image $\pi\left[T^{m,\alpha}_{k^2,-}\right]$
of
\[
T^{m,\alpha}_{k^2,-}
\equiv
{\mathrm{Im}}\, \left(
\frac{1}{2}I+W_\Omega^t[\tilde{S}_{n,k;r},\cdot]
\right)
+\frac{\partial}{\partial\nu_{\Omega^-}}w_\Omega^-\left[\tilde{S}_{n,k;r},
C^{m,\alpha}(\partial\Omega)
\right]
\]
is a   subspace of the finite dimensional Hausdorff topological     quotient vector space  $ 
Q^{m,\alpha}_{k^2,-}$ and is accordingly closed in $Q^{m,\alpha}_{k^2,-}$. Since the canonical projection $\pi$ is continuous, then the preimage
\[
\pi^{\leftarrow}\left[\pi\left[T^{m,\alpha}_{k^2,-} \right]\right]
\]
is weakly$^\ast$ closed in $(S^{m,\alpha}_{k^2,-} )^\perp$.  Now a simple computation based on the definition of $\pi$  shows that
\[
\pi^{\leftarrow}\left[\pi\left[T^{m,\alpha}_{k^2,-} \right]\right]=T^{m,\alpha}_{k^2,-}
\]
and accordingly $T^{m,\alpha}_{k^2,-}$ is a closed subspace of $(S^{m,\alpha}_{k^2,-})^\perp$. Assume by contradiction that $T^{m,\alpha}_{k^2,-}\subsetneq (S^{m,\alpha}_{k^2,-})^\perp$. Then   there exists
\[
\tau_-\in (S^{m,\alpha}_{k^2,-})^\perp\setminus T^{m,\alpha}_{k^2,-}\,.
\]
By the Hahn-Banach Theorem in the dual space
\[
\left(C^{\max\{1,m\},\alpha}(\partial\Omega)\right)'_s
\]
with the weak$^\ast$ topology, there exists
\[
\tilde{\mu}_-\in\left(\left(C^{\max\{1,m\},\alpha}(\partial\Omega)\right)'_s\right)'
\]
such that
\begin{equation}\label{mldc_thm:codasumnwv-2}
\langle \tilde{\mu}_-,\tau_-\rangle=1\,,\qquad \langle \tilde{\mu},x'\rangle=0\qquad
\forall x' \in T^{m,\alpha}_{k^2,-}\,.
\end{equation}
Now the canonical isometry of $C^{\max\{1,m\},\alpha}(\partial\Omega)$ is an isomorphism of the space $C^{\max\{1,m\},\alpha}(\partial\Omega)$ onto the bi-dual
\[
\left(\left(C^{\max\{1,m\},\alpha}(\partial\Omega)\right)'_s\right)'\,,
\]
which is a subspace of the strong bidual $\left(C^{\max\{1,m\},\alpha}(\partial\Omega)\right)''$
(cf.~\textit{e.g.}, Brezis \cite[Prop.~3.14, p.~64]{Br11}). Hence, there exists $\mu_-\in C^{\max\{1,m\},\alpha}(\partial\Omega)$ such that
\[
\langle \tilde{\mu}_-,x'\rangle=\langle x',\mu_-\rangle\qquad\forall x'\in 
\left(C^{\max\{1,m\},\alpha}(\partial\Omega)\right)'
\]
and we must have
\begin{equation}\label{mldc_thm:codasumnwv-3}
\langle\tau_-,\mu_- \rangle=1\,,\qquad \langle x',\mu_-\rangle=0\qquad
\forall x' \in T^{m,\alpha}_{k^2,-}\,.
\end{equation}
In particular, $\mu_-$ annihilates the subspace ${\mathrm{Im}}\, \left(
\frac{1}{2}I+W_\Omega^t[\tilde{S}_{n,k;r},\cdot]
\right)$ of $T^{m,\alpha}_{k^2,-}$ and accordingly Proposition \ref{mldc_prop:isomquotn-} implies that
\begin{eqnarray}\label{mldc_thm:codasumnwv-4}
\lefteqn{\mu_-\in
\left(
{\mathrm{Ker}}\,\left(\frac{1}{2}I+W_\Omega[\tilde{S}_{n,k;r},\cdot]\right)_{|C^{\max\{1,m\},\alpha}(\partial\Omega)}
\right)^{\perp\perp}
}
\\ \nonumber
&&\qquad\qquad\qquad\qquad
=
{\mathrm{Ker}}\,\left(\frac{1}{2}I+W_\Omega[\tilde{S}_{n,k;r},\cdot]\right)_{|C^{\max\{1,m\},\alpha}(\partial\Omega)}\,,
\end{eqnarray}
where the annihilators are taken with respect to the natural duality 
 (\ref{mldc_rem:wtnotation1}) of $C^{\max\{m,1\},\alpha}(\partial\Omega)$ 
(cf.~\textit{e.g.}, Rudin, \cite[Thm.~4.7, Chap. IV]{Ru91}). 
Hence,  the jump formulas for the acoustic double layer   potential imply that 
\begin{equation}\label{mldc_thm:codasumnwv-4a}
w_\Omega^+\left[\tilde{S}_{n,k;r},\mu_-\right]=0\qquad\text{on}\ \partial\Omega
\end{equation}
and accordingly that 
$w_\Omega^+\left[\tilde{S}_{n,k;r},\mu_-\right]$ solves the Helmholtz equation in $\Omega$ that corresponds to $k^2$ with zero Dirichlet boundary condition (cf.~\textit{e.g.},  \cite[Thm.~7.1 (ii)]{DoLa17}). 
By known regularity results for the  acoustic  double layer   potential (cf.~\textit{e.g.}, 
 \cite[Thm.~7.3]{DoLa17}), we have $w_\Omega^+[\tilde{S}_{n,k;r},\mu_-]\in 
 C^{\max\{1,m\},\alpha}(\overline{\Omega})$. Since $w_\Omega^+\left[\tilde{S}_{n,k;r},\mu_-\right]$ solves the Helmholtz equation in $\Omega$ that corresponds to $k^2$ with zero Dirichlet boundary condition 
Corollary \ref{mldc_corol:eipdd} (i) implies that
\begin{equation}\label{mldc_thm:codasumnwv-4b}
\frac{\partial}{\partial\nu_{\Omega}}w_\Omega^+\left[\tilde{S}_{n,k;r},
\mu_-
\right]\in {\mathrm{Ker}}\,\left(\frac{1}{2}I+W_\Omega^t[\tilde{S}_{n,k;r},\cdot]
\right)_{|C^{\max\{1,m\}-1,\alpha}(\partial\Omega)}
\end{equation}
Then the jump formula for the   normal derivative of the acoustic double layer potential implies that
\begin{eqnarray}\label{mldc_thm:codasumnwv-5}
\lefteqn{-\frac{\partial}{\partial\nu_{\Omega^-}}w_\Omega^-\left[\tilde{S}_{n,k;r},
\mu_-
\right]
}
\\ \nonumber
&&\qquad\qquad 
=
\frac{\partial}{\partial\nu_{\Omega}}w_\Omega^+\left[\tilde{S}_{n,k;r},
\mu_-
\right]\in {\mathrm{Ker}}\,\left(\frac{1}{2}I+W_\Omega^t[\tilde{S}_{n,k;r},\cdot]\right)_{|C^{\max\{1,m\}-1,\alpha}(\partial\Omega)}
\end{eqnarray}
(cf.~\textit{e.g.}, \cite[Thm.~6.5]{La25a}). Moreover, the jump formula for the  acoustic double layer potential implies that
\begin{equation}\label{mldc_thm:codasumnwv-6}
\mu_-=w_\Omega^+\left[\tilde{S}_{n,k;r},\mu_-\right]-w_\Omega^-\left[\tilde{S}_{n,k;r},\mu_-\right]
=-w_\Omega^-\left[\tilde{S}_{n,k;r},\mu_-\right]\qquad\text{on}\ \partial\Omega\,,
\end{equation}
(cf.~(\ref{mldc_thm:dlay1a}), (\ref{mldc_thm:codasumnwv-4a})). 
Since $\mu_-$ annihilates the subspace $\frac{\partial}{\partial\nu_{\Omega^-}}w_\Omega^-\left[\tilde{S}_{n,k;r},
C^{m,\alpha}(\partial\Omega)
\right]$ of $T^{m,\alpha}_{k^2,-}$, we have
\[
\langle
\frac{\partial}{\partial\nu_{\Omega^-}}w_\Omega^-\left[\tilde{S}_{n,k;r},\eta\right]
,\mu_- \rangle=0
\qquad\forall \eta\in C^{m,\alpha}(\partial\Omega)\,.
\]
Hence, Remark \ref{mldc_rem:helpo} (iii) and
 the second Green Identity in distributional form for solutions of the Helmholtz equation of \cite[Thm.~8.18]{La25b} implies that
\begin{eqnarray*}
\lefteqn{
0=\langle
\frac{\partial}{\partial\nu_{\Omega^-}}w_\Omega^-\left[\tilde{S}_{n,k;r},\eta\right]
,\mu_- \rangle
}
\\ \nonumber
&&\qquad
=-\langle
\frac{\partial}{\partial\nu_{\Omega^-}}w_\Omega^-\left[\tilde{S}_{n,k;r},\eta\right]
,w_\Omega^-\left[\tilde{S}_{n,k;r},\mu_-\right]_{|\partial\Omega} \rangle
\\ \nonumber
&&\qquad
=-\int_{\partial\Omega}w_\Omega^-\left[\tilde{S}_{n,k;r},\eta\right]
\frac{\partial}{\partial\nu_{\Omega^-}}w_\Omega^-\left[\tilde{S}_{n,k;r},\mu_-\right]\,d\sigma
\qquad\forall \eta\in C^{m,\alpha}(\partial\Omega)\,.
\end{eqnarray*}
In particular, the same equality holds for all $\eta\in  C^{\max\{m,1\},\alpha}(\partial\Omega)$. 
Hence, the jump formula (\ref{mldc_thm:dlay1a}) for the acoustic double layer  potential implies that
\begin{eqnarray}\label{mldc_thm:codasumnwv-7}
\lefteqn{\frac{\partial}{\partial\nu_{\Omega^-}}w_\Omega^-\left[\tilde{S}_{n,k;r},\mu_-\right]
\in 
\left({\mathrm{Im}}\,
 \left(-\frac{1}{2}I+W_\Omega[\tilde{S}_{n,k;r},\cdot]\right)_{|C^{\max\{1,m\},\alpha}(\partial\Omega)}\right)^\perp
}
\\ \nonumber
&&\qquad\qquad\qquad\qquad\qquad\qquad\qquad\qquad\qquad\qquad
={\mathrm{Ker}}\,\left(-\frac{1}{2}I+W_\Omega^t[\tilde{S}_{n,k;r},\cdot]\right) \,,
 \end{eqnarray}
 where the annihilator is taken with respect to the natural duality (\ref{mldc_rem:wtnotation1}) of the space $C^{\max\{m,1\},\alpha}(\partial\Omega)$ (cf.~\textit{e.g.}, Rudin, \cite[Thm.~4.12, Chap. IV]{Ru91}). 
Then by combining  (\ref{mldc_thm:codasumnwv-5}) and the membership of (\ref{mldc_thm:codasumnwv-7}), we deduce that
\[
\frac{\partial}{\partial\nu_{\Omega^-}}w_\Omega^-\left[\tilde{S}_{n,k;r},\mu_-\right]
=0\,.
\]
 Since $w_\Omega^-\left[\tilde{S}_{n,k;r},\mu_-\right]$ solves the Helmholtz equation in $\Omega^-$ with zero
Neumann boundary conditions and  the outgoing $k$-radiation condition  (cf.~\textit{e.g.}, \cite[Thm.~6.20 (ii)]{La25b}),  then we have
\[
\mu_-=-w_\Omega^-\left[\tilde{S}_{n,k;r},\mu_-\right]_{|\partial\Omega}\in S^{m,\alpha}_{k^2,-}\,.
\]
Since $\tau_-\in (S^{m,\alpha}_{k^2,-})^\perp$, we have $\langle\tau_-,\mu_-\rangle=0$, a contradiction (see (\ref{mldc_thm:codasumnwv-3})). Hence, $T^{m,\alpha}_{k^2,-} = (S^{m,\alpha}_{k^2,-})^\perp$ and thus the proof is complete.\hfill  $\Box$ 

\vspace{\baselineskip}

Then by combining Proposition \ref{mldc_prop:-1/2+w,dnw-} (ii), Theorem \ref{mldc_thm:codasumnwv-} and the regularity Theorem \ref{mldc_thm:pm1/2wtreg}, we deduce the validity of the following statement.
\begin{theorem}\label{mldc_thm:codamasumnwv-}
Let $m\in{\mathbb{N}}$, $\alpha\in]0,1[$. Let $\Omega$ be a bounded open subset of ${\mathbb{R}}^{n}$ of class $C^{\max\{1,m\},\alpha}$.  Let $k\in {\mathbb{C}}\setminus ]-\infty,0]$, ${\mathrm{Im}}\,k\geq 0$. Then
\begin{eqnarray}\label{mldc_thm:codamasumnwv-1}
\lefteqn{
(S^{m,\alpha}_{k^2,-})^\perp\cap V^{m-1,\alpha}(\partial\Omega)
}
\\ \nonumber
&&\qquad=
{\mathrm{Im}}\, \left(
\frac{1}{2}I+W_\Omega^t[\tilde{S}_{n,k;r},\cdot] \right)_{|V^{m-1,\alpha}(\partial\Omega)}
+\frac{\partial}{\partial\nu_{\Omega^-}}w_\Omega^-\left[\tilde{S}_{n,k;r},
C^{m,\alpha}(\partial\Omega)
\right]\,.
\end{eqnarray}
\end{theorem}
{\bf Proof.} Since $W_\Omega^t[\tilde{S}_{n,k;r},\cdot]$ maps  $V^{m-1,\alpha}(\partial\Omega)$ to itself (cf.~\textit{e.g.}, \cite[Cor.~10.1]{DoLa17} in case $m\geq 1$ and \cite[Cor.~8.5]{La25a} in case $m=0$), Proposition \ref{mldc_prop:-1/2+w,dnw-} (ii) implies that the right hand side of equality (\ref{mldc_thm:codamasumnwv-1}) is contained in the left hand side of  equality (\ref{mldc_thm:codamasumnwv-1}). On the other hand, if $\tau\in (S^{m,\alpha}_{k^2,-})^\perp
\cap V^{m-1,\alpha}(\partial\Omega)$, then Theorem \ref{mldc_thm:codasumnwv-} implies that there exist $\phi\in  \left(C^{\max\{1,m\},\alpha}(\partial\Omega)\right)'$ and 
$\psi\in C^{m,\alpha}(\partial\Omega)$ such that
\begin{equation}\label{mldc_thm:codamasumnwv-2}
\tau=\frac{1}{2}\phi+W_\Omega^t[\tilde{S}_{n,k;r},\phi] 
+\frac{\partial}{\partial\nu_{\Omega^-}}w_\Omega^-\left[\tilde{S}_{n,k;r},\psi\right]\,.
\end{equation}
By Proposition \ref{mldc_prop:-1/2+w,dnw-} (ii), we have $\frac{\partial}{\partial\nu_{\Omega^-}}w_\Omega^-\left[\tilde{S}_{n,k;r},\psi\right]\in (S^{m,\alpha}_{k^2,-})^\perp\cap V^{m-1,\alpha}(\partial\Omega)$ and accordingly the above equality  (\ref{mldc_thm:codamasumnwv-2}) implies that $\frac{1}{2}\phi+W_\Omega^t[\tilde{S}_{n,k;r},\phi] $ belongs to
$(S^{m,\alpha}_{k^2,-})^\perp\cap V^{m-1,\alpha}(\partial\Omega)$. Hence, the regularity Theorem \ref{mldc_thm:pm1/2wtreg} implies that $\phi\in V^{m-1,\alpha}(\partial\Omega)$. Then equality (\ref{mldc_thm:codamasumnwv-2}) implies that $\tau$ belongs to the right hand side of  equality (\ref{mldc_thm:codamasumnwv-1}).\hfill  $\Box$ 

\vspace{\baselineskip}

Then by combining the necessary condition of Theorem \ref{mldc_thm:nexextneuhe}, Theorem \ref{mldc_thm:codamasumnwv-} and the jump formulas for the normal derivative of the acoustic  single and double  layer potential, we deduce the validity of the following statement.
\begin{theorem}[of existence for the exterior Neumann problem]\label{mldc_thm:exextneuhe}
  Let   $m$ belong to ${\mathbb{N}}$. Let    $\alpha\in]0,1[$. Let $\Omega$ be a bounded open subset of ${\mathbb{R}}^{n}$ of class $C^{\max\{1,m\},\alpha}$. 
  
  Let $k\in {\mathbb{C}}\setminus ]-\infty,0]$, ${\mathrm{Im}}\,k\geq 0$.  If  $g\in V^{m-1,\alpha}(\partial\Omega)$, then the exterior Neumann problem
  \begin{equation}
\label{mldc_thm:exextneuhe1}
\left\{
\begin{array}{ll}
\Delta u+k^2 u=0 &{\text{in}}\ \Omega^-\,,
\\
-\frac{\partial u}{\partial\nu_{\Omega^-}}=g
&{\text{on}}\ \partial\Omega\,,
\\
u\ \text{satisfies\ the\ outgoing}\  k-\text{radiation\ condition}\,,
 \end{array}
\right.
\end{equation}
has   a  solution 
\begin{equation}\label{mldc_thm:exextneuhe1a}
u\in\left\{
\begin{array}{ll}
C^{0,\alpha}_{{\mathrm{loc}}}(\overline{\Omega^-})_\Delta &\text{if}\  m=0\,,
 \\
C^{m,\alpha}_{{\mathrm{loc}}}(\overline{\Omega^-})& \text{if}\ m\geq 1\,,
\end{array}
\right.  
\end{equation}
 if and only if 
\begin{equation}
\label{mldc_thm:exextneuhe2}
\langle g,v_{|\partial\Omega}\rangle =0
\end{equation}
for all $v\in C^{\max\{1,m\},\alpha}_{{\mathrm{loc}}}(\overline{\Omega^-})$ such that
\begin{eqnarray}\label{mldc_thm:exextneuhe3}
&&\Delta v+k^2 v=0 \qquad\text{in}\ \Omega^-\,,\qquad \frac{\partial v}{\partial\nu_{\Omega^-}}=0\qquad\text{on}\ \partial\Omega\,,
\\\nonumber
&& v \ \text{satisfies\ the\ outgoing\ $k$-radiation\ condition.	}
\end{eqnarray}
\end{theorem}
{\bf Proof.} The necessity follows by Theorem \ref{mldc_thm:nexextneuhe}. If $g$ satisfies condition (\ref{mldc_thm:exextneuhe2}), then $g\in (S^{m,\alpha}_{k^2,-})^\perp$ and Theorem \ref{mldc_thm:codamasumnwv-} implies the existence of
$\phi_-\in V^{m-1,\alpha}(\partial\Omega)$ and $\mu_-\in C^{m,\alpha}(\partial\Omega)$ such that
\[
-g=-\frac{1}{2}\phi_--W_\Omega^t[\tilde{S}_{n,k;r},\phi_-] +\frac{\partial}{\partial\nu_{\Omega^-}}
w_\Omega^-[\tilde{S}_{n,k;r},\mu_-]\,.
\]
Thus if we set
\[
u=v_\Omega^-[\tilde{S}_{n,k;r},\phi_-] +w_\Omega^-[\tilde{S}_{n,k;r},\mu_-]\,,
\]
then $u$ satisfies the memberships of (\ref{mldc_thm:exextneuhe1a}) by known  regularity results for the acoustic single layer potential (cf.~\textit{e.g.}, \cite[Thm.~7.1]{DoLa17} in case $m\geq 1$, \cite[Thm.~5.3]{La25a} in case $m=0$) and by known regularity results for the acoustic double   layer  potential (cf.~\textit{e.g.}, 
 \cite[Thm.~7.3]{DoLa17} in case $m\geq 1$, 
 \cite[Thm.~5.1]{La25a} in case $m=0$). Then   $v_\Omega^-[\tilde{S}_{n,k;r},\phi_-]$ and $w_\Omega^-[\tilde{S}_{n,k;r},\mu_-]$ are known to satisfy the outgoing $k$-radiation condition for $m\geq 1$ and the same holds in case $m=0$ by \cite[Thm.~6.20]{La25b}.  Also, the known  jump formulas for the acoustic single layer potential of \cite[Thm.~6.3]{La25a} imply that $-\frac{\partial u}{\partial\nu_{\Omega^-}}=g$. Hence, $u$ satisfies the exterior Neumann problem (\ref{mldc_thm:exextneuhe1}).\hfill  $\Box$ 

\vspace{\baselineskip}

\section{A representation formula for the solutions of the Helmholtz equation}\label{mldc_sec:refosoh}

\begin{theorem}
\label{mldc_thm:reprhelfun}
Let  $m\in{\mathbb{N}}$, $\alpha\in]0,1[$. Let $\Omega$ be a bounded open subset of ${\mathbb{R}}^{n}$ of class $C^{\max\{1,m\},\alpha}$.  Let $k\in {\mathbb{C}}\setminus ]-\infty,0]$, ${\mathrm{Im}}\,k\geq 0$. 
Assume that $k^{2}$ is not a Dirichlet eigenvalue for $-\Delta$ in $(\Omega^-)_{j}$ for each $j\in {\mathbb{N}}\setminus\{0\}$ such that $j\leq \varkappa^{-}$. Let $A^{m-1,\alpha}$ be a topological supplement of the finite dimensional subspace
\[
{\mathrm{Ker}}\left(-\frac{1}{2}I+W_\Omega^t[\tilde{S}_{n,k;r},\cdot]
\right)_{|C^{\max\{1,m\}-1,\alpha}(\partial\Omega)} 
\]
of $V^{m-1,\alpha}(\partial\Omega)$, \textit{i.e.}, 
\[
V^{m-1,\alpha}(\partial\Omega)=A^{m-1,\alpha}\oplus {\mathrm{Ker}}\left(-\frac{1}{2}I+W_\Omega^t[\tilde{S}_{n,k;r},\cdot]
\right)_{|C^{\max\{1,m\}-1,\alpha}(\partial\Omega)} \,,
\]
where the direct sum is topological (cf.~(\ref{mldc_eq:vm-1a}) and Remark \ref{mldc_rem:wtnotation}). Then the map $\Xi$ from
\[
A^{m-1,\alpha}
\times\left\{
\psi \in C^{m,\alpha}(\partial\Omega):\,-\frac{1}{2}\psi+W_\Omega[\tilde{S}_{n,k;r},\psi]=0
\right\}
\]
 to the subspace
 \[
 \{
 u\in C^{m,\alpha}(\overline{\Omega}):\,\Delta u+k^2 u=0\ \text{in}\ \Omega
 \}
 \]
of $C^{m,\alpha}(\overline{\Omega})$ that takes $(\phi,\psi)$ to
\[
\Xi[\phi,\psi]\equiv v_\Omega^+[\tilde{S}_{n,k;r},\phi]
+w_\Omega^+[\tilde{S}_{n,k;r},\psi]
\]
is a linear homeomorphism.
\end{theorem}
{\bf Proof.}   Since the operator $W_\Omega^t[\tilde{S}_{n,k;r},\cdot]$ is compact in $C^{\max\{1,m\}-1,\alpha}(\partial\Omega)$
 (cf.~\textit{e.g.}, \cite[Cor.~10.1]{DoLa17}), the space 
\[
{\mathrm{Ker}}\left(
-\frac{1}{2}I+W_{\Omega}^t[\tilde{S}_{n,k;r},\cdot]
\right)_{|C^{\max\{1,m\}-1,\alpha}(\partial\Omega)}
\]
has finite dimension and is a closed finite dimensional subspace of   $V^{m-1,\alpha}(\partial\Omega)$. Then a topological supplement as $A^{m-1,\alpha}$ is well known to exist (cf.~\textit{e.g.}, Taylor and Lay \cite[Chap.~IV, \S 12,   Thm.~12.3]{TaLa80}).  

If $(\phi,\psi)$ belongs to the domain of $\Xi$, then classical  regularity results for the  acoustic  single and double layer potentials imply that
\[
 v_\Omega^+[\tilde{S}_{n,k;r},\phi]+w_\Omega^+[\tilde{S}_{n,k;r},\psi]\in C^{m,\alpha}(\overline{\Omega})
 \]
(cf.~\cite[Thms.~5.1, 5.3]{La25a} in case $m=0$, \cite[Thms.~7.1, 7.3]{DoLa17}  in case $m\geq 1$). We now show that $\Xi$ is injective. To do so, we now assume that $(\phi,\psi)$ belongs to the domain of $\Xi$ and that
\begin{equation}\label{mldc_thm:reprhelfun1}
v_\Omega^+[\tilde{S}_{n,k;r},\phi]+w_\Omega^+[\tilde{S}_{n,k;r},\psi]=0\qquad\text{in}\  \overline{\Omega}
\end{equation}
and we prove that $\phi=0=\psi$. By taking the distributional normal derivative of Definition \ref{mldc_defn:conoderdedu}, we have
\begin{equation}\label{mldc_thm:reprhelfun1a}
\frac{\partial}{\partial \nu_\Omega}v_\Omega^+[\tilde{S}_{n,k;r},\phi]+\frac{\partial}{\partial \nu_\Omega}
w_\Omega^+[\tilde{S}_{n,k;r},\psi]=0\qquad\text{on}\ \partial\Omega \,.
\end{equation}
 Since 
\[
-\frac{1}{2}\psi+W_{\Omega}[\tilde{S}_{n,k;r},\psi]=0\,,
\]
then our assumption that $k^{2}$ is not a Dirichlet eigenvalue for $-\Delta$ in $(\Omega^-)_{j}$ for each $j\in {\mathbb{N}}\setminus\{0\}$ such that $j\leq \varkappa^{-}$ and 
\cite[Thm.~9.2 (ii)]{La25b}  imply that $w_\Omega^+[\tilde{S}_{n,k;r},\psi]$ belongs to
$C^{\max\{1,m\},\alpha}(\overline{\Omega})$ and satisfies the equality $\frac{\partial}{\partial \nu_\Omega}
w_\Omega^+[\tilde{S}_{n,k;r},\psi]=0$ on $\partial\Omega$. Then the above equality  (\ref{mldc_thm:reprhelfun1a}) implies that $\frac{\partial}{\partial \nu_\Omega}v_\Omega^+[\tilde{S}_{n,k;r},\phi]=0$  on $\partial\Omega$ and the jump formula for the normal derivative of the acoustic single layer potential implies that
\[
-\frac{1}{2}\phi+W_{\Omega}^t[\tilde{S}_{n,k;r},\phi]=0\qquad\text{on}\ \partial\Omega\,,
\]
(cf.~\cite[Thm.~6.3]{La25a} in case $m=0$, \cite[Thm.~7.1]{DoLa17}  in case $m\geq 1$).  Then Corollary \ref{mldc_corol:eipedd} (i) implies that $\phi\in C^{\max\{1,m\}-1,\alpha}(\partial\Omega)$ and accordingly that
\begin{equation}\label{mldc_thm:reprhelfun1aa}
\phi\in {\mathrm{Ker}}\left(
-\frac{1}{2}I+W_{\Omega}^t[\tilde{S}_{n,k;r},\cdot]
\right)_{|C^{\max\{1,m\}-1,\alpha}(\partial\Omega)} 
 \,.
\end{equation}
Since $(\phi,\psi)$ belongs to the domain of $\Xi$, then  $\phi$ belongs to the topological supplement  $A^{m-1,\alpha}$ of the space in (\ref{mldc_thm:reprhelfun1aa}) 
and thus $\phi=0$. Then we go back to equation (\ref{mldc_thm:reprhelfun1}) and obtain
\[
w_\Omega^+[\tilde{S}_{n,k;r},\psi]=0\qquad\text{in}\  \overline{\Omega}\,.
\]
Since $-\frac{1}{2}\psi+W_{\Omega}[\tilde{S}_{n,k;r},\psi]=0$, the jump formula for the  acoustic  double layer potential implies that 
\[
\psi=w_\Omega^+[\tilde{S}_{n,k;r},\psi]_{|\partial\Omega}-w_\Omega^-[\tilde{S}_{n,k;r},\psi]_{|\partial\Omega}=0\,,
\]
(cf.~(\ref{mldc_thm:dlay1a})).
 Hence, $\Xi$ is injective. We now turn to show that $\Xi$ is surjective. Let $u\in C^{m,\alpha}(\overline{\Omega})$ be such that $\Delta u+k^2 u=0$ in $\Omega$.   
 Then we set
\[
g\equiv\frac{\partial u}{\partial\nu_\Omega}\,,
\]
we note that $g\in V^{-1,\alpha}(\partial\Omega)$ in case $m=0$, $g\in C^{m-1,\alpha}(\partial\Omega)$ in case $m\geq 1$
and we consider the  integral equation
\begin{equation}\label{mldc_thm:reprhelfun2}
-\frac{1}{2}\phi+W_{\Omega}^t[\tilde{S}_{n,k;r},\phi]=g\qquad\text{on}\ \partial\Omega\,.
\end{equation}
By classical results, $W_\Omega[\tilde{S}_{n,k;r},\cdot]$ is compact in $C^{m,\alpha}(\partial\Omega)$ in case $m\geq 1$ (cf.~\textit{e.g.}, \cite[Thm.~7.4]{DoLa17}), $W_\Omega^t[\tilde{S}_{n,k;r},\cdot]$ is compact in
$C^{m-1,\alpha}(\overline{\Omega})$ in case $m\geq 1$ (cf.~\textit{e.g.}, \cite[Cor.~10.1]{DoLa17}) and $W_\Omega^t[\tilde{S}_{n,k;r},\cdot]$ is compact in $V^{-1,\alpha}(\partial\Omega)$ in case $m=0$
(cf.~\cite[Cor.~8.5]{La25a}). Then  the Fredholm Alternative Theorem of Wendland in the duality pairing (\ref{mldc_rem:wtnotation2})
  implies that there exists 
$
 \phi\in  
 V^{m-1,\alpha}(\partial\Omega) 
 $
that satisfies equation (\ref{mldc_thm:reprhelfun2}) if and only if
\begin{equation}\label{mldc_thm:reprhelfun3}
\langle
g, \zeta \rangle=0\qquad\forall \zeta\in {\mathrm{Ker}}\left(
-\frac{1}{2}I+W_{\Omega}[\tilde{S}_{n,k;r},\cdot]
\right)_{|C^{\max\{1,m\},\alpha}(\partial\Omega)}
\,.
\end{equation}
(cf.~\textit{e.g.}, Kress~\cite[Thm.~4.17]{Kr14}). By our assumption that $k^{2}$ is not a Dirichlet eigenvalue for $-\Delta$ in $(\Omega^-)_{j}$ for each $j\in {\mathbb{N}}\setminus\{0\}$ such that $j\leq \varkappa^{-}$ and \cite[Cor.~9.9]{La25b}  (see also \cite[Thm.~7]{La12} in case $m\geq 1$), we know that
\begin{eqnarray*}
\lefteqn{
{\mathrm{Ker}}\left(
-\frac{1}{2}I+W_{\Omega}[\tilde{S}_{n,k;r},\cdot]
\right)_{|C^{\max\{1,m\},\alpha}(\partial\Omega)}
}
\\ \nonumber
&&\qquad
=\left\{
v_{|\partial\Omega}:\,v\in C^{m,\alpha}(\overline{\Omega}), \Delta v+k^2 v=0,\ \frac{\partial v}{\partial\nu_\Omega}=0
\right\}
\\ \nonumber
&&\qquad
=\left\{
v_{|\partial\Omega}:\,v\in C^{\max\{1,m\},\alpha}(\overline{\Omega}), \Delta v+k^2 v=0,\ \frac{\partial v}{\partial\nu_\Omega}=0
\right\}
\,.
\end{eqnarray*}
 Hence, equation (\ref{mldc_thm:reprhelfun2}) can have a solution $\phi\in V^{m-1,\alpha}(\partial\Omega) $ if and only if
\[
\langle
g, v_{|\partial\Omega} \rangle=0
\]
for all $v\in C^{\max\{1,m\},\alpha}(\overline{\Omega})$ such that $ \Delta v+k^2 v=0$, $\frac{\partial v}{\partial\nu_\Omega}
=0$. 
Next we note that if $m=0$, equality $\Delta u +k^2 u=0$ implies that $u\in C^{0,\alpha}(\overline{\Omega})_\Delta$ and equality $ \Delta v+k^2 v=0$
implies that $v\in C^{1,\alpha}(\overline{\Omega})_\Delta$.
 Then  the second Green Identity in distributional form of \cite[Thm.~8.1]{La25b} and equality (\ref{mldc_prop:nschext3}),  imply that 
\begin{eqnarray*}
\lefteqn{
\langle
g, v_{|\partial\Omega}\rangle=\langle\frac{\partial u}{\partial\nu_\Omega},v_{|\partial\Omega}\rangle
}
\\ \nonumber
&&\qquad
=\int_{\partial\Omega}u\frac{\partial v}{\partial\nu_\Omega}\,d\sigma
+
\langle E^\sharp[\Delta u+k^2 u],v\rangle-\int_{\Omega}u\ (\Delta v+k^2 v)\,dx=0 
\end{eqnarray*}
for all $v\in C^{\max\{1,m\},\alpha}(\overline{\Omega})$ such that $ \Delta v+k^2 v=0$, $\frac{\partial v}{\partial\nu_\Omega}
=0$.  Hence, there exists $\phi\in V^{m-1,\alpha}(\partial\Omega) $ such that equation  (\ref{mldc_thm:reprhelfun2}) is satisfied. By our assumption on $A^{m-1,\alpha}$, there exist
\[
\phi_1\in A^{m-1,\alpha}\,,\quad \phi_2\in {\mathrm{Ker}}\left(
-\frac{1}{2}I+W_{\Omega}^t[\tilde{S}_{n,k;r},\cdot]
\right)
_{|C^{\max\{1,m\}-1,\alpha}(\partial\Omega)}
\]
such that
\[
\phi=\phi_1+\phi_2\,.
\]
Since $ \phi\in  
 V^{m-1,\alpha}(\partial\Omega)$   solves equation (\ref{mldc_thm:reprhelfun2}), we have 
\[
-\frac{1}{2}\phi_1+W_{\Omega}^t[\tilde{S}_{n,k;r},\phi_1]=g\,.
\]
Next we note that
\[
w\equiv u-v_\Omega^+[\tilde{S}_{n,k;r},\phi_1]
\]
satisfies the interior Neumann problem
\[
\left\{
\begin{array}{ll}
 \Delta w +k^2 w=0 & \text{in}\ \Omega\,,
 \\
 \frac{\partial w}{\partial\nu_\Omega}=0 & \text{in}\ \partial\Omega\,.
\end{array}
\right.
\]
By our assumption that $k^{2}$ is not a Dirichlet eigenvalue for $-\Delta$ in $(\Omega^-)_{j}$ for each $j\in {\mathbb{N}}\setminus\{0\}$ such that $j\leq \varkappa^{-}$ and by \cite[Cor.~9.9 (ii)]{La25b}  (see also \cite[Thm.~7]{La12} in case $m\geq 1$), there exists a unique $\psi\in C^{m,\alpha}(\partial\Omega)$ such that $-\frac{1}{2}\psi+W_\Omega[\tilde{S}_{n,k;r},\psi]=0$ and
\[
w=w^{+}_\Omega[\tilde{S}_{n,k;r},\psi]\,.
\]
Then we have
\[
u-v^{+}_\Omega[\tilde{S}_{n,k;r},\phi_1]=w^{+}_\Omega[\tilde{S}_{n,k;r},\psi]\qquad\text{in}\ \overline{\Omega}\,.
\]
Since $(\phi_1,\psi)$ belongs to the domain of $\Xi$, we have 
  $\Xi[\phi_1,\psi]=u$. Hence, $\Xi$ is surjective.\hfill  $\Box$ 

\vspace{\baselineskip}

\begin{remark}\label{mldc_rem:reprhelfun}
 Under the assumptions of Theorem \ref{mldc_thm:reprhelfun}, if $m\geq 1$ (a variational case), then we have $ V^{m-1,\alpha}(\partial\Omega)= C^{m-1,\alpha}(\partial\Omega)$ and we can choose
\begin{eqnarray*}
\lefteqn{
A^{m-1,\alpha}= \biggl\{\phi\in C^{m-1,\alpha}(\partial\Omega):\,\int_{\partial\Omega}\phi\overline{\omega}\,d\sigma=0 
}
\\ \nonumber
&&\qquad\qquad\qquad\qquad
 \forall\omega\in {\mathrm{Ker}}\left(
-\frac{1}{2}I+W_{\Omega}^t[\tilde{S}_{n,k;r},\cdot]
\right)_{|C^{\max\{1,m\}-1,\alpha}(\partial\Omega)}  \biggr\}\,.
\end{eqnarray*}
\end{remark}
Next we introduce the following corollary, that generalizes to case $m\geq 0$ the variational case $m\geq 1$ of reference  \cite[Thm.~A1 of the Appendix]{AkLa22} with Akyel.
\begin{corollary}\label{mldc_corol:reprhelfun}
Let the same assumptions of  the representation Theorem \ref{mldc_thm:reprhelfun} hold. If we further assume that $k^2$ is not a Neumann eigenvalue for $-\Delta$ in $\Omega$, then the map from 
$V^{m-1,\alpha}(\partial\Omega)$ (cf.~(\ref{mldc_eq:vm-1a})) to the subspace
\[
\left\{
u\in C^{m,\alpha}(\overline{\Omega}):\,\Delta u+k^2 u=0\ \text{in}\ \Omega
\right\}
\]
 of $C^{m,\alpha}(\overline{\Omega})$, which takes $\phi$ to $v^{+}_\Omega[\tilde{S}_{n,k;r},\phi]$ is a linear homeomorphism. 
\end{corollary}
{\bf Proof.} Under the assumptions  that $k^2$ is not a Neumann eigenvalue for $-\Delta$ in $\Omega$ and that $k^{2}$ is not a Dirichlet eigenvalue for   for $-\Delta$ in $(\Omega^-)_{j}$ for each $j\in {\mathbb{N}}\setminus\{0\}$ such that $j\leq \varkappa^{-}$,  Remark \ref{mldc_rem:dimker-12i+wk2t} implies that 
\[
{\mathrm{Ker}}\left(
-\frac{1}{2}I+W_{\Omega}^t[\tilde{S}_{n,k;r},\cdot]
\right)_{|C^{\max\{1,m\}-1,\alpha}(\partial\Omega)}
=\{0\}\,.
\]
Under the assumptions  that $k^2$ is not a Neumann eigenvalue for $-\Delta$ in $\Omega$ and that $k^{2}$ is not a Dirichlet eigenvalue for   for $-\Delta$ in $(\Omega^-)_{j}$ for each $j\in {\mathbb{N}}\setminus\{0\}$ such that $j\leq \varkappa^{-}$, \cite[Cor.~9.9]{La25b} implies that
\[
{\mathrm{Ker}}\left(
-\frac{1}{2}I+W_{\Omega}[\tilde{S}_{n,k;r},\cdot]
\right)_{|C^{m,\alpha}(\partial\Omega)}
=\{0\}\,.
\]
 and accordingly the representation Theorem \ref{mldc_thm:reprhelfun}
with $A^{m-1,\alpha}= V^{m-1,\alpha}(\partial\Omega)$  implies the validity of the statement. \hfill  $\Box$ 

\vspace{\baselineskip}

\section{Appendix: a   consequence of a Holmgren u\-niqueness theorem}
As pointed out by Lupo and Micheletti \cite[Proof of Theorem 2]{LuMi93}, one can deduce
the uniqueness Theorem \ref{mldc_thm:holmgr} below by a Theorem of Holmgren.  Solely for the convenience of the reader we include a proof. If $\Omega$ is an open subset of ${\mathbb{R}}^n$, then 
$W^{1,2}(\Omega)$ denotes the set of (equivalence classes of) functions of class $L^2(\Omega)$ such that the first order distributional partial derivatives of order one are in 
$L^2(\Omega)$ with its usual norm
\[
\|u\|_{W^{1,2}(\Omega)}\equiv\sum_{|\eta|\leq 1}\|D^\eta u\|_{L^2(\Omega)}
\qquad\forall u\in W^{1,2}(\Omega)\,,
\]
and $W^{1,2}_0(\Omega)$ denotes the closure of the space of test functions $C^\infty_c(\Omega)$ in $W^{1,2}(\Omega)$.
\begin{theorem}[of Holmgren]\label{mldc_thm:holmgr}
 Let $\Omega $ be a bounded open Lipschitz subset of ${\mathbb{R}}^n$.
 Let $u\in W^{1,2}_0(\Omega) $ be a Neumann eigenfunction of $-\Delta$ in $\Omega$, \textit{i.e.}, there exists $\lambda\in {\mathbb{C}}$ such that
\begin{equation}\label{mldc_thm:holmgr1}
\int_\Omega -\nabla u\nabla v+\lambda uv\,dx=0\qquad\forall v\in W^{1,2}(\Omega)\,.
\end{equation}
Then $u=0$ (almost everywhere) in $\Omega$. 
\end{theorem}
{\bf Proof.} The function $U$ defined by
\[
U(x)\equiv
\left\{
\begin{array}{ll}
 u(x) & \text{if}\ x\in\Omega
 \\
0  & \text{if}\ x\in {\mathbb{R}}^n\setminus\Omega
\end{array}
\right.
\]
solves $\Delta U+\lambda U=0$ in the sense of distributions in  ${\mathbb{R}}^n$.  
 Indeed, $U$ is locally integrable in ${\mathbb{R}}^n$  and thus
\begin{eqnarray*}
\lefteqn{
\langle\Delta U+\lambda U,\varphi\rangle=\int_{{\mathbb{R}}^n }U(x)\Delta\varphi(x)+\lambda U(x)\varphi(x)\,dx
}
\\ \nonumber
&&\qquad\qquad\qquad
=\int_{\Omega }U(x)\Delta\varphi(x)+\lambda U(x)\varphi(x)\,dx
\\ \nonumber
&&\qquad\qquad\qquad
=\int_{\Omega }u(x)\Delta\varphi(x)+\lambda u(x)\varphi(x)\,dx
\end{eqnarray*}
for all  $\varphi\in {\mathcal{D}}( {\mathbb{R}}^n )$ 
and the membership of $u\in W^{1,2}_0(\Omega) $ implies that
\[
\int_{\Omega }u(x)\Delta\varphi(x)\,dx=-\int_{\Omega }\nabla u(x)\nabla\varphi(x)\,dx
\]
for all  $\varphi\in {\mathcal{D}}({\mathbb{R}}^n )$ and thus
\begin{eqnarray*}
\lefteqn{
\langle\Delta U+\lambda U,\varphi \rangle=\int_{\Omega }u(x)\Delta\varphi(x)+\lambda u(x)\varphi(x)\,dx
}
\\ \nonumber
&&\qquad\qquad\qquad
=  \int_{\Omega }- \nabla u(x)\nabla\varphi(x)+\lambda u(x)\varphi(x)\,dx =0
\end{eqnarray*}
for all  $\varphi\in {\mathcal{D}}( {\mathbb{R}}^n )$. Since $U$ solves an elliptic equation with constant coefficients in the sense of the distributions  ${\mathbb{R}}^n$, then $U$ is (almost everywhere) equal to a real analytic function in ${\mathbb{R}}^n$. Since $U$ equals zero on the nonempty open set ${\mathbb{R}}^n\setminus\over-line{\Omega}$, we have $U=0$ (almost everywhere) in ${\mathbb{R}}^n$ and thus $u=0$ (almost everywhere) in $\Omega$.\hfill  $\Box$ 

\vspace{\baselineskip}

We can now deduce the following (nonvariational) consequence of the Holmgren Theorem \ref{mldc_thm:holmgr}.
 \begin{corollary}\label{mldc_corol:holmgr}
 Let  $\alpha\in]0,1[$. Let $\Omega$ be a bounded open subset of ${\mathbb{R}}^{n}$ of class $C^{1,\alpha}$. Let $\lambda\in {\mathbb{C}}$. If $u\in C^{0,\alpha}(\overline{\Omega})$ satisfies the Neumann problem (\ref{mldc_thm:nexintneuhe3}) and $u_{|\partial\Omega}=0$, then $u=0$ in $\overline{\Omega}$.
\end{corollary}
{\bf Proof.} Since $\Delta u=-\lambda u\in C^{0,\alpha}(\overline{\Omega})\subseteq C^{-1,\alpha}(\overline{\Omega})$, we have $u\in C^{0,\alpha}(\overline{\Omega})_\Delta$. Then condition    $u_{|\partial\Omega}=0$ and the   regularization theorem of  \cite[Thm.~1.2 (iii)]{La25} implies that $u\in C^{1,\alpha}(\overline{\Omega})$. Moreover, Remark \ref{mldc_rem:helpo} (i)
 and \cite[(2.21)]{La25b} imply that 
   $u\in C^{1,\alpha}(\overline{\Omega})_\Delta  \subseteq C^2(\Omega)$. Then the first Green Identity
   implies that $u$ solves the weak form (\ref{mldc_thm:holmgr1}) for all $v\in C^\infty(\overline{\Omega})$
   (cf.~\textit{e.g.}, \cite[Thm.~4.2]{DaLaMu21}).  Since   $C^\infty(\overline{\Omega})$ is dense in $W^{1,2}(\Omega)$ (cf.~\textit{e.g.},  Brezis \cite[Cor.~9.8]{Br11}),  then the H\"{o}lder inequality implies that
     $u$ solves the weak form (\ref{mldc_thm:holmgr1}) for all $v\in W^{1,2}(\Omega)$. Since  $u_{|\partial\Omega}=0$ and $u\in C^{1,\alpha}(\overline{\Omega})\subseteq W^{1,2}(\Omega)$, we have   $u\in W^{1,2}_0(\Omega)$ (cf.~\textit{e.g.}, Brezis \cite[Thm.~9.17]{Br11}). Hence,  the Holmgren   Theorem \ref{mldc_thm:holmgr} implies that $u=0$ in $\overline{\Omega}$.\hfill  $\Box$ 

\vspace{\baselineskip}

 \noindent
{\bf Statements and Declarations}\\

 \noindent
{\bf Competing interests:} This paper does not have any  conflict of interest or competing interest.

 \noindent
{\bf Acknowledgement.}   The author  acknowledges  the support of   GNAMPA-INdAM,   of the Project funded by the European Union – Next Generation EU under the National Recovery and Resilience Plan (NRRP), Mission 4 Component 2 Investment 1.1 - Call for tender PRIN 2022 No. 104 of February, 2 2022 of Italian Ministry of University and Research; Project code: 2022SENJZ3 (subject area: PE - Physical Sciences and Engineering) ``Perturbation problems and asymptotics for elliptic differential equations: variational and potential theoretic methods'' and the support from the EU through the H2020-MSCA-RISE-2020 project EffectFact, Grant agreement ID: 101008140.

\end{document}